\documentclass[12pt]{article}
\usepackage[all]{xy}
\usepackage{amsmath}
\usepackage{amsfonts}
\usepackage{amssymb}
\usepackage{amscd}
\usepackage{amsthm}
\usepackage{latexsym}
\usepackage{amsbsy}
\newtheorem{lem}{Lemma}[section]
\newtheorem{prop}{Proposition}[section]
\newtheorem{theorem}{Theorem}[section]

\newtheorem{cor}{Corollary}[section]
\newtheorem{definition}{Definition}[section]
\newtheorem{example}{Example}[section]

\newtheorem{remark}{Remark}

\begin{document}
\begin{titlepage}
\title{Very Basic Noncommutative Geometry}
\author{ \\   \\  \\ \\ \\ \\ Masoud Khalkhali\\Mathematics Department, University of Western Ontario\\
London ON, Canada}
\date{ }
\maketitle

\end{titlepage}
\tableofcontents

 \section{Introduction}

One of the major advances of science in the 20th century  was the discovery of a
mathematical formulation of  quantum mechanics by Heisenberg in
1925 \cite{hei}. From a mathematical point of view, transition from  classical
mechanics to quantum mechanics amounts to, among other things, passing from the {\it commutative
algebra} of {\it classical observables} to the {\it noncommutative algebra} of
{\it quantum mechanical observables}. Recall that in classical mechanics
an observable (e.g. energy, position, momentum, etc.) is a
function on a manifold called the phase space of the system. Immediately
after Heisenberg's work,
ensuing papers by Dirac \cite{dir} and Born-Heisenberg-Jordan \cite{bhj},
made it clear that a quantum mechanical observable is a
(selfadjoint) operator on a Hilbert space called the state
space of the system. Thus the commutative algebra
of functions on a space is replaced by the noncommutative algebra of
operators on a Hilbert space.

A little more than fifty years after these developments,  Alain Connes realized that a
similar procedure can in fact be applied to areas of mathematics where the classical
notions of  space (e.g. measure space, locally compact space, or a smooth space)
looses its applicability and pertinence and can be replaced by a new idea of space, represented
by a noncommutative algebra.

Connes' theory, which is generally known as  {\it noncommutative
geometry}, is a rapidly growing new area of mathematics that interacts
with and contributes to many disciplines in  mathematics and physics.
For a recent survey see Connes' article \cite{ac00}.
 Examples of such
interactions and contributions include:
 theory of
operator algebras, index theory of elliptic operators, algebraic and
differential topology, number theory, standard model of elementary
particles, quantum Hall effect, renormalization in quantum field theory, and string theory. (For a description of these relations in more details see the report below.)
 To understand the basic ideas of
noncommutative geometry one should perhaps first come to grips with the
idea of a {\it noncommutative space}.

The inadequacy of the classical notions of space manifests itself for example when one deals with
highly singular \textquoteleft \textquoteleft bad
quotients";  spaces such as the quotient of a nice
space by the ergodic action of a group or the space of leaves of a
foliation. In all these examples the quotient space is typically ill behaved even
as a topological  space. For example it may fail to be even
Hausdorff, or have enough open sets, let alone being a reasonably smooth
space. The unitary dual of a noncompact (Lie) group, except when the group is abelian or almost
abelian,  is another example of an ill behaved space.

One of Connes' key
observations is that in all these situations one can attach a
noncommutative algebra, through a {\it noncommutative quotient construction}, that captures most of
the information. Examples of this noncommutative quotient construction include crossed product by action
of a group, or by action of a groupoid. In general the noncommutative quotient is the groupoid
algebra of a topological groupoid.

Noncommutative geometry has as its limiting case the classical
geometry, but geometry  expressed in algebraic terms.  Thus to
understand its relation with classical geometry one should first
understand one of the most important ideas of mathematics which
can be expressed as a {\it duality} between commutative algebra
and geometry. This is by no  means a new observation or a new
trend. To the contrary, this duality has always existed and been
utilized in mathematics and its applications. The earliest example
is perhaps the use of numbers in counting!
 It is, however, the case that throughout the  history each new generation of mathematicians find new ways of
formulating this principle and at the same time broaden its scope.
Just to mention a few highlights of this rich history we mention
Descartes (analytic geometry), Hilbert (affine varieties and commutative
algebras), Gelfand-Naimark (locally compact spaces and commutative
$C^*$-algebras), Grothendieck (Schemes and topos theory), and Connes (noncommutative geometry).

A key idea here is the well-known relation between a  space and
the commutative algebra of functions on that space. More precisely
there is a duality between certain categories of geometric spaces
and categories of algebras representing those spaces.
Noncommutative geometry builds on, and vastly extends, this
fundamental duality between geometry and commutative algebras.

 For example, by a
celebrated theorem of Gelfand and Naimark \cite{gelnai} one knows that the category of locally compact
Hausdorff
spaces is equivalent to the dual of the category of commutative  $C^*$-algebras. Thus one can think of
not necessarily commutative $C^*$-algebras  as the dual of a  category of {\it noncommutative locally
compact spaces}. What makes this a successful proposal is first of all a rich supply of
examples and secondly the possibility of extending many of the topological and geometric
invariants to this new class of spaces.  Let us briefly recall a few
other
 examples from  a long list of results in mathematics that put in duality certain categories of geometric
objects with a corresponding  category of algebraic objects.

To wit, Hilbert's Nullstellensatz states that the category of
algebraic varieties over an algebraically closed  field is
equivalent to the dual  of the category of finitely generated
commutative algebras without nilpotent elements (so called reduced
algebras). This is a perfect analogue of the Gelfand-Naimark
theorem in the world of commutative algebras.

Similarly, the Serre-Swan theorem states that
the category of vector bundles
over a compact Hausdorff space (resp. affine algebraic variety) X is equivalent to the category of
finitely generated projective modules over the algebra of continuous functions
(resp. regular functions) on X.

Thus a pervasive idea in noncommutative geometry is to treat
(certain classes) of noncommutative algebras as noncommutative
spaces and  try to extend tools of geometry, topology, and
analysis to this new setting.  It should be emphasized, however,
that, as a rule, this extension is never straightforward and
always involve surprises and new phenomena. For example the theory
of the flow of weights and the corresponding modular automorphism
group in von Neumann algebras has no counterpart in classical
measure theory, though the theory of von Neumann algebras is
generally regarded as noncommutative measure theory.  Similarly
the extension of de Rham homology for manifolds to cyclic
cohomology for noncommutative algebras was not straightforward and
needed some highly nontrivial considerations.

Of all the topological invariants for spaces, topological $K$-theory has
the most  straightforward extension to the noncommutative realm. Recall
that topological $K$-theory classifies vector bundles on a topological space. Using the above mentioned
Serre-Swan theorem, it is natural to define, for a not necessarily
commutative ring $A$,  $K_0(A)$ as the group defined by the semi-group of isomorphism classes of finite
projective $A$-modules. The definition of $K_1(A)$ follows the same
pattern as in the commutative case, provided $A$ is a Banach algebra  and the main theorem of
topological $K$-theory,  the Bott periodicity theorem, extends
to all Banach algebras.

The situation is much less clear for $K$-homology, the theory dual
to $K$-theory.  By the work of Atiyah, Brown-Douglas-Fillmore, and
Kasparov, one can say, roughly speaking, that $K$-homology cycles
on a space $X$ are represented by abstract elliptic operators on
$X$ and while $K$-theory classifies vector bundles on $X$,
$K$-homology classifies the abstract elliptic operators on $X$.
The pairing between  $K$-theory and $K$-homology takes the form
$<[D], [E]>=$ the Fredholm index of the elliptic operator $D$ with
coefficients in the vector bundle $E$. Now one good thing about
this way of formulating $K$-homology is that it almost immediately
extends to noncommutative $C^*$-algebras. The two theories are
unified in a single theory called $KK$-theory due to G. Kasparov.

Cyclic cohomology was discovered by Connes in 1981 \cite{ac81,
ac85} as the right noncommutative analogue of de Rham homology of
currents and as a target space for noncommutative Chern character
maps from both $K$-theory and $K$-homology. One of the main
motivations of Connes seems to be transverse index theory on
foliated spaces. Cyclic cohomology can be used to identify the
$K$-theoretic index of transversally  elliptic operators which lie
in the $K$-theory of the noncommutative algebra of the foliation.
The formalism of cyclic cohomology and Chern-Connes character maps
form an indispensable part of noncommutative geometry. In a
different direction, cyclic homology also appeared in the 1983
work of Tsygan \cite{tsy} and was used, independently, also by
Loday and Quillen \cite{lq} in their study of the Lie algebra
homology of the Lie algebra of stable matrices over an associative
algebra.
   We won't pursue this aspect of
 cyclic homology in these notes.

A very interesting recent develpoment in cyclic cohomology theory is the
{\it Hopf-cyclic cohomology} of Hopf algebras and Hopf module
(co)algebras in general. Motivated by the original work of Connes and Moscovici \cite{cm1, cm2} this
theory is now  extended and elaborated on by several  authors \cite{ak1,
ak2, hkrs1, hkrs2, kr1, kr2, kr3, kr4}. There are also very interesting
relations between cocycles for Hopf-cyclic cohomology theory of the Connes-Moscovici
Hopf algebra $\mathcal{H}_1$  and operations on spaces
of modular forms and modular Hecke algebras \cite{cm3, cm4}, and spaces
of $\mathbb{Q}$-lattices \cite{cma}. We will say nothing about these
develpments in these notes. Neither we shall discuss the approach of  Cuntz and Quillen
 to cyclic cohomology theory and their cellebrated proof of excision property for periodic
 (bivariant) cyclic cohomology \cite{cu,  cst, cq1, cq2, cq3}.   

The following \textquoteleft \textquoteleft dictionary"
illustrates   noncommutative analogues of some of the classical
theories and concepts originally conceived for spaces. In these
notes we deal only with a few items of this dictionary. For a much
fuller account and explanations, as well as applications of noncommutative geometry,
the reader should consult Connes'
beautiful book  \cite{acb}.

$$
\begin{tabular}{|l|l|}
\hline
{\bf commutative} & {\bf noncommutative}\\
\hline
measure space & von Neumann algebra \\
locally compact  space & $C^\ast$ - algebra \\
vector bundle & finite projective module\\
complex variable & operator on a Hilbert space\\
real variable & sefadjoint operator\\
infinitesimal & compact operator\\
range of a function & spectrum of an operator\\
$K$-theory & $K$-theory\\
vector field & derivation\\
integral & trace\\
closed de Rham current & cyclic cocycle\\
de Rham complex & Hochschild homology \\
de Rham cohomology & cyclic homolgy \\
Chern character & Chern-Connes character \\
Chern-Weil thoery & noncommutative Chern-Weil thoery\\
elliptic operator & $K$-cycle\\
spin Riemannian manifold & spectral triple\\
index theorem & local index formula\\
group, Lie algebra & Hopf algebra, quantum group\\
symmetry & action of Hopf algebra\\
\hline
\end{tabular}
$$

Noncommutative geometry is already a vast subject. These notes are
just  meant to be an introduction to a few aspects of this
fascinating  enterprize. To get a much better sense of the beauty
and depth of the subject the reader should consult Connes'
magnificent book \cite{acb} or his recent survey \cite{ac00} and
references therein.  Meanwhile, to give a
sense of the state of the subject at the present time, its
relation with other fields of mathematics,
 and its most pressing issues, we reproduce here
  part of the text of  the final report prepared by the organizers of a conference on noncommutative
  geometry in 2003
  \footnote{BIRS Workshop on Noncommutative Geometry, Banff International Research Station, Banff, Alberta,
  Canada, April 2003, Organized by Alain Connes, Joachim Cuntz, George Elliott, Masoud Khalkhali, and
  Boris Tsygan. Full report available at: www.pims.math.ca/birs.}:
   \\

{\it \textquoteleft \textquoteleft
{\bf 1. The Baum-Connes conjecture}\\

  This conjecture, in its simplest form, is
 formulated for any locally compact topological group. There are more
 general Baum-Connes conjectures with coefficients
 for groups acting on C*-algebras,
  for groupoid C*-algebras, etc., that for the sake of brevity we don't
 consider here. In a nutshell the Baum-Connes conjecture
 predicts that the K-theory of the group C*-algebra of a  given
 topological group is isomorphic, via an explicit map called the Baum-Connes map, to an appropriately defined
 K-homology of the classifying space of the group. In other words
 invariants of groups defined through noncommutative geometric tools
 coincide with invariants defined through classical algebraic
 topology
 tools. The Novikov conjecture on the homotopy invariance of higher signatures of non-simply
 connected manifolds is a consequence of the Baum-Connes conjecture (the relevant
 group here is the fundamental group of the manifold). Major advances were made in this problem in
 the past seven years by
  Higson-Kasparov,   Lafforgue, Nest-Echterhoff-Chabert,  Yu,
  Puschnigg and others. \\

{\bf 2. Cyclic cohomology and KK-theory}\\

 A major discovery  made by
     Alain Connes in 1981, and independently  by Boris Tsygan in 1983, was the discovery of cyclic cohomology as the
     right noncommutative analogue of de Rham homology and a natural target for a Chern character map from K-theory
     and K-homology.  Coupled with
     K-theory, K-homology and KK-theory, the formalism of cyclic
     cohomology fully extends many aspects of classical differential
     topology like Chern-Weil theory to noncommutative spaces. It is an indispensable tool in noncommutative geometry.
     In recent years Joachim Cuntz and Dan Quillen have formulated an
     alternative
     powerful new approach to cyclic homology theories which brings with it many new insights as well as a successful
     resolution of an old open problem in
     this area, namely establishing the excision property of periodic cyclic cohomology.

 For applications of noncommutative geometry to problems of index theory, e.g. index theory
on foliated spaces, it is necessary to extend the formalism of
cyclic cohomology to a bivariant cyclic theory for topological
algebras and to extend Connes's Chern character to a fully
bivariant setting. The most general approach to this problem is
due to Joachim Cuntz.  In fact the approach of Cuntz made it possible  to
extend the domain (and definition) of KK-theory to very general
categories of topological algebras (rather than just C*-algebras).
The fruitfulness of this idea  manifests itself  in the V. Lafforgue's proof
of the Baum-Connes conjecture for groups with property T, where the
extension of KK functor to Banach algebras plays an important role.

 A new trend in cyclic cohomology theory is the study of the cyclic cohomology of
 Hopf algebras and quantum groups. Many  noncommutative spaces, such as quantum spheres and quantum homogeneous spaces, admit a quantum group of
 symmetries. A remarkable discovery of Connes and Moscovici in the past
 few years is the fact that diverse structures, such as the space of leaves
 of a (codimension one) foliation  or the space of modular forms, have  a unified quantum symmetry.
   In their study of transversally
 elliptic operators on foliated manifolds  Connes and Moscovici came up with a new noncommutative and
 non-cocommutative Hopf algebra denoted by $\mathcal{H}_n$ (the Connes-Moscovici Hopf algebra). $\mathcal{H}_n$ acts on the transverse foliation algebra of
 codimension $n$ foliations and thus appears as the quantized symmetries of a foliation. They noticed that if one extends the noncommutative Chern-Weil theory
 of Connes from group and Lie algebra  actions to actions of Hopf algebras, then the characteristic
 classes defined via the local index formula are in the image
 of this  new characteristic map. This extension of Chern-Weil theory involved the introduction of
 cyclic cohomology for Hopf algebras.  \\

{\bf 3. Index theory and noncommutative geometry}\\

 The index theorem of Atiyah and Singer and its various generalizations
 and ramifications
 are at the core of noncommutative geometry and its applications. A
 modern abstract index theorem in the noncommutative setting is the
 local index formula of Connes and Moscovici. A key ingredient of
 such an abstract index formula is the idea of an spectral triple due to
 Connes. Broadly speaking, and neglecting the parity, a spectral triple $(A, H, D)$ consists of an
 algebra $A$ acting by
 bounded operators on the Hilbert space $H$ and a self-adjoint operator $D$ on $H$. This data must
 satisfy certain
 regularity properties which constitute an abstraction of basic elliptic estimates for elliptic PDE's acting
 on sections
 of vector bundles on compact manifolds.  The local index formula replaces the old non-local
 Chern-Connes cocycle by a  new Chern character form $Ch (A, H, D)$ of the given spectral triple
  in the cyclic complex of the algebra $A$. It is a local formula in the sense
 that the cochain $Ch (A, H, D)$ depends, in the classical case,  only on the germ of the heat kernel of $D$ along the diagonal  and in
 particular is independent of smooth perturbations. This makes the formula extremely attractive
 for practical calculations. The challenge now is to apply this formula to diverse situations beyond
 the cases considered so far, namely transversally elliptic operators on foliations
 (Connes and Moscovici) and the Dirac operator on quantum $SU_2$
 (Connes). \\

{\bf 4. Noncommutative geometry and number theory}\\

Current applications and connections of noncommutative geometry to
number theory can be divided into four categories. (1) The work of
Bost and Connes, where they construct a noncommutative dynamical
system $(B, \sigma_t)$ with partition function the Riemann zeta
function $\zeta (\beta )$, where $\beta$ is the inverse
temperature. They show that at the pole $\beta =1$ there is an
spontaneous symmetry breaking. The symmetry group of this system
is the group of id\'{e}les which is isomorphic to the Galois group
$Gal (Q^{ab}/Q)$. This gives a natural interpretation of the zeta
function as the partition function of a quantum statistical
mechanical system. In particular the class field theory
isomorphism appears very naturally in this context.  This approach
has been extended to the Dedekind zeta function of an arbitrary
number field by Cohen, Harari-Leichtnam, and
Arledge-Raeburn-Laca. All these results concern  abelian
extensions of number fields and their generalization to
non-abelian extensions is still lacking. (2) The work of Connes on
the Riemann hypothesis. It starts by producing a conjectural trace
formula which refines  the Arthur-Selberg trace formula. The main
result of this theory states that this trace formula is valid if
and only if the Riemann hypothesis is satisfied by all
$L$-functions with Gr\"{o}ssencharakter on the given number
field $k$. (3) The work of Connes and Moscovici on quantum
symmetries of the modular Hecke algebras $\mathcal{A}(\Gamma)$
where they show that this algebra admits a natural action of the
transverse Hopf algebra $\mathcal{H}_1$. Here $\Gamma$ is a
congruence subgroup of $SL (2, Z)$ and the algebra
$\mathcal{A}(\Gamma)$ is the crossed product of the algebra of
modular forms of level $\Gamma $ by the action of the Hecke
operators. The action of the generators $X, Y$ and $\delta_n$ of
$\mathcal{H}_1$ corresponds to the Ramanujan operator, to the weight or
number operator, and to the action of certain group cocycles on
$GL^+ (2, Q)$, respectively.  What is very surprising is that the
same  Hopf algebra  $\mathcal{H}_1$ also acts naturally on the
(noncommutative) transverse space of codimension one foliations.
(4) Relations with arithmetic algebraic geometry and Arakelov
theory. This is currently being pursued by Consani, Deninger,
Manin, Marcolli and others. \\

{\bf 5. Deformation quantization and quantum geometry }\\

The noncommutative algebras that appear in noncommutative geometry
usually are obtained either as the  result of a process called
noncommutative quotient construction or by deformation
quantization of some algebra of functions on a classical space.
These two constructions are not mutually exclusive. The starting
point of deformation quantization is an algebra of functions on a
Poisson manifold where the Poisson structure gives the
infinitesimal direction of quantization. The existence of
deformation quantizations for all Poisson manifolds was finally
settled by M. Kontsevich in 1997 after a series of partial results
for symplectic manifolds. The algebra of pseudodifferential
operators on a manifold is a deformation quantization of the
algebra of classical symbols on the cosphere bundle of the
manifold. This simple observation is the beginning of an approach
to the proof of the index theorem, and its many generalizations  by
Elliott-Natsume-Nest and Nest-Tsygan, using cyclic cohomology theory. The same can be said
about Connes's groupoid approach to index theorems. In a different direction, quantum geometry
also consists of the study of noncommutative metric spaces and noncommutative complex structures." \\}

Let us now briefly describe the contents of these notes. In Section 2 we describe some
of the fundamental algebra-geometry correspondences at work in mathematics. The most basic ones
for noncommutative geometry are the  Gelfand-Naimark  and the
Serre-Swan theorems. In Section 3 we describe the noncommutative quotient construction and give several
examples. This is one of the most universal methods of constructing noncommutative spaces directly
related to classical geometric examples.  Section 4 is devoted to cyclic
cohomology and its various definitions. In Section 5 we define the Chern-Connes character map, or
the noncommutative Chern character map,  from
$K$-theory to cyclic cohomology. In an effort to make these notes as self contained as possible,
we have added
three appendices covering very basic material on $C^*$-algebras, projective modules, and category theory
language.

These notes are partly based on series of lectures I gave at the Fields
Institute in Toronto, Canada,  in Fall 2002 and at the Institute for
Advanced Studies in Physics and Mathematics (IPM),
Tehran, Iran,  in Spring 2004. I also used  part of these notes in my lectures
at the {\it second annual spring institute and workshop on
noncommutative geometry} in Spring 2004,
Vanderbilt University, USA.  It is a great pleasure to thank the organizers of this event,
Alain Connes (director), to whom I owe much more than I can adequately express, Dietmar
Bisch, Bruce Hughes, Gennady Kasparov, and  Guoliang Yu.  I would also like to thank
Reza Khosrovshahi the director of the mathematics division of IPM in Tehran whose encouragement and
support  was instrumental
in bringing these notes to existence.

\section{ Some examples of geometry-algebra correspondence}
We give several examples of geometry-commutative algebra
correspondences. They all put into correspondence, or duality,
certain categories of geometric objects with a category of
algebraic objects. Presumably, the more one knows about these
relations the better one is prepared to pursue  noncommutative
geometry.

\subsection{Locally compact  spaces and commutative
$C^*$-algebras}

In functional analysis the celebrated {\it Gelfand-Naimark Theorem} \cite{gelnai} states that the
category of locally compact  Hausdorff spaces is anti-equivalent to the category
of commutative $C^*$-algebras:
$$\{\text{\bf locally compact Hausdorff spaces}\} \simeq \{\text{\bf  commutative $C^*$-algebras} \}^{op}.$$

Let $\mathcal{S}$ be the category whose objects are locally compact Hausdorff spaces and whose
morphisms are continuous and {\it proper} maps. (Recall that a map $f:X\to Y$ is called proper if
for any compact $K\subset Y,\; f^{-1}(K)$ is compact; of course, if $X$ is compact and $f$
is continuous, then
$f$ is proper).

Let $\mathcal{C}$ be the category whose objects are commutative
$C^\ast$-algebras and whose morphisms are {\it proper}
$\ast$-homomorphisms. (A $\ast$ -homomorphism $f: A \rightarrow B$
is called proper if for any approximate identity $(e_i)$ in $A$,
$f(e_i)$ is an approximate identity in $B$. See Appendix A for
definitions.)

Define two contravariant functors
$$C_0 :\mathcal{C} \to \mathcal{S}, \quad
\Omega : \mathcal{S}\to \mathcal{C},$$
as follows. For a  locally compact Hausdorff space $X$, let $C_0(X)$ denote  the algebra of
complex valued continuous functions on $X$ that ``vanish at $\infty$".
This means for any $\epsilon >0$ there is a compact subset $K\subset X$
such that $|f(x)|< \epsilon$ for $x\notin K$:
$$ C_0(X)= \{f: X \rightarrow \mathbb{C}, \text{ f is continuous and $f(\infty)=0$}\}.$$
Under pointwise addition and scalar multiplication $C_0(X)$ is
obviously an algebra over the field of complex numbers
$\mathbb{C}$. Endowed with the sup-norm
$$\| f \| = \| f \|_\infty = \sup \{ |f(x) | ; \; x\in X \}, $$
and $\ast$-operation
$$f\mapsto f^\ast,\; f^\ast (x) = \bar{f}(x),$$
one  checks that
$C_0(X)$ is a commutative $C^\ast$-algebra.
If $f:X\to Y$ is a continuous and proper map, let
$$ C_0(f)=f^*: C_0(Y) \longrightarrow C_0(X), \quad f^*(g)=g\circ f,$$
 be the pullback of $f$.   It is a proper  $\ast$-
homomorphism of $C^\ast$-algebras. We have thus defined the functor $C_0$.

To define $\Omega$, called the {\it functor of points} or the {\it spectrum functor}, let $A$ be a commutative
$C^\ast$-algebra. Let
$$\Omega(A) = \text{set of  characters of $A$} =Hom_{\mathcal{C}}(A, \mathbb{C}),$$
 where a {\it character} is simply a nonzero algebra map $A \rightarrow
 \mathbb{C}$. (it turns out that they are also $*$-morphisms).
$\Omega (A)$  is a locally compact Hausdorff space under the topology of pointwise
convergence.
 Given    a proper morphism of $C^*$-algebras $f:A\to B$, let
 $$\Omega (f) :\Omega (B) \to \Omega (A), \quad \Omega (f)=f^*,$$
 where $f^\ast (\varphi )=\varphi \circ f$. It can be shown that  $\Omega (f)$ is a proper and
 continuous map.

To show that $C_0$ and $\Omega$ are equivalences of categories, quasi-inverse to each
other, one shows that   for any locally compact Hausdorff space $X$ and any
commutative $C^\ast$-algebra $A$, there are natural isomorphisms
$$X\overset{\sim}{\longrightarrow}\Omega (C_0 (X)), \quad \quad x\mapsto e_x,$$
$$A\overset{\sim}{\longrightarrow} C_0 (\Omega (A)), \quad \quad a\mapsto \hat{a}.$$
Here  $e_x$  is the {\it evaluation at $x$} map defined by $e_x(f)=f(x)$, and
$a\mapsto \hat{a}$ is the celebrated {\it Gelfand transform}
defined by $\hat{a}(\varphi)=\varphi (a).$
The first isomorphism is much easier to establish and does not require the theory
of Banach algebras. The second isomorphism is what is proved by Gelfand
and Naimark in 1943 \cite{gelnai} using Gelfand's theory of commutative
Banach algebras. We sketch a proof of this result is Appendix A.

Under the Gelfand-Naimark correspondence compact Hausdorff spaces
correspond to unital $C^*$-algebras. We therefore have
a duality, or equivalence of categories
$$\{ \mbox{\bf compact Hausdorff spaces} \} \simeq\{
\mbox{\bf commutative unital
$C^\ast$-algebras} \}^{op}.$$

Based on Gelfand-Naimark theorem, we can think of  the dual of the
category  of not necessarily commutative $C^*$-algebras as  the
category of {\it noncommutative locally compact Hausdorff spaces}.
Various operations and concepts for spaces can be paraphrased in
terms of  algebras of functions on  spaces and can then  be
immediately generalized to noncommutative spaces. This is the easy
part of noncommutative geometry! Here is a dictionary suggested by
the Gelfand-Naimark theorem:

$$
\begin{matrix}
\mbox{\bf space } & \mbox{\bf algebra }\\
\mbox{compact} & \mbox{unital}\\
\mbox{1-point compactification} & \mbox{unitization}\\
\mbox{Stone-Cech compactification} & \mbox{multiplier algebra}\\
\mbox{closed subspace; inclusion } & \mbox{closed ideal; quotient algebra}\\
\mbox{surjection} & \mbox{injection}\\
\mbox{injection} & \mbox{surjection}\\
\mbox{homoemorphism} & \mbox{automorphism}\\
\mbox{Borel measure} & \mbox{ positive functional}\\
\mbox{probability measure} & \mbox{state}\\
\mbox{disjoint union}& \mbox{direct sum}\\
\mbox{cartesian product} & \mbox{minimal tensor product}\\
\end{matrix}
$$

\subsection{ Vector bundles and finite projective modules}
{\it Swan's Theorem} \cite{swa} states that the category of complex vector bundles
on a compact Hausdorff space $X$ is equivalent to the category of
finite (i.e. finitely generated)  projective modules over the algebra $C(X)$ of
continuous functions on $X$:
$$\{\text{\bf vector bundles on X}\} \simeq \{ \text{\bf finite projective C(X)-modules}\}.$$
 There are similar results for real and quaternionic vector bundles \cite{swa}. This result was motivated
 and in fact is the topological counterpart of a  an
analogous result, due to Serre, which characterizes algebraic vector bundles over an affine
algebraic variety   as finite  projective modules over
the coordinate ring of the variety. Swan's theorem sometimes is called the Serre-Swan theorem.

Recall that a right module $P$ over a unital  algebra $A$ is called {\it projective}
if there exists a right
$A$-module $Q$ such that
$$P\oplus Q\simeq \bigoplus_I A,$$
 is a free module. Equivalently, $P$ is projective if every module
surjection $P \to Q \to 0$ splits as a    right $A$-module map. $P$ is
called {\it finite} if there exists a surjection $A^n \to P \to 0$ for some integer $n$.

We describe the Serre-Swan correspondence between vector bundles and finite projective
modules. Given a vector bundle $p: E\rightarrow
X$ , let
$$P= \Gamma (E)=\{s: X\rightarrow E;\;  ps=id_X \}$$
be the set of all continuous {\it global sections}  of $E$. It is clear that under
fiberwise scalar multiplication and addition, $P$ is a $C(X)$ module. If
$f: E \to F$ is a bundle map, we define a module map $\Gamma (f): \Gamma
(E) \to \Gamma (F)$ by $\Gamma (f)(s)(x)=f(s(x))$ for all $s\in \Gamma
(E)$ and $x\in X$. We have thus defined a functor $\Gamma $, called the global section
functor, from the
category of vector bundles over $X$ and continuous bundle maps to the category of $C(X)$-modules and module maps.

Using compactness of $X$ and a partition of unity one shows
that there is a vector bundle $F$ on $X$ such that $E\oplus F\simeq X\times
\mathbb{C}^n$ is a trivial bundle. Let  $Q$ be  the space of global sections of
$F$.   We have
$$P \oplus Q\simeq A^n,$$
which shows that  $P$ is finite projective.

To show that all finite projective   $C(X)$-modules  arise in this way
we
proceed as follows.
 Given  a finite projective $C(X)$-module $P$, let $Q$ be a $C(X)$-module
such that $P \oplus Q\simeq  A^n,$ for some integer $n$. Let $e:A^n\rightarrow A^n$ be the right $A$-linear
map corresponding to the projection into first coordinate:  $(p, q)\mapsto (p, 0)$. It is
obviously an
idempotent in $M_n(C(X))$. One defines a vector  bundle $E$ as the image
of this idempotent $e$:
$$E=\{(x, v); \; e(x)v=v, \mbox{for all}\, x \in X, \,v\in \mathbb{C}^n\}\subset X\times \mathbb{C}^n.$$
Now it is easily shown that $\Gamma (E)\simeq P$. With some more work it is shown that the functor $\Gamma$
is full and faithful and hence defines  an equivalence of categories.

Based on the Serre-Swan theorem, one usually thinks of finite
projective modules over  noncommutative algebras as {\it
noncommutative vector bundles}. We give a few examples starting
with a commutative one.\\

{\bf Examples}\\

\noindent 1. The {\it Hopf line bundle} on the two sphere $S^2$,
also known as  {\it magnetic monopole bundle},  can be defined in
various ways. (It was discovered, independently, by Hopf and Dirac
in 1931, motivated by very different considerations).  Here is an
approach that lends itself to noncommutative generalizations. Let
$\sigma_1$, $\sigma_2$, $\sigma_3$, be three matrices in
$M_2(\mathbb{C})$ that satisfy the {\it canonical anticommutation
relations}:
$$\sigma_i \sigma_j +\sigma_j \sigma_i =2\delta_{ij},$$
for all $i, j=1, 2, 3.$ Here $\delta_{ij}$ is the Kronecker symbol. A
canonical choice is the so called {\it Pauli spin matrices}:
$$
\sigma_1=\left(
\begin{matrix}
0 & 1 \\
1 & 0
\end{matrix}
\right) \; , \quad \sigma_2 = \left(
\begin{matrix}
0 & i \\
-i & 0
\end{matrix}
\right) \; , \quad \sigma_3 = \left(
\begin{matrix}
1 & 0 \\
0 & -1
\end{matrix}
\right) \; .$$
Define a function
$$F: S^2 \rightarrow M_2(\mathbb{C}), \quad F(x_1, x_2,
x_3)=x_1\sigma_1 +x_2\sigma_2 +x_3\sigma_3,$$ where $x_1, x_2,
x_3$ are coordinate functions on $S^2$  so that $x_1^2+ x_2^2+
x_3^2=1$. Then $F^2 (x)=I_{2\times 2}$ for all $x\in S^2$ and
therefore
$$e=\frac{1+F}{2}$$ is an
idempotent in $M_2(C(S^2))$. It thus defines a complex vector bundle on $S^2$. We have,
$$ e(x_1, x_2, x_3)= \frac{1}{2} \left(
\begin{matrix}
1 + x_3 & x_1 + ix_2 \\
x_1-ix_2 & 1-x_3
\end{matrix}
\right).
$$

 Since
$$ rank\, F(x)= trace \,F(x)=1$$
for all $x\in S^2$, we have in fact a complex line bundle over $S^2$.
 It can be shown that it is the line bundle associated to the Hopf fibration
$$ S^1 \to S^3 \to S^2.$$

Incidentally, $e$ induces  a map $f: S^2 \to P^1 (\mathbb{C} )$, where $f(x)$
is the 1-dimensional subspace defined by the image of $F(x)$, which is 1-1 and onto. Our
line bundle is just the pull back of the canonical line bundle over $P^1 (\mathbb{C} )$.

This  example can be generalized to higher dimensional spheres.
One can construct matrices  $\sigma_1, \cdots, \sigma_{2n+1}$ in $M_{2^n} (\mathbb{C})$
satisfying  the {\it Clifford
algebra relations} \cite{kar}
$$\sigma_i\sigma_j + \sigma_j\sigma_i = 2\delta_{ij},$$
for all $i, j=1, \cdots, 2n+1.$
Define a matrix valued function $F$ on the $2n$ dimensional sphere $S^{2n}$,
$F\in M_{2^n} (C(S^{2n}))$ by
$$F=\sum_{i=1}^{2n+1} x_i \sigma_i.$$
Then $F^2(x) =1$ for all $x\in S^{2n}$ and hence $e=\frac{1+F}{2}$ is an
idempotent.  \\

\noindent 2. (Hopf line bundle on quantum spheres)\\
The Podle\'{s}  quantum sphere $S^2_q$ is the
$\ast$-algebra generated over $\mathbb{C}$ by the elements $a, a^{\ast}$ and $b$ subject to the
relations
\[
a a^{\ast}+q^{-4}b^2=1,\;a^{\ast}a+b^2=1,\;ab=q^{-2}ba,\;a^{\ast}b=q^2ba^{\ast}.
\]

The quantum analogue of the Dirac(or Hopf) monopole line bundle over $S^2$ is given by the following idempotent
in $M_2(S^2_q)$ \cite{bm}
$$\mathbf{e}_q=\frac{1}{2}\left[\begin{array}{cc} 1+q^{-2}b & qa \\ q^{-1} a^{\ast} & 1-b \end{array}\right].$$
It can be directly checked that $\mathbf{e}_q^2=\mathbf{e}_q.$\\

\noindent 3. (Projective modules on the noncommutative torus)\\
Let us first recall the definition of the {\it smooth noncommutative
torus} $\mathcal{A}_{\theta}$ \cite{ac80}. Among several possible
definitions the following is the most direct one. Let $\theta \in
\mathbb{R}$ be a fixed parameter and let
$$\mathcal{A}_{\theta}=\{\sum_{m, n \in \mathbb{Z}}a_{mn}U^mV^n;\;
(a_{mn})\in \mathcal{S}(\mathbb{Z}^2)\},$$
where $\mathcal{S}(\mathbb{Z}^2)$ is the Schwartz space of rapidly
decreasing
sequences $(a_{mn})\in \mathbb{C}$ indexed by $\mathbb{Z}^2$. The
relation
$$VU=e^{2\pi i\theta}UV,$$
defines an algebra structure on $\mathcal{A}_{\theta}$.

Let $E=\mathcal{S}(\mathbb{R})$ be the Schwartz space of rapidly
decreasing  functions on $\mathbb{R}$. It is easily checked that
the following formulas define an $\mathcal{A}_{\theta}$ module
structure on $E$:
$$ (Uf)(x)=f(x+\theta), \quad (Vf)(x)=e^{-2 \pi i \theta}f(x).$$
It can shown that $E$ is finite and projective \cite{ac80}.

This construction can be generalized \cite{ac80, cr, rie2}. Let $n, m$ be integers with $m>0$
and let $E_{n, m}=\mathcal{S}(\mathbb{R}\times \mathbb{Z}_m)$, where
$\mathbb{Z}_m$  is the cyclic group of order $m$. The following formulas
define an $\mathcal{A}_{\theta}$ module structure on $E_{n, m}$:
\begin{eqnarray*}
(Uf)(x, j)&=&f(x+\theta -\frac{n}{m}, j-1),\\
(Vf)(x, j)&=&e^{2 \pi i (x-j\frac{n}{m})}f(x, j).
\end{eqnarray*}

It can be shown that when $n-m\theta \neq 0$, the module $E_{n, m}$ is
finite and projective. In particular for irrational $\theta$ it is
always finite and projective.

For more examples of noncommutative vector bundles see \cite{cl}.

\subsection{Affine varieties and finitely generated commutative reduced algebras}
In commutative algebra, {\it Hilbert's Nullstellensatz} \cite{cg}
states that the category of affine algebraic varieties over an
algebraically closed field $\mathbb{F}$ is anti-equivalent to the
category of finitely generated commutative reduced unital
$\mathbb{F}$ algebras:
$$\{ \text{\bf affine algebraic varieties}\} \simeq $$
$$\{\text{\bf finitely generated  commutative
reduced  algebras}\}^{op}.$$

Recall that an {\it affine algebraic variety} (sometimes called an {\it algebraic set}) over a field
 $\mathbb{F}$ is a
subset of an affine space $\mathbb{F}^n$ which is the set of zeros
of a set of polynomials in $n$ variables over $\mathbb{F}$. A
morphism between  affine varieties $V \subset \mathbb{F}^n$ and
$W\subset \mathbb{F}^m$ is a map $f: V \longrightarrow W$  which
is the restriction of a polynomial map $\mathbb{F}^n \rightarrow
\mathbb{F}^m.$ It is clear that affine varieties and morphisms
between them form a category.

A {\it reduced algebra}
is by definition an algebra with no {\it nilpotent elements}, i.e. if $x^n=0$ for some $n$  then $x=0$.

The above correspondence associates to a variety $V\subset \mathbb{F}^n$ its {\it coordinate
ring} $\mathbb{F}[V]$ defined by
$$\mathbb{F}[V]:=Hom_{Aff}(V, \{pt\})\simeq  \mathbb{F}[x_1, \cdots, x_n]/I,$$
where $I$ is the {\it vanishing ideal} of $V$ defined by
$$I=\{f\in \mathbb{F}[x_1, \cdots, x_n]; f(x)=0 \;\mbox{for all}\; x\in V\}.$$
Obviously $\mathbb{F}[V]$ is a finitely generated commutative
unital reduced algebra. Moreover, Given  a morphism of varieties
$f: V \to W$, its pull-back defines an algebra homomorphism $f^*:
\mathbb{F}[W] \to \mathbb{F}[V]$. We have thus defined the
contravariant {\it coordinate ring functor}
 from
 affine varieties to finitely generated reduced commutative unital algebras.

 Given a  finitely generated commutative unital algebra $A$ with $n$ generators we can
 write it as
 $$A\simeq \mathbb{F}[x_1, \cdots, x_n]/I,$$
 where the ideal $I$ is a {\it radical ideal} if and only if $A$ is a reduced algebra. Let
 $$V:=\{x \in \mathbb{F}^n; f(x)=0\; \mbox{for all} \; f\in I\},$$
 denote the variety defined by the ideal $I$. The classical
 form of  Nullstellensatz \cite{har} states that if $\mathbb{F}$ is algebraically closed and $A$ is reduced
 then
 $A$ can be recovered
 as the
 coordinate ring of the variety $V$:
 $$\mathbb{F}[V]\simeq A=\mathbb{F}[x_1, \cdots, x_n]/I.$$
 This is the main step in showing that the coordinate ring functor is an anti-equivalence of categories.
 Showing that the functor is full and faithful is easier.
 In Section 6 we sketch a proof of this fact when  $\mathbb{F}$ is the
 field of complex numbers.

\subsection{Affine schemes and commutative rings}

The above correspondence between finitely generated reduced
commutative algebras and affine varieties in not an ideal result.
One  is naturally interested in larger classes of algebras, like
algebras  with nilpotent elements as well as algebras over fields
which are not algebraically closed or algebras over arbitrary
rings; this last case in particularly important in number theory.
In general one wants to know what kind of geometric objects
correspond to a commutative ring and how this correspondence goes.
{\it Affine schemes} are exactly defined to address this question.
We follow the exposition in \cite{har}.

Let $A$ be a commutative unital ring. The {\it prime spectrum} (or
simply the {\it spectrum}) of $A$ is a pair $(\mbox{Spec} A,
\mathcal{O}_A)$ where $\mbox{Spec} A$ is a topological space  and
$\mathcal{O}_A$ is a sheaf of rings on $\mbox{Spec} A$ defined as
follows. As a set $\mbox{Spec} A$ consists of all {\it prime
ideals} of $A$ (an ideal $I\subset A$ is called {\it prime } if
for all $a, b$ in $A$, $ab \in A$ implies that either $a\in I$, or
$b\in I$).  Given an ideal $I\subset A$, let $V(I)\subset
\mbox{Spec} A$ be the set of all prime ideals which contain $I$.
We can define a topology on $\mbox{Spec} A$, called the {\it
Zariski topology},  by declaring sets of the type $V(I)$ to be
closed (this  makes sense since the easily established relations
$V(IJ)=V(I) \cup V(J)$ and $V(\sum I_i)= \cap V(I_i)$ show that
the intersection of a family of closed sets is closed and the
union of two closed sets is closed as well). One checks that
$\mbox{Spec} A$ is always compact but is not necessarily
Hausdorff.

For each prime ideal $p\subset A$, let $A_p$ denote the {\it
localization} of $A$ at $p$. For an open set $U\subset \mbox{Spec} A$, let $\mathcal{O}_A(U)$
be  the set of all continuous sections $s: U \to \cup_{p\in U} A_p$. (By definition a section $s$ is called
continuous if locally around any point $p\in U$ it is of the form
$\frac{f}{g}$, with $g\notin p$). One checks that
 $\mathcal{O}_A$ is a sheaf of commutative rings on $\mbox{Spec} A$.

Now $(\mbox{Spec} A, \mathcal{O}_A)$ is a so called {\it ringed space}
and
$A \mapsto (\mbox{Spec} A, \mathcal{O}_A)$ is functor called the {\it spectrum functor}.
 A unital ring homomorphism $f: A \to B$ defines a continuous map $f^*:
\mbox{Spec} B \to \mbox{Spec} A$ by $f^*(p)= f^{-1}(p)$ for all prime
ideals $p\subset B$. (note that if $I$ is a maximal ideal $f^{-1}(I)$ is not necessarily maximal.
This is one of the reasons one considers, for arbitrary rings,  the
prime spectrum  and not
the maximal spectrum.)

An {\it affine scheme} is a ringed space $(X, \mathcal{O})$ such that
$X$ is homeomorphic to $\mbox{Spec} A$ for a commutative ring $A$ and $\mathcal{O}$ is isomorphic to
$\mathcal{O}_A$. The spectrum functor defines an equivalence of
categories
$$\{\mbox{\bf affine schemes}\} \simeq \{\mbox{ \bf commutative rings}\}^{op} .$$

The inverse equivalence is given by the {\it global section
functor} that sends an affine scheme  to the ring  of its global
sections.

In the same vein categories of modules over a ring can be
identified with categories of sheaves of modules over the spectrum of the ring.
Let $A$ be a commutative ring and let $M$ be an $A$-module. We define a sheaf of
modules $\mathcal{M} $ over $\mbox{Spec} A$ as follows.  For each prime  ideal
$p\subset A$, let $M_p$ denote the localization of $M$ at $p$. For any
open set $U\subset \mbox{Spec} A$ let $\mathcal{M}(U)$ denote the set of
continuous  sections $s: U\longrightarrow \cup_p M_p$ (this means that $s$ is locally a
fraction $\frac{m}{f}$ with $m\in M$ and  $f\in A_p$). One can recover $M$
from $\mathcal{M}$ by showing that $M\simeq \Gamma \mathcal{M}$ is the space of global
sections of $M$. The functors $M\mapsto \mathcal{M}$ and $\mathcal{M}\mapsto  \Gamma \mathcal{M}$
define equivalence of categories \cite{har}:
$$\{\mbox{\bf modules over $A$}\} \simeq \{\mbox{\bf quasi-coherent sheaves on Spec $A$}\}.$$

\subsection{Compact Riemann surfaces and algebraic function fields}

It can be shown that the category of compact Riemann surfaces is
anti-equivalent to the category of algebraic function fields:
$$\{ \text{\bf compact Riemann surfaces} \}\simeq \{ \text{\bf algebraic function fields}\}^{op}.$$

Recall that a {\it Riemann surface} is a complex manifold of complex
dimension one. A morphism between Riemann surfaces $X$ and $Y$ is a holomorphic
map $f: X \rightarrow Y$.

An {\it algebraic function field} is a finite extension of the
field $\mathbb{C}(x)$ of rational functions in one variable. A
morphism of function fields is simply an algebra map.

To a compact Riemann surface one associates the field  $M(X)$ of
meromorphic functions on $X$. For example the field of meromorphic
functions on the Riemann sphere is the field of rational functions
$\mathbb{C}(x)$. In the other direction, to    a finite extension of $\mathbb{C}(x)$ one associates
the compact Riemann surface of the algebraic function $p(z, w)=0.$ Here
$w$ is a generator of the field over $\mathbb{C}(x)$. This
correspondence is essentially due to Riemann.

\subsection{Sets and Boolean algebras}
Perhaps the simplest notion of space, free of any extra structure, is
the concept of  a set. In a sense set theory can be regarded as the
geometrization of logic. There is a duality between the category of sets and the
category of complete atomic  Boolean algebras (see, e.g., M. Barr's Acyclic Models,
 CRM Monograph Series, Vol 17, AMS
publications, 2002):
$$ \{\text{\bf{sets}}\}\simeq \{ \text{\bf{complete atomic Boolean algebras}}\}^{op}.$$

Recall that  a {\it Boolean algebra} is a unital ring $B$ in which
$x^2=x$ for all $x$ in $B$.  A Boolean algebra  is necessarily
commutative as can be easily shown. One defines an order relation
on $B$ by declaring $x\leq y$ if there is an $y'$ such that
$x=yy'$. It can be checked that this is in fact a partial  order
relation on $B$. An {\it Atom} in a Boolean algebra is an element
$x$ such that $x>0$ and there is no $y$ with $0<y<x$. A Boolean
algebra is {\it atomic} if every element $x$ is the supremos of
all the atoms smaller than $x$. A Boolean algebra is {\it
complete} if every subset has a supremum and infimum. A morphism
of complete Boolean algebras is a unital ring map which preserves
all infs and sups. (Of course any unital ring map between  Boolean
algebras preserves finite sups and infs).

Now, given a set $S$ let
$$B={\bf 2}^S= \{f: S \longrightarrow {\bf 2}\},$$
where ${\bf 2}:=\{0, 1\}.$ Note that $B$ is a complete atomic Boolean algebra.
Any map $f: S \rightarrow T$ between sets defines a morphism of
complete atomic Boolean algebras via pullback: $ f^*(g):=g\circ f$,
and
$ S \mapsto {\bf 2}^S$
is a contravariant functor.

In the opposite direction, given a Boolean algebra $B$, one defines its {\it spectrum} $\Omega (B)$
by
$$\Omega (B)= Hom_{Boolean} (B, {\bf 2}),$$
where we now think of {\bf 2} as a Boolean algebra with two elements. It can be shown  that the two
functors that
we have defined are anti-equivalences
of categories, quasi-inverse to each other. Thus once again we have a duality between
a certain category of geometric
objects, namely sets, and
a category of commutative algebras, namely complete atomic Boolean algebras.

\section{ Noncommutative quotients}

In this section we recall the method of noncommutative quotients as advanced by Connes in \cite{acb}.
This is a technique that allows one to replace ``bad quotients'' by nice noncommutative
spaces, represented by noncommutative algebras. In some cases, like noncommutative quotients for group actions,
the noncommutative quotient can be defined as a crossed product
algebra. In general, however, noncommutative quotients are defined
as groupoid algebras. In Section 3.1 we recall the definition of a
groupoid together with its various refinements like topological,
smooth and \'{e}tale groupoids. In Section 3.2 we define the
groupoid algebra of a groupoid and  give several examples. An
important concept is Morita equivalence of algebras. We treat both
the purely algebraic theory as well as the notion of strong Morita
equivalence for $C^*$-algebras in Section 3.3. Finally
noncommutative quotients are defined in Section 3.4.

\subsection{ Groupoids}
\begin{definition}
A groupoid is a small category in which every morphism is an
isomorphism.
\end{definition}

Let $\mathcal{G}$ be groupoid. We denote the set of objects of
$\mathcal{G}$ by $\mathcal{G}^{(0)}$ and, by a small abuse of notation,  the set of morphisms of $\mathcal{G}$
by $\mathcal{G}$. Every morphism has a {\it source}, has a {\it
target} and has an {\it inverse}. They define maps, denoted by $s$, $t$,
and $i$, respectively,
$$ s:\mathcal{G} \longrightarrow \mathcal{G}^{(0)}, \quad
t:\mathcal{G} \longrightarrow \mathcal{G}^{(0)},$$
$$i: \mathcal{G} \longrightarrow \mathcal{G}.$$

Composition $\gamma_1 \circ \gamma_2$ of morphisms $\gamma_1 $ and $
\gamma_2$ is only defined if $ s(\gamma_1)=t(\gamma_2)$. Composition
defines a map
$$ \circ:  \mathcal{G}^{(2)}=\{(\gamma_1, \gamma_2); s(\gamma_1)=
t(\gamma_2)\}\longrightarrow \mathcal{G}.$$

\noindent {\bf Examples}\\
1. (Groups). Every group $G$ defines at least two groupoids in  a
natural way:\\
1.a  Define a category $\mathcal{G}$ with one object * and
$$Hom_{\mathcal{G}}(*, *)=G,$$
where the composition of morphisms is simply the group multiplication.
This is obviously a groupoid. \\

2.b Define a category $\mathcal{G}$ with
$$obj \;  \mathcal{G}=G, \quad
Hom_{\mathcal{G}} (s, t)=\{g\in G; \; gsg^{-1}=t\}.$$
Again, with composition defined by group multiplication, $\mathcal{G}$ is a groupoid.

2. (Equivalence relations) Let $\sim$ denote an equivalence relation on
a set $X$. We define a groupoid $\mathcal{G}$, called the {\it graph of $\sim$}, as
follows. Let
\begin{eqnarray*}
obj \; \mathcal{G}= X, \quad Hom_{\mathcal{G}}(x, y)&=&{*} \quad
\mbox{if}\;
x\sim y, \\
&=& \varnothing \quad  \mbox{otherwise.}
\end{eqnarray*}

Note that the set of morphisms of $\mathcal{G}$ is identified with the
graph of the relation $\sim$ in the usual sense:
$$\mathcal{G}=\{ (x, y); \; x\sim y \}\subset X\times X.$$

Two extreme cases of this graph construction are particularly important.
When the equivalence relation reduces to equality, i.e.,
$x\sim y$ iff $x=y$, we have
$$ \mathcal{G}= \Delta (X)=\{(x, x); \;x\in X\}.$$
On the other extreme when $x\sim y$ for all $x$ and $y$, we obtain the
{\it groupoid of pairs} where
$$ \mathcal{G}= X\times X.$$

3. (Group actions). Example 1) can be generalized as follows. Let
$$G \times X \longrightarrow X, \quad (g, x)\mapsto gx,$$
denote the action of a group $G$ on a set $X$. We define a groupoid
$\mathcal{G}=X\rtimes G$, called the {\it transformation  groupoid} of the action, as
follows. Let $obj \;\mathcal{G} =X$, and
$$Hom_{\mathcal{G}}(x, y)=\{g\in G; \; gx=y\}.$$
Composition of morphisms is  defined via group multiplication. It is
easily checked that $\mathcal{G}$ is a groupoid. Its set of
morphisms can be identified as
$$\mathcal{G}\simeq X \times G,$$
where the  composition of morphisms is given by
$$(gx, h) \circ (x, g)= (x, hg).$$

Note that Example 1.a corresponds to the action of a group on a point
and example 1.b corresponds to the action of a group on itself via conjugation.

As we shall see, one can not get very far with just discrete
groupoids. To get really interesting examples like the groupoids
associated to continuous actions of topological groups and to
foliations, one needs to consider topological as well as smooth
groupoids, much in the same way as one studies topological and Lie
groups.

A {\it topological groupoid} is a groupoid such that its
set of morphisms $\mathcal{G}$, and set of objects
$\mathcal{G}^{(0)}$ are topological spaces, and its  composition, source, target and
inversion maps are continuous.

A special class of topological groupoids, called \'{e}tale or r-discrete groupoids,
are  particularly convenient to work with. An {\it \'{e}tale groupoid} is a topological groupoid
such that its set of morphisms $\mathcal{G}$ is a locally compact topological space
and the fibers of the target map $t:\mathcal{G}\to \mathcal{G}^{(0)}$
$$ \mathcal{G}^x= t^{-1}(x), \quad x\in \mathcal{G}^{(0)},$$
are discrete. 

A {\it Lie groupoid} is a groupoid such that $\mathcal{G}$ and
$\mathcal{G}^{(0)}$ are smooth manifolds, the inclusion
$\mathcal{G}^{(0)} \to \mathcal{G}$ as well as the maps $s, t, i$ and the
composition map $\circ$ are
smooth, and $s$ and $t$ are submersions. This last condition will gaurantee that the domain of the
composition map $\mathcal{G}^{(2)}=\{(\gamma_1, \gamma_2); s(\gamma_1)=
t(\gamma_2)\}$ is a smooth manifold.\\

\noindent {\bf Examples:}\\

\noindent 1. Let $G$ be a discrete group acting by homeomorphisms on a locally
compact Hausdorff space $X$. The transformation groupoid $X\rtimes G$ is
naturally an \'{e}tale groupoid. If $G$ is a Lie group acting smoothly on a smooth manifold $X$, then
the transformation groupoid $X\rtimes G$ is a Lie groupoid. \\

\noindent 2. Let $(V, \mathcal{F})$ be a foliated manifold and let $T$
be a complete transversal for the foliation. This means that $T$ is transversal to the leaves
of the foliation and each  leaf has at least one intersection with $T$.
One defines an (smooth) \'{e}tale groupoid $\mathcal{G}$ as follows. The objects of $\mathcal{G}$  is the
transversal $T$ with its smooth structure. For any two points $x$ and $y$ in $T$ let
$Hom_{\mathcal{G}} (x, y)=\varnothing$ if $x$ and $y$ are not in the same leaf. When
they are in the same leaf, say $L$,  let  $Hom_{\mathcal{G}} (x, y)$ denote the set of all continuous
paths in $L$ connecting $x$ and $y$ modulo the equivalence relation defined by {\it holonomy}.  It can be shown that
$\mathcal{G}$ is an smooth \'{e}tale groupoid (cf. \cite{acb} for details and
many examples).

\subsection{Groupoid algebras}

The notion of {\it groupoid algebra} of a groupoid is a
generalization of the notion of {\it group algebra} (or {\it
convolution algebra}) of a group and it reduces to group algebras
for groupoids with one object. To define the groupoid algebra of a
locally compact topological groupoid in general one needs the analogue of a
Haar measure for groupoids called a Haar system. While we won't recall its general definition here,
we should
mention that,
 unlike locally compact groups, an arbirtary locally compact groupoid need not have a Haar system
 \cite{ren}.  For discrete
groupoids as well as \'{e}tale groupoids and Lie groupoids,  however, the convolution
product can be easily defined. We start by recalling the
definition of the groupoid algebra of a discrete groupoid. As we
shall see, in the discrete case the groupoid algebra can be easily
described in terms of matrix algebras and group algebras.

Let $\mathcal{G}$ be a discrete groupoid and let
$$ \mathbb{C} \mathcal{G}=\bigoplus_{\gamma \in \mathcal{G}} \mathbb{C} \gamma,$$
denote the vector space generated by the set of morphisms of
$\mathcal{G}$ as its basis. The formula
\begin{eqnarray*}
\gamma_1 \gamma_2 &=&\gamma_1 \circ \gamma_2, \quad \mbox{if} \; \gamma_1
\circ \gamma_2 \quad \mbox{is defined},\\
 \gamma_1 \gamma_2&=&0, \quad \mbox{ otherwise,}
\end{eqnarray*}
  defines an associative product on $\mathbb{C} \mathcal{G}$. The
  resulting algebra is called the {\it groupoid algebra} of the groupoid
  $\mathcal{G}$. Note that $\mathbb{C} \mathcal{G}$ is unital if and
  only if the set $\mathcal{G}^{(0)}$ of objects of $\mathcal{G}$ is finite. The unit then  is given by
  $$ 1=\sum_{x \in \mathcal{G}^{0}}  id_x.$$

An alternative description of the groupoid algebra $\mathbb{C}
\mathcal{G}$ which is more appropriate for generalization to
\'{e}tale and topological groupoids is as follows. Note that
$$\mathbb{C} \mathcal{G} \simeq \{ f: \mathcal{G} \to \mathbb{C}; \;f
\mbox{ has finite support}\},$$ and the product is given by the
{\it convolution product}
$$ (fg) (\gamma)=\sum_{\gamma_1 \circ \gamma_2
=\gamma}f(\gamma_1)g(\gamma_2).$$

Given an \'{e}tale  groupoid $ \mathcal{G}$, Let
$$C_c(\mathcal{G})=\{ f: \mathcal{G} \to \mathbb{C}; \;
\mbox{ f is continuous and has compact support} \}.$$
Under the above convolution product  $C_c(\mathcal{G})$ is an
algebra called the {\it  convolution algebra} of $\mathcal{G}.$ Note
that for any $\gamma$ the above sum is finite (why?) and the convolution
of two functions with compact support has compact support.

When the groupoid is the transformation groupoid of
a group action, the  groupoid algebra reduces to a {\it crossed product
algebra}. We recall its general definition.
Let $G$ be a discrete group, $A$ be an algebra and let $Aut (A)$ denote the group
of automorphisms of $A$.
An {\it action} of $G$ on $A$ is   a group
homomorphism
$$\alpha :G\longrightarrow {\text Aut} (A).$$
Sometimes one refers to the triple $(A, G, \alpha )$ as a {\it
noncommutative dynamical system}. We use the simplified notation
$g(a):=\alpha (g)(a).$
The {\it crossed product} or {\it semidirect product} algebra $A\rtimes G$ is
defined as follows. As a vector space
$$A\rtimes G =A\otimes \mathbb{C}G.$$
Its  product is defined by
$$(a\otimes g)(b\otimes h)=ag(b)\otimes gh.$$
It is easily checked that endowed with the above product, $A\rtimes G$ is
an associative algebra. It is unital if and only if $A$ is unital and
$G$ acts by unital automorphisms.

One  checks that $A\rtimes G$ is the universal algebra generated by
subalgebras $A$ and $\mathbb{C}G$ subject to the relation
$$gag^{-1}=g(a),$$
for all $g$ in $G$ and $a$ in $A$.

We need a $C^*$-algebraic  analogue of the above construction.
 In the above situation assume that $A$ is a
$C^*$-algebra and let $Aut (A)$ denote the group of $C^*$-automorphisms of $A$.  We define a pre-$C^*$ norm on the algebraic crossed
product $A\rtimes G$ as follows. We choose a faithful representation
$$\pi : A\longrightarrow \mathcal{L}(H)$$
of
$A$ on a Hilbert space $H$
and then define a faithful representation of the algebraic crossed product $A\rtimes G$ on the
Hilbert space $l^2 (G, H)$ by
$$(a\otimes g)(\psi)(h)= \pi(a)\psi(g^{-1}h). $$
The {\it reduced $C^*$-crossed product} $A\rtimes_r G$ is by definition the completion of the
algebraic crossed product with  respect to the norm induced from the above faithful representation.
It can be shown that this definition is independent of the choice of the
faithful representation $\pi$ \cite{D}.\\

\noindent{\bf Examples:}\\

\noindent 1. It is not difficult to determine the general form of the groupoid
algebra of a discrete groupoid. Let
$$\mathcal{G}^0=\bigcup_i \mathcal{G}^0_i$$
denote the decomposition of the set of objects of $\mathcal{G}$ into its
``connected components". By definition two objects $x$ and $y$ belong to
the same connected component if there is a morphism $\gamma$ with source
$(\gamma)=x$ and target $(\gamma)=y$. We have a direct sum
decomposition of the  groupoid algebra $\mathbb{C}\mathcal{G}$:
$$\mathbb{C}\mathcal{G}\simeq \bigoplus_i \mathbb{C}\mathcal{G}_i.$$

Thus suffices to consider only groupoids with a connected set of objects. Choose an
 object $x_0\in \mathcal{G}^0$, and let
 $$G=Hom_{\mathcal{G}}(x_0, x_0)$$
 be the {\it isotropy group} of $x_0$. Assume first that $G=\{1\}$  is the trivial group. For
 simplicity assume that $\mathcal{G}^0$ is a finite set with $n$ elements. In other words our
 groupoid $\mathcal{G}$
 is the groupoid of pairs  on a set of $n$ elements. Recall that
 $$\mathcal{G}=\{(i, j); \; i, j=1, \cdots n\}$$
 with composition given by
 $$(l, k)\circ (j, i)= (l, i) \quad  \mbox{if}\; k=j.$$
 (Composition is not defined otherwise).

  We claim that
 $$\mathbb{C}\mathcal{G}\simeq M_n(\mathbb{C}).$$
 Indeed it is easily checked that the map
 $$ (i, j)\mapsto E_{i, j},$$
 where $E_{i, j}$ denote the matrix units defines an algebra isomorphism.

 \begin{remark} As is emphasized by Connes in the opening section of \cite{acb}, this is in fact the way Heisenberg
 discovered matrices
 in the context of quantum mechanics \cite{hei}. In other words
 noncommutative algebras appeared first in quantum mechanics as a
 groupoid algebra! We recommend the reader carefully examine the
 arguments of \cite{hei} and \cite{acb}.
 \end{remark}

In general, but still assuming that $\mathcal{G}^{(0)}$ is connected, it
is easy to see that
$$\mathbb{C}\mathcal{G}\simeq \mathbb{C} G\otimes M_n(\mathbb{C}),$$
where $G=Hom_{\mathcal{G}}(x_0, x_0)$ is the isotropy group of $\mathcal{G}$. \\

\noindent 2. We look at groupoid algebras of certain
\'{e}tale groupoids. \\

\noindent a. We start by an example from \cite{acb}:  an \'{e}tale
groupoid defined by an equivalence relation. Let
$$X=[0, 1]\times \{1\} \cup [0, 1]\times \{2\}$$
 denote the disjoint
union of two copies of the interval $[0, 1].$ Let $\sim$ denote
the equivalence relation that identifies $(x, 1)$ in the first
copy  with $(x, 2) $ in the second copy for $0<x<1$. Let
$\mathcal{G}$ denote the corresponding groupoid  with its topology
inherited from $X\times X$. It is clear that $\mathcal{G}$ is an
\'etale groupoid.  The elements of the groupoid algebra
$C_c(\mathcal{G})$ can be identified as continuous matrix valued
functions on $[0, 1]$ satisfying some boundary condition:
$$ C_c(\mathcal{G})=\{f: [0, 1]\to M_2(\mathbb{C}); \quad f(0)\;
\mbox{and}\; f(1)\; \mbox{are diagonal} \}.$$

\noindent b. Let $X$ be a locally compact Hausdorff space and assume a discrete group $G$ acts on
$X$ by homemorphisms. The
induced action on $A=C_c(X)$, the algebra of continuous $\mathbb{C}$-valued functions on $X$ with
compact support,   is defined by
$$(gf)(x)=f(g^{-1}x)$$
for all $f\in A$, $g\in G$, and $x\in X$.
Let $\mathcal{G}=X\rtimes G$ denote the transformation
groupoid defined by the action of  $G$ on
 $X$.\\
 Exercise: Show that we have an algebra isomorphism
 $$C_c(X\rtimes G)\simeq C_c(X)\rtimes G.$$

 The groupoid $C^*$-algebra $C^*_r(\mathcal{G})$, though we have not defined it as such, turns out
 to be   isomorphic with the reduced crossed
product algebra $C_0(X)\rtimes_r G$.\\

\noindent 3.  Let $X$ be a locally compact Hausdorff space and let $\mathcal{G}$ denote the groupoid
of pairs on $X$. Its groupoid $C^*$-algebra
$C^* \mathcal{G}$ is defined as follows. Let $\mu$ be a positive Borel
measure on $X$. We define a convolution product and an $*$-operation on the space of
morphisms of $\mathcal{G}$,
$C_c( X\times X)$, by
\begin{eqnarray*}
f*g(x, y)&=&\int_X f(x,z) g(z, y)d\mu (z),\\
f^*(x, y)&=& \overline{f(y,x)}.
\end{eqnarray*}
These operations turn $C_c(X\times X)$ into an $*$-algebra. Convolution
product defines a canonical $*$-representation of this algebra on
$L^2(X, \mu)$ by
$$C_c(X\times X)\longrightarrow \mathcal{L}(L^2(X, \mu)), \quad f\mapsto
f*-.$$
The integral operator associated to a continuous function with compact
support is a compact operator and it can be shown that the completion of the image of this
map is the space of compact operators on $L^2(X, \mu).$ Thus we have 
$$C^* \mathcal{G} \simeq \mathcal{K}(L^2 (X, \mu).$$

In the other extreme, for the 
equivalence relation defined by equality the groupoid $C^*$-algebra is given by
$$C^* \mathcal{G}\simeq C_0(X).$$

These two examples  are continuous $C^*$-analogues of our  discrete example 1 above.

\subsection{ Morita equivalence}
\subsubsection{Algebraic theory}
In this section algebra means an associative unital algebra over a
commutative ground ring $k$. All modules are assumed to be unitary in
the sense that the unit of the algebra acts as identity operator on the module. Let $A$ be an
algebra. We denote by
$\mathcal{M}_A$ the category of right $A$-modules.

Algebras $A$ and $B$ are called {\it Morita equivalent} if there is an equivalence of categories
$$ \mathcal{M}_A\simeq \mathcal{M}_B.$$

In general there are many ways to define  a functor $F: \mathcal{M}_A\to \mathcal{M}_B$. By a simple observation of Eilenberg-Watts,
however, if $F$ preserves finite limits and colimits, then there exists a unique $A-B$ bimodule $X$
such that
$$F(M)= M\otimes_A X, \quad \text{for all}\quad M\in \mathcal{M}_A.$$
Composition of functors obtained in this way simply correspond to the balanced tensor product of the defining
bimodules.
It is therefore clear that $A$ and $B$ are Morita equivalent if and only
if there exists an $A-B$ bimodule $X $ and a $B-A$ bimodule $Y$ such that we have isomorphisms of bimodules
$$ X\otimes_B Y \simeq A, \quad Y\otimes_A X\simeq B,$$
where the $A$-bimodule structure on $A$ is defined by $a(b)c=abc$, and similarly for $B$.
Such
bimodules are called {\it invertible (or equivalence) bimodules}.

Given an $A-B$ bimodule $X$, we define algebra homomorphisms
$$A \longrightarrow End_B(X), \quad B\longrightarrow End_A(X),$$
$$a\mapsto L_a, \quad b\mapsto R_b,$$
where $L_a$ is the operator of left multiplication by $a$ and $R_b$ is the operator of
right multiplication by $b$.

In general it is rather hard to characterize the invertible bimodules.
The following theorem is one of the main results of Morita:
\begin{theorem}
An $A-B$ bimodule $X$ is invertible if and only if $X$ is finite and
projective both as a left $A$-module and as a right $B$-module and the
natural maps
$$A \rightarrow End_B(X), \quad B\rightarrow End_A(X),$$
are algebra isomorphisms.
\end{theorem}

{\bf Example}. Any unital algebra $A$ is Morita equivalent to the
algebra $M_n(A)$ of $n\times n$ matrices over $A$. The $(A, M_n(A))$
equivalence bimodule is $X=A^n$ with obvious left $A$-action and   right $M_n(A)$-action. This example
can be generalized as follows.

{\bf Example}. Let $P$ be a finite projective left $A$-module and  let
$$B=End_A(P).$$
Then the algebras $A$ and $B$ are Morita equivalent. The
equivalence $A-B$ bimodule is given by $X=P$ with obvious $A-B$
bimodule structure. As a special case, we obtain the following
geometric example.

{\bf Example} Let $X$ be a compact Hausdorff space and $E$ a
complex vector bundle on $X$. Then the algebras $A= C(X)$ of
continuous functions on $X$ and $B= \Gamma (End (E))$ of global
sections of the endomorphisms bundle of $E$ are Morita equivalent.
In fact, in view of Swan's theorem  this is a special case of the
last example with $P= \Gamma (E)$ the global sections of $E$.
There are analogous results for real as well as quaternionic
vector bundles. If $X$ happens to be a smooth manifold we can let
$A$ to be the algebra of smooth functions on $X$ and $B$ be the
algebra of smooth sections of $End(E)$.

Given a category $\mathcal{C}$ we can consider the category $Fun
(\mathcal{C})$ whose objects are  functors from $\mathcal{C} \rightarrow
\mathcal{C}$ and whose morphisms are natural transformations between
functors. The {\it center} of a category $\mathcal{C}$ is by definition
the set of natural transformations from the identity functor to itself:
$$\mathcal{Z}(\mathcal{C}):=Hom_{Fun (\mathcal{C})} (Id,Id).$$
Equivalent categories obviously have  isomorphic  centers.

Let $\mathcal{Z} (A)=\{a\in A; \,ab=ba\,
\mbox{for all}\, b\in A\}$ denote the center
of an algebra $A$. It is easily
seen that for $\mathcal{C}=\mathcal{M}_A$ the natural map
$$\mathcal{Z}(A) \rightarrow \mathcal{Z}(\mathcal{C}), \quad a\mapsto R_a, $$
where $R_a(m)=ma$ for any module $M$ and any $m\in M$, is one to one and onto.
It follows that  Morita equivalent algebras have isomorphic centers:
$$ A \overset{M.E.}{\sim}B \; \Rightarrow  \; \mathcal{Z}(A)\simeq \mathcal{Z}(B). $$
 In particular two commutative algebras are Morita equivalent if and only
if they are isomorphic. We say that commutativity is not a {\it Morita invariant property}. \\

 Exercise: Let $A$ be a unital $k$-algebra. Show that there is a 1-1
 correspondence between space of traces on $A$ and $M_n(A)$. Extend this fact to arbitrary Morita
 equivalent algebras.

We will see in   Section 4 that Morita equivalent algebras have
isomorphic Hochschild and cyclic (co)homology groups. They have isomorphic algebraic $K$-theory as
well.

\subsubsection{ Strong Morita equivalence}
Extending the Morita theory to non-unital algebras and to topological
algebras needs  more work and is not an easy task.  For
(not
necessarily unital) $C^*$-algebras we have Rieffel's notion of {\it strong
Morita equivalence} that we recall below.

For $C^*$-algebras one is mostly interested in their $*$-representations on  a Hilbert space.
Thus one must consider equivalence $A-B$ bimodules $X$ such that if $H$ is a Hilbert
space and a right $A$-module, then $H\otimes_A X$ is a Hilbert space as
well. This leads naturally to the concepts of Hilbert module and Hilbert
bimodule that we recall below.

Let $B$ be a not necessarily unital $C^*$-algebra. A right {\it Hilbert module} over $B$ is a
right $B$-module $X$ endowed with a $B$-valued inner product such that
$X$ is complete with respect to its natural norm. More precisely, we
have a sesquilinear map
$$X\times X \longrightarrow B, \quad (x, y) \mapsto <x, y>,$$
such that for all $x, y$ in $X$ and $b$ in $B$ we have
$$<x, y> =<y, x>^*, \, <x, yb>=<x, y>b, \,\mbox{and}\, <x, x> >0 \,\mbox{for}\, x\neq 0.$$
It can be shown that $\| x\|:=\| <x, x>\|^{1/2}$ is a norm on $X$. We
assume $X$ is complete with respect to this norm.

Of course, for $B=\mathbb{C}$, a Hilbert $B$-module is just a Hilbert space. A very simple geometric example to keep in  mind is the following. Let
$M$ be a compact Hausdorff space and let $E$ be a complex vector bundle on $M$ endowed with a Hermitian
inner product. One defines a Hilbert module structure on the space
$X=\Gamma (E)$ of continuous sections of $E$ by
$$ <s, t>(m)= <s(m), t(m)>_m$$
for continuous sections $s$ and $t$ and $m\in M$.

A morphism of Hilbert $B$-modules $X$ and $Y$ is a bounded $B$-module map $X\to Y$.
Every bounded operator on a Hilbert space has an adjoint. This is not
the case for Hilbert modules. (This is simply because, even purely algebraically,
a submodule of a module need not have a complementary submodule).
A bounded $B$-linear map $T: X\rightarrow
X$ is called {\it adjointable} if there is a bounded $B$-linear map $T^*:
X\rightarrow X$ such that for all $x$ and $y$ in $X$ we have
$$ <Tx, y>=<x, T^*y>.$$
Let $\mathcal{L}_B(X)$ denote the algebra of bounded  adjointable
$B$-module maps $X\rightarrow X$. It is a $C^*$-algebra.

Let $A$ and $B$ be $C^*$-algebras. A {\it Hilbert
$A-B$ bimodule} consists of a right Hilbert $B$-module $X$ and a
 $C^*$ map
$$ A \longrightarrow \mathcal{L}_B(X).$$

$C^*$-algebras $A$ and $B$ are called {\it strongly Morita equivalent}
if there is a  Hilbert $A-B$ bimodule $X$ and a Hilbert $B-A$ bimodule
$Y$ such that we have isomorphisms of bimodules
$$ X\otimes_B Y \simeq A, \quad \quad Y\otimes_A X\simeq B.$$
(The tensor products are completions of algebraic tensor products with respect to thier natural
pre-Hilbert module structures.)

Two  unital $C^*$-algebras are strongly Morita equivalent if and
only if they are Morita equivalent as algebras. Any $C^*$-algebra $A$ is
Morita equivalent to its {\it stabilization} $A\otimes \mathcal{K}$,
where $\otimes$ is the $C^*$ tensor product and $\mathcal{K}$ is the
algebra of compact operators on a Hilbert space. Strongly Morita equivalent algebras have naturally
isomorphic topological $K$-theory. (cf. \cite{gvf} for more details, and a proof of these statements or references to
original sources.)

\subsection{ Noncommutative quotients}

>From a purely set theoretic point of view, all one needs to form a
quotient space $X/\sim$ is an equivalence relation $\sim $ on a
set $X$. The equivalence relation however is usually obtained from
a much richer structure by forgetting part of this structure. For
example, $\sim$ may arise from an action of a group $G$ on $X$
where $x\sim y$ if and only if $gx=y$ for some $g$ in $G$ ({\it
orbit equivalence}). Note that there may be, in general, many $g$
with this property. That is  $x$ may be identifiable with $y$ in
more than one way. Of course when we form the equivalence relation
this extra information is lost. The key idea in dealing with bad
quotients in Connes' theory  is to keep track of this extra
information!

We call, rather vaguely, this extra structure the {\it quotient data}. Now Connes's dictum in forming
noncommutative quotients can be summarized as follows:
$$ \mbox{\bf quotient data} \rightsquigarrow \mbox{\bf groupoid} \rightsquigarrow \mbox{\bf groupoid
algebra}, $$
where the noncommutative quotient is defined to be the groupoid algebra itself.

Why is this a reasonable approach? The answer is that first of
all, by a theorem of M. Rieffel (see below) when the classical
quotient, defined by a group action, is a reasonable space, the
algebra of continuous functions on the classical quotient is
strongly Morita equivalent to the groupoid algebra.  Now it is known that Morita
equivalent algebras have isomorphic $K$-theory, Hochschild  and
cyclic (co)homology groups. Thus the topological invariants
defined via noncommutative geometry are the same for the two
constructions and no information is lost.

For bad quotients there is no reasonable space but we think of the
noncommutative algebra defined as a groupoid algebra as
representing a noncommutative quotient space. Thanks to
noncommutative geometry, tools like K-theory, K-homology,  cyclic
cohomology and the local index formula, etc., can be applied to
great advantage in the study of
these noncommutative spaces. \\

\noindent {\bf Example 1.}\\
 a) We start with a simple example from \cite{acb}. Let $X=\{a,
b\}$ be a set with two elements and define an equivalence relation on
$X$ that identifies $a$ and $b$, $a\sim b$:
$$ \overset{a}{\bullet} \leftrightsquigarrow \overset{b}{\bullet}$$
The corresponding groupoid here is the groupoid of pairs on the set $X$.
By Example 1 in Section 3.3  its groupoid algebra is the algebra of 2 by 2 matrices
$M_2(\mathbb{C})$. The identification is given by
$$f_{aa}(a, a)+ f_{ab}(a, b)+f_{ba}(b, a)+f_{bb}(b, b)\mapsto \left(
\begin{matrix}
f_{aa} & f_{ab} \\
f_{ba} & f_{bb}
\end{matrix}
\right).  $$

The algebra of functions on the classical quotient,
on the other hand, is given  by
$$\{f: X \to \mathbb{C}; \, f(a)=f(b)\}\simeq \mathbb{C}.$$
Thus the classical quotient and the noncommutative quotient are Morita
equivalent.\\

$$ M_2(\mathbb{C}) \xleftarrow{\text{noncommutative quotient}}\; \overset{a}{\bullet} \leftrightsquigarrow \overset{b}{\bullet} \; \xrightarrow{\text{classical quotient}} \mathbb{C}$$\\

\noindent b) The above  example can be generalized. For example let $X$ be
 a finite set with $n$ elements with the
equivalence relation $x\sim y$ for all $x, y$ in $X$. The
corresponding groupoid $\mathcal{G}$ is the groupoid of pairs
 and its groupoid algebra, representing the noncommutative quotient,  is
$$  \mathbb{C}\mathcal{G} \simeq M_n(\mathbb{C}).$$

The algebra of functions on the classical quotient  is given by
$$\{f: X \to \mathbb{C}; \, f(a)=f(b)\; \mbox{for all a, b in $X$}\}\simeq \mathbb{C}.$$
Again  the classical quotient is obviously Morita equivalent to the noncommutative quotient.\\

\noindent c) Let  $G$ be a group (not necessarily finite) acting on a finite set $X$.
The algebra of functions
on  the classical quotient is
$$ C(X/G)=\{f: X \to \mathbb{C}; \, f(x)=f(gx) \, \mbox{for all}\, g \in G, x\in X\}
\simeq \bigoplus_{\mathcal{O}} \mathbb{C},$$
where $\mathcal{O}$ denotes the set of orbits of $X$ under the action of $G$.

The noncommutative quotient, on the other hand,  is defined to be
the groupoid algebra of the transformation groupoid $ \mathcal{G}=X\rtimes G$. Note that as we saw before
this algebra is isomorphic to the crossed product algebra  $C(X)\rtimes
G$.  From Section 3.3 we have,
$$ \mathbb{C}\mathcal{G}\simeq C(X)\rtimes G \simeq \bigoplus_{ i\in \mathcal{O}}G_i \otimes
M_{n_i}(\mathbb{C}),$$
where $G_i$ is the isotropy group of the i-th orbit, and $n_i$  is the
size of the i-th orbit. Comparing the classical quotient with the
noncommutative quotient we see that:\\
i) If the action of $G$ is free then $G_i=\{1\}$  for all orbits
$i$ and therefore the  two algebras are Morita equivalent:
$$ C(X/G)\simeq \bigoplus_{\mathcal{O}} \mathbb{C}\overset{M. E.}{\sim}
\bigoplus_{ i\in \mathcal{O}} M_{n_i}(\mathbb{C})\simeq \mathbb{C}\mathcal{G}.$$

\noindent ii) The information about the isotropy groups is not lost in  the  noncommutative quotient
construction, while  the classical quotient totally neglects the
isotropy groups. \\

\noindent d) Let $X$ be a locally compact Hausdorff space and
consider the equivalence relation $x\sim y$ for all $x$ and $y$ in
$X$. The corresponding groupoid is again the groupoid of pairs. It
is a locally compact topological groupoid and its groupoid
$C^*$-algebra as we saw in Section 3.3 is  the algebra of compact
operators $\mathcal{K} (L^2(X, \mu))$. This algebra is obviously
strongly Morita
equivalent to the classical quotient algebra $\mathbb{C}$.\\

\noindent {\bf Example 2.} Let $\theta \in \mathbb{R}$ be a fixed real number.
Consider the action of the group of integers $G=\mathbb{Z}$ on the unit circle
$\mathbb{T}=\{z\in \mathbb{C}; \; |z|=1\}$ via  {\it rotation by $\theta$}:
$$ (n, z)\mapsto  e^{2\pi i n \theta}z.$$

For $\theta =\frac{p}{q}$  a rational number, the quotient space
$\mathbb{T}/\mathbb{Z}$ is a circle and the classical quotient algebra
$$``C(\mathbb{T}/\mathbb{Z})":=\{f \in C(\mathbb{T}); \;  f(gz)=f(z), \mbox{ for all}\; g, z\}
\simeq C(\mathbb{T}).$$
The noncommutative quotient $A_{\theta}=C(\mathbb{T})\rtimes_r \mathbb{Z}$, for any $\theta$, is
the unital $C^*$-algebra generated by two unitaries $U$ and $V$ subject to the relation
$$VU=e^{2\pi i \theta}UV.$$
It can be
shown that when $\theta$ is a rational number $A_{\theta}$  is isomorphic to
the space of continuous sections of the
endomorphism bundle $End (E)$ of a complex vector bundle $E$ over the 2-torus $\mathbb{T}^2$. Thus
for $\theta$ rational, $A_{\theta}$  is Morita equivalent to $C(\mathbb{T}^2)$:
$$C(\mathbb{T})\rtimes_r \mathbb{Z} \, \overset{M.E.}{\sim} \, C(\mathbb{T}^2).$$

If $\theta$ is an irrational number then each orbit is dense in
$\mathbb{T}$ and the quotient space $\mathbb{T}/\mathbb{Z}$ has
only two open set. It is an uncountable set with a trivial
topology. In particular it is not Hausdorff. Obviously, a
continuous function on the circle which is constant on each orbit
is necessarily constant since orbits are dense. Therefore
$$``C(\mathbb{T}/\mathbb{Z})" \simeq \mathbb{C}.$$
The noncommutative quotient $A_{\theta}$ in this case is  a simple $C^*$-algebra and is
not Morita equivalent to
$C(\mathbb{T}^2)$.

Let $G$ be a discrete group acting by homeomorphisms on a locally
compact Hausdorff space $X$. Recall that the action is called {\it free} if for all $g\neq e,$
we have $gx\neq x$ for
all $x\in G$.  The action is called {\it proper} if the map $G\times X \to X, (g, x)\mapsto gx$
is a proper map in the sense that the inverse image of a compact set is
compact. One shows that when the action is free and proper the orbit space $X/G$ of a locally compact and
Hausdorff space is  again a locally compact and Hausdorff space.
Similarly, if $X$ is a smooth manifold and the action is free and proper then there exists a unique smooth structure on $X/G$ such
that the quotient map $X\to X/G$ is smooth.

The following result of M. Rieffel \cite{rie} clarifies the relation between the classical
quotients and noncommutative quotients for group actions:

\begin{theorem} Assume $G$ acts  freely and  properly
on a locally compact Hausdorff space $X$. Then we have a strong Morita
equivalence between the $C^*$-algebras
$C_0(X/G)$ and $C_0(X)\rtimes_r G.$
\end{theorem}

\section{Cyclic cohomology}
Cyclic cohomology was discovered by Alain Connes in 1981
\cite{ac81}. One of Connes'  main motivations came from index
theory on foliated spaces. The $K$-theoretic index of a
transversally  elliptic operator on a foliated manifold is an
element of the $K$-theory group of a noncommutative algebra,
called the foliation algebra of the given foliated manifold.
Connes realized that to identify this class it would be desirable
to have a noncommutative analogue of the Chern character with
values in a, as yet unknown, cohomology theory for noncommutative
algebras that would play the role of de Rham homology of smooth
manifolds.

   Now to define a noncommutative de Rham theory for
noncommutative algebras is a highly nontrivial matter. This is in sharp
contrast with the situation in $K$-theory where extending the
topological $K$-theory to Banach algebras  is essentially a routine
matter. Note that  the usual algebraic
formulation of de Rham theory starts with the module of Kaehler
differentials  and its exterior algebra which  does not make sense for
noncommutative algebras.

Instead the answer was found by Connes by
analyzing the algebraic structures hidden in {\it traces of products of
commutators}. These expressions are directly defined in terms of an
elliptic operator and its parametix and were shown, via an index formula, to give the index of the
operator when paired with a $K$-theory class.

Let us read what Connes wrote in the Oberwolfach conference notebook after his
 talk, summarizing his discovery and how he arrived at it \cite{ac81}:\\

{\it ``The transverse elliptic theory for foliations requires as a
preliminary step a purely algebraic work, of computing for a
noncommutative algebra $\mathcal{A}$ the cohomology of the following
complex: $n$-cochains are multilinear functions $\varphi (f^0, \cdots,
f^n)$ of $f^0, \cdots ,f^n \in  \mathcal{A}$ where
$$\varphi (f^1, \cdots,
f^0)=(-1)^n \varphi (f^0, \cdots,
f^n)$$
and the boundary is
\begin{eqnarray*}
b\varphi (f^0, \cdots,f^{n+1}) &=&\varphi (f^0f^1, \cdots,f^{n+1})-\varphi (f^0, f^1f^2,\cdots,
f^{n+1})+\cdots \\& &+(-1)^{n+1} \varphi (f^{n+1}f^0, \cdots,
f^n).
\end{eqnarray*}
The basic class associated to a transversally elliptic operator,
for $\mathcal{A}=$ the algebra of the foliation is given by:
$$ \varphi (f^0, \cdots,f^n)=Trace  (\varepsilon F[F, f^0][F, f^1]\cdots[F, f^n]), \quad
f^i\in \mathcal{A}$$
where
$$F=\left(
\begin{matrix}
0& Q \\
P & 0
\end{matrix}
\right), \quad  \varepsilon= \left(
\begin{matrix}
1& 0 \\
0& -1
\end{matrix}
\right),
$$
and $Q$ is a parametrix of $P$. An operation
$$S: H^n(\mathcal{A})\to
H^{n+2}(\mathcal{A})$$
is constructed as well as a pairing
$$K(\mathcal{A})
\times H(\mathcal{A})\to \mathbb{C}$$ where $K(\mathcal{A})$ is
the algebraic $K$-theory of $A$. It gives the index of the
operator from its associated class $\varphi$. Moreover $<e,
\varphi>=<e, S \varphi>$ so that the important group to determine
is the inductive limit $H_g=\underset{\to}{Lim} H^n(\mathcal{A})$
for the map $S$. Using the tools of homological algebra the groups
$H^n(\mathcal{A}, \mathcal{A}^*)$ of Hochschild cohomology with
coefficients in the bimodule $\mathcal{A}^*$ are easier to
determine and the solution of the problem is obtained
in two steps,\\
1) the construction of a map
$$B: H^n(\mathcal{A},
\mathcal{A}^*) \to H^{n-1}(\mathcal{A})$$
and the proof of a long exact
sequence
$$
\cdots\to H^n(\mathcal{A}, \mathcal{A}^*) \overset{B}{\to} H^{n-1}(\mathcal{A})\overset{S}{\to}
H^{n+1}(\mathcal{A})\overset{I}{\to}H^{n+1}(\mathcal{A},
\mathcal{A}^*)\to \cdots
$$
where $I$ is the obvious map from the cohomology of the above
complex to
the Hochschild cohomology.\\
2) The construction of a spectral sequence with $E_2$ term given
by the cohomology of the degree $-1$ differential $I\circ B$ on
the Hochschild groups $H^n(\mathcal{A}, \mathcal{A}^*)$ and which
converges strongly to a graded group associated to the inductive
limit.

This purely algebraic theory is then used. For
$\mathcal{A}=C^{\infty}(V)$ one gets the de Rham homology of
currents, and for the pseudo torus, i.e. the algebra of the
Kronecker foliation, one finds that the Hochschild cohomology
depends on the Diophantine nature of the rotation number while the
above theory gives $H^0_g$ of dimension $1$, $H^1_g$ of dimension
$2$, and $H^0_g$ of dimension $1$ as expected but from some
remarkable
cancellations". }\\

In a different direction, cyclic homology also appeared in the
1983 work of Tsygan \cite{tsy} and was used also, independently,
by Loday and Quillen \cite{lq}. The Loday-Quillen-Tsygan theorem
states that the  cyclic homology of an algebra $A$ is the
primitive part (in the sense of Hopf algebras) of the Lie algebra
homology of the Lie algebra $gl(A)$ of stable matrices.
Equivalently, the Lie algebra homology of $gl(A)$ is isomorphic
with the exterior algebra over the cyclic homology of $A$ with
dimension shifted by 1:
$$H^{Lie}_{\bullet}(gl(A))\simeq \wedge (HC_{\bullet}(A)[-1]).$$
We won't pursue this connection in these notes.

In section 4.1 we recall basic notions of Hochschild (co)homology
theory and give several computations. Theorems of Connes
\cite{ac85} (resp. Hochschild-Kostant-Rosenberg \cite{hkr}) on the
Hochschild cohomology of the algebra of smooth functions on a
manifold (resp. algebra of regular functions on  a smooth affine
variety) are among the most important results of this theory.

In Section 4.2 we define cyclic cohomology via Connes' cyclic complex. The easiest approach, perhaps,
to introduce the map $B$  is to introduce first a
  bicomplex called
cyclic bicomplex. This leads to a new definition of cyclic (co)homology,
 a definition of the operator $B$, and a proof
of the long exact sequence of Connes, relating Hochschild, and cyclic
cohomology groups. A third definition of cyclic cohomology is via
Connes's $(b, B)$-bicomplex. The equivalence of these three
definitions is established by explicit maps. Finally we recall Connes'
computation of the cyclic (co)homology of the algebra of smooth
functions on a manifold, and Burghelea's result on the cyclic homology of group rings.

\subsection{Hochschild (co)homology}
Hochschild cohomology of associative algebras was defined by G.
Hochschild through an explicit complex in \cite{hoc}. This complex
is a generalization of the standard complex for group cohomology.
One of the original motivations was to give a cohomological
criterion for   separability of algebras as well as a
classification of (simple types) of algebra extensions in terms of
second  cohomology.  Once it was realized, by Cartan and Eilenberg
\cite{ce}, that Hochschild cohomology is an example of their newly
discovered theory of derived functors, tools of homological
algebra like resolutions became available.

The Hochschild-Kostant-Rosenberg theorem \cite{hkr} and its smooth
version by Connes\cite{ac85}  identifies the Hochschild homology
of the algebra of regular functions on a smooth affine variety or
smooth functions on a  manifold with differential forms and   is
among the most important results of this theory. Because of this
result one usually thinks of the Hochschild homology of an algebra
$A$ with coefficients in $A$ as a noncommutative analogue of
differential forms on $A$.

As we shall see later in this section Hochschild (co)homology is
related to cyclic (co)homology through Connes' long exact
sequence. For this reason computing the Hochschild (co)homology is
often the first step in computing the cyclic (co)homology of a
given algebra.

Let $A$ be an algebra and let $M$ be an $A$-bimodule. Thus $M$ is a left and
right $A$-module and the two actions are compatible in the sense that
$a(mb)=(am)b,$
for all $a,b$ in $A$ and $m$ in $M$. The {\it Hochschild cochain complex of $A$
with coefficients in $M$}, denoted
$( C^{\bullet}(A, M), \delta )$, is defined as
$$ C^0(A, M)=M, \quad C^n(A, M)=Hom(A^{\otimes n}, M), \quad n\geq 1,$$
\begin{eqnarray*}
(\delta m)(a)&=&ma-am,\\
(\delta f)(a_1, \cdots, a_{n+1})&=& a_1f(a_2, \cdots, a_{n+1})+\sum_{i=1}^n (-1)^{i+1}
f(a_1, \cdots, a_ia_{i+1},\cdots , a_{n+1})\\
& & + (-1)^{n+1}f(a_1, \cdots, a_{n})a_{n+1},
\end{eqnarray*}
where $m \in M=C^0(A, M)$, and $f \in C^n(A,M), n\geq 1$.

One  checks that
$$\delta^2 =0.$$
The cohomology of the complex $(C^{\bullet}(A,M), \delta )$ is by
definition the {\it Hochschild cohomology} of $A$ with
coefficients in $M$ and will be denoted by $H^{\bullet} (A, M).$

Among all  bimodules $M$ over an  algebra $A$, the following  two bimodules play
an important role:\\

\noindent 1) M=A, with bimodule structure $a(b)c=abc$, for all $a, b, c$ in $A$.  The Hochschild complex
$C^{\bullet} (A, A)$ is also known as the {\it deformation complex }, or {\it Gerstenhaber complex}
 of
$A$. It plays an important role in deformation theory of associative
algebras  pioneered
by Gerstenhaber \cite{ger}. For example
it is easy to see that $H^2(A, A)$ is the space of {\it infinitesimal deformations} of $A$ and
$H^3(A, A)$
is the {\it space of obstructions} for deformations of $A$.\\

\noindent 2) $M= A^*=Hom (A, k)$ with bimodule structure defined by
$$(afb)(c)=f(bca),$$
for all $a, b, c$ in $A$, and $f$ in $A^*$. This
bimodule is relevant to cyclic cohomology. Indeed as we shall see
the Hochschild groups $H^{\bullet}(A, A^*)$ and the cyclic cohomology groups
$HC^{\bullet}(A)$ enter into a long
exact sequence (Connes's long sequence). Using the identification
$$Hom(A^{\otimes n}, A^*) \simeq Hom(A^{\otimes (n+1)}, k), \quad f\mapsto \varphi,$$
$$\varphi(a_0, a_1, \cdots ,a_n)=f(a_1, \cdots, a_n)(a_0),$$
the Hochschild differential $\delta$ is transformed into the
differential $b$ given by
\begin{eqnarray*}
b\varphi(a_0, \cdots ,a_{n+1})&=&\sum_{i=0}^n (-1)^i
 \varphi (a_0, \cdots a_ia_{i+1}, \cdots, a_{n+1})\\
 & &+(-1)^{n+1}\varphi (a_{n+1} a_0, a_1, \cdots ,a_n).
\end{eqnarray*}
Thus for  $n=0, 1, 2$ we have the following formulas for $b$:
\begin{eqnarray*}
b\varphi (a_0, a_1)&=&\varphi (a_0a_1)-\varphi (a_1a_0),\\
b\varphi (a_0, a_1, a_2)&=&\varphi (a_0a_1, a_2)-\varphi (a_0, a_1a_2)+
\varphi (a_2a_0, a_1),\\
b\varphi (a_0, a_1, a_2, a_3)&=&\varphi (a_0a_1, a_2, a_3) -\varphi
(a_0, a_1a_2, a_3)\\
& &+\varphi (a_0, a_1, a_2a_3) -\varphi (a_3a_0, a_1, a_2).
\end{eqnarray*}
We give a few examples of Hochschild cohomology in low dimensions. \\

\noindent{\bf Examples}\\
1. $n=0$. It is clear that
$$H^0(A, M)=\{m\in M; ma=am \quad \text{for all} \quad a\in A\}.$$
In particular for $M=A^*$,
$$H^0(A, A^*)=\{f: A \to k; \; f(ab)=f(ba)\; \mbox{ for all}\; a, b \in
A\},$$
is the space of traces on $A$. \\

\noindent Exercise: For $A=k[x, \frac{d}{dx}]$, the algebra of
differential operators with
polynomial coefficients, show that $H^0(A, A^*)=0.$\\

\noindent 2. $n=1$. A Hochschild {\it 1-cocycle} $f \in C^1(A, M)$ is  simply a {\it derivation}, i.e. a
linear map  $f:A \rightarrow M$
such that
$$f(ab)=af(b) +f(a)b,$$
for all $a, b$ in $A$. A cocycle is a {\it coboundary} if and only
if the corresponding derivation is  {\it inner}, that is there
exists $m$ in $M$ such that $f(a)=ma-am$ for all $a$ in $A$.
Therefore
$$H^1(A, M)=\frac{\text{derivations}}{\text{inner derivations}}$$
Sometimes this is called the space of {\it outer derivations} of $A$
to $M$. \\

\noindent Exercise: 1) Show that any derivation on the algebra $C(X)$ of
continuous functions on a compact Hausdorff space $X$ is zero. (Hint: If
$f=g^2$ and $g(x)=0$ then $f'(x)=0$.)\\
2) Show that any derivation on the matrix algebra $M_n(k)$ is inner.
(This was proved by Dirac in \cite{dir} where derivations are called
{\it quantum differentials}. \\
3) Show that any derivation on the {\it Weyl algebra} $A=k[x, \frac{d}{dx}]$ is inner as well. \\

\noindent 3. $n=2.$ We show, following Hochschild \cite{hoc}, that  $H^2(A, M)$ classifies
{\it abelian  extensions} of $A$ by $M$. Let $A$ be a unital algebra over a field
$k$. By definition, an abelian
extension is an exact sequence of algebras
$$0\longrightarrow M\longrightarrow B \longrightarrow A\longrightarrow 0,$$
such that $B$ is unital, $M$ has trivial multiplication ($M^2=0$),
and the induced $A$-bimodule structure on $M$ coincides with the
original bimodule structure. Let $E(A,M)$ denote the set of
isomorphism classes of such extensions. We define a natural
bijection
$$E(A, M) \simeq H^2(A, M)$$
as follows. Given an extension as above, let $s: A \rightarrow B$
be a linear splitting for the projection $B\rightarrow A$, and let
$f:A \otimes A \rightarrow M$ be its {\it curvature} defined by,
$$ f(a, b)=s(ab)-s(a)s(b),$$
for all $a, b$ in $A$. One can easily check that $f$ is a
Hochschild 2-cocycle and its class is independent of the choice of
splitting $s$. In the other direction, given a 2-cochain $f:A
\otimes A \rightarrow M$, we try to define a multiplication on
$B=A\oplus M$ via
$$(a,m)(a',m')= (aa', am'+ma'+f(a, a')).$$
It can be checked that this defines an associative multiplication if and
only if $f$ is a 2-cocycle. The extension associated to a 2-cocycle $f$
is the extension
$$ 0\longrightarrow M \longrightarrow A\oplus M \longrightarrow
A\longrightarrow 0.$$
It can be checked that these two maps are bijective and inverse to each
other.

We show that Hochschild cohomology is a derived functor. Let
$A^{op}$ denote the {\it opposite algebra} of $A$, where
$A^{op}=A$ and the new multiplication is defined by $a.b:=ba$.
 There is a one to one
correspondence between $A$-bimodules and left $A\otimes A^{op}$-modules defined by
$$(a\otimes b^{op})m=amb.$$
Define a functor from the category of
left $A\otimes A^{op}$  modules to $k$-modules by
$$ M\mapsto Hom_{A\otimes A^{op}}(A, M)=\{m\in M; ma=am \quad \text{for all} \quad a\in A\}=H^0(A, M).$$
To show that Hochschild cohomology is the derived functor of the functor
$ Hom_{A\otimes A^{op}}(A, -)$, we introduce the {\it bar resolution} of $A$. It is
defined by
$$ 0\longleftarrow A\overset{b'}{\longleftarrow} B_1(A)\overset{b'}{\longleftarrow} B_2(A) \cdots,$$
where $B_n(A)=A\otimes A^{op}\otimes A^{\otimes n}$ is the free left $A\otimes A^{op}$ module
generated by $A^{\otimes n}$. The differential $b'$ is defined by
\begin{eqnarray*}
b'(a\otimes b\otimes a_1\otimes \cdots \otimes a_n)&=&aa_1\otimes b\otimes a_2 \cdots \otimes a_n\\
& &+\sum_{i=1}^{n-1} (-1)^i(a\otimes b\otimes a_1\otimes \cdots
a_ia_{i+1}\cdots
\otimes a_n)\\
& &+(-1)^{n}(a\otimes a_nb\otimes a_1\otimes \cdots
\otimes a_{n-1}).
\end{eqnarray*}
Define the operators $s:B_n(A)\rightarrow B_{n+1}(A)$, $n\geq 0$,  by
$$s(a\otimes b\otimes a_1\otimes \cdots \otimes a_n)=1\otimes b\otimes a\otimes
a_1\otimes \cdots \otimes a_n.$$
One checks that
$$b's+sb'=id$$
which shows that $(B_{\bullet}(A), b')$ is acyclic. Thus
$(B_{\bullet}(A), b')$ is a projective resolution of $A$ as an $A$-bimodule. Now, for any $A$-bimodule
$M$  we have
$$ Hom_{A\otimes A^{op}} (B_{\bullet}(A), M)\simeq (C^{\bullet}(A, M),
\delta),$$
which shows that Hochschild cohomology is a derived functor.

One can therefore use resolutions to compute Hochschild cohomology
groups. Here are a few exercises\\
1. Let
$$A=T(V)=k\oplus V \oplus V^{\otimes 2}\oplus \cdots,$$
be the tensor algebra of a vector space $V$. Show that
$$ 0\longleftarrow T(V) \overset{\delta}{\longleftarrow} T(V)\otimes T(V)
 \overset{\delta}{\longleftarrow}
T(V)\otimes V \otimes T(V) \longleftarrow 0,$$
$$\delta (x\otimes y)=xy, \quad \delta (x\otimes v \otimes y)=xv \otimes
y -x\otimes vy,$$
is a free resolution of $T(V)$. Conclude that $A$ has Hochschild
cohomological dimension 1 in the sense that $H^n(A, M)=0$ for all $M$
and all $n\geq 2.$ Compute $H^0$ and $H^1$ \cite{lod}.\\

2. Let $A=k[x_1, \cdots ,x_n]$ be the polynomial algebra in $n$
variables over a field  $k$ of characteristic zero. Let $V$ be an $n$ dimensional
vector space over $k$. Define a resolution of the form
$$ 0\leftarrow A \leftarrow A\otimes A \leftarrow A\otimes V\otimes A
\leftarrow \cdots A\otimes \wedge^i V\otimes A \cdots \leftarrow A \otimes
\wedge^n V \otimes A \leftarrow 0$$
by tensoring resolutions in 1) above for one dimensional vector spaces.

Conclude that for any symmetric $A$-bimodule $M$,
$$ H^i(A, M) \simeq M\otimes \wedge^i V, \quad i=0, 1,\cdots. $$

Before proceeding further let us recall the definition of the {\it
Hochschild homology} of an algebra with coefficients in a bimodule
$M$. The {\it Hochschild complex  of $A$ with coefficients in
$M$},  $(C_{\bullet}(A, M), \delta)$, is defined by
$$C_0(A,M)=M, \quad \text{ and} \quad C_n(A,M)=M \otimes A^{\otimes n}, \; n=1, 2,\cdots$$
and the {\it Hochschild boundary} $\delta :C_n(A, M)\longrightarrow C_{n-1}(A,
M)$ is defined by
\begin{eqnarray*}
\delta (m\otimes a_1\otimes \cdots a_n)&=&m a_1\otimes a_1\cdots a_n
 + \sum_{i=1}^{n-1}(-1)^i m\otimes a_1\otimes a_ia_{i+1}\cdots a_n\\
& & +(-1)^n a_nm \otimes a_1\otimes \cdots a_n.
\end{eqnarray*}
The  Hochschild homology of $A$ with coefficients in $M$ is, by definition,
the homology of the complex $(C_{\bullet}(A, M), \delta)$.
We denote this homology by $H_{\bullet}(A, M)$. It is clear that
$$H_0(A, M)=M/[A, M],$$
where $[A, M]$ is the subspace of $M$ spanned by commutators $am-ma$ for
$a$ in $A$ and $m$ in $M$.

The following  facts are easily established:\\
1. Hochschild homology $H_{\bullet}(A, M)$ is the derived functor of the
functor
$$ A\otimes A^{op}-Mod \longrightarrow k-Mod, \quad M\mapsto
A\otimes_{A\otimes A^{op}} M,$$
i.e.
$$H_n(A, M)=Tor_n^{A\otimes A^{op}}(A, M).$$
For the proof one uses the bar resolution as before.\\

2. (Duality) Let $M^*=Hom(M, k)$. It is an $A$-bimodule via
$(afb)(m)=f(bma).$ One checks that  the natural isomorphism
$$Hom (A^{\otimes n}, M^*)\simeq Hom (M\otimes A^{\otimes n}, k),\quad  n=0, 1,\cdots$$
is compatible with differentials. Thus if $k$ is field of
characteristic zero, we have
$$ H^{\bullet} (A, M^*) \simeq (H_{\bullet}(A, M))^*.$$

>From now on we denote by $HH^n(A)$ the Hochschild group $H^n(A, A^*)$
and by $HH_n(A)$ the Hochschild group $H_n(A, A)$.

For applications of Hochschild and cyclic (co)homology to noncommutative geometry,
 it is crucial to consider
topological algebras, topological bimodules  and continuous chains and
cochains on them. For example while the algebraic Hochschild groups of
the algebra of smooth functions on a smooth manifold are not known,
its  topological Hochschild (co)homology
is computed by Connes  as we recall below. We will give only a brief
outline of the definitions and refer the reader to \cite{ac85, acb}
for more details.

Let $A$ be a locally convex topological
algebra and $M$ be a locally convex topological $A$-bimodule. Thus $A$
is a locally convex topological vector space and the multiplication map $A\times A \to A$ is continuous.
Similarly $M$ is a locally convex topological vector space such that both module maps $A\times M \to M$
and $M\times A \to M$ are continuous. 
In the
definition of continuous Hochschild homology one uses the {\it projective tensor
product}
$M\hat{\otimes}A^{\hat{\otimes}n}$
of locally convex spaces. The algebraic Hochschild boundary, being continuous, naturally extends to
 topological completions.

 For cohomology one should use {\it jointly
 continuous} multilinear maps
 $$\varphi : A\times \cdots \times A \to M.$$
 With these provisions, the rest of the algebraic formalism remains the same
 and carries over to the topological set up. In using projective resolutions, one
 should use only those topological resolutions that admit a continuous linear splitting.  This gaurantees
 that the comparison theorem for projective resolutions remain true in the continuous setting.

We give a few examples of Hochschild (co)homology computations. In
particular we shall see  that group homology and Lie algebra homology are
instances of Hochschild homology. We start by recalling the classical
results of Connes \cite{ac85}
and Hochschild-Kostant-Rosenberg \cite{hkr} on the Hochschild homology of smooth
commutative algebras.\\

\noindent{\bf Example} (Commutative Algebras)\\
 Let $A$ be a commutative unital algebra over a ring $k$. We recall the definition of the {\it algebraic
de Rham complex} of $A$. The module of 1-forms over $A$, denoted by $\Omega^1 A$, is defined to be a left
$A$-module $\Omega^1 A$ with a universal derivation
$$ d: A\longrightarrow \Omega^1 A.$$
This means that any other derivation $\delta : A \rightarrow M$
into a left $A$-module $M$, uniquely factorizes through $d$. One
usually defines $\Omega^1 A=I/I^2$ where $I$ is the kernel of the
multiplication map $A\otimes A \rightarrow A$. Note that since $A$
is commutative this map is an algebra homomorphism. $d$ is defined
by
$$d(a)= a\otimes 1-1\otimes
a \quad \text{mod} (I^2).$$
One  defines the space of n-forms on $A$ as the n-th exterior power of the
$A$-module
$\Omega^1 A$:
$$\Omega^n A: =\wedge^n_A \Omega^1A.$$
There is a unique extension of $d$ to a graded derivation
$$d:\Omega^{\bullet}A \longrightarrow \Omega^{\bullet +1}A.$$
It satisfies the relation $d^2=0$. The {\it algebraic de Rham
cohomology} of $A$ is defined to be
 the cohomology of the complex  $(\Omega^{\bullet}A, d)$.

 Let $M$ be a symmetric $A$-bimodule. We compare the complex $(M\otimes_A
 \Omega^{\bullet}A, 0)$ with the Hochschild complex of $A$ with
 coefficients in $M$. Consider the {\it antisymmetrization map}
 $$\varepsilon_n: M\otimes_A \Omega^nA \longrightarrow M\otimes
 A^{\otimes n}, \quad n=0, 1, 2, \cdots,$$
 $$\varepsilon_n (m\otimes da_1\wedge  \cdots \wedge da_n)=\sum_{\sigma \in S_n} sgn (\sigma)
 m\otimes a_{\sigma (1)}\otimes \cdots \otimes a_{\sigma (n)},$$
 where $S_n$ denotes the symmetric group on $n$ letters.
We also have a  map
 $$\mu_n :M\otimes A^{\otimes n} \longrightarrow M\otimes_A \Omega^n A, \quad n=0, 1, \cdots$$
 $$\mu_n (m\otimes a_1\otimes \cdots a_n)= m\otimes da_1 \wedge \cdots \wedge da_n.$$
 One checks that both maps are morphisms of complexes, i.e.
 $$\delta \circ \varepsilon_n =0,  \quad \mu_n \circ \delta =0.$$
 Moreover, one has
 $$\mu_n \circ \varepsilon_n =n! \, Id_n.$$
 It follows that if $k$ is a field of characteristic zero then the
 antisymmetrization
 map
 $$ \varepsilon_n: M\otimes_A \Omega^nA \longrightarrow H_n(A, M),$$
 is an inclusion. For $M=A$ we obtain a natural inclusion
 $$\Omega^nA \longrightarrow HH_n(A).$$

 The celebrated Hochschild-Kostant-Rosenberg  theorem \cite{hkr} states
 that if $A$ is the algebra of regular functions on an smooth affine variety
 the above map is an isomorphism.

Let $M$ be a smooth closed manifold and let $A=C^{\infty}(M)$ be the
algebra of smooth complex valued functions on $M$. It is a locally convex (in fact, Frechet)
topological  algebra. Fixing a finite atlas on $M$, one defines a family
of seminorms
$$p_n(f)=sup \{|\partial^{I}(f)|; \; |I|\leq n\},$$
 where the supremum is
 over all coordinate charts. It is easily seen that the induced topology
 is independent of the choice of atlas. In \cite{ac85}, using an
 explicit resolution, Connes shows that the canonical map
 $$HH_n^{cont}(A) \to \Omega^n M, \quad f_0\otimes \cdots \otimes f_n \mapsto
 f_0df_1 \cdots df_n,$$
 is an isomorphism. In fact the original, equivalent,  formulation of Connes in \cite{ac85} is for
   continuous Hochschild cohomology $HH^n (A)$ which is shown to be isomorphic to
    the continuous dual of $\Omega^n M$ (space of n-dimensional {\it de Rham
    currents}). \\

\noindent{\bf Example} (Group Algebras)\\
 It is clear from the original definitions that group (co)homology is an
example of Hochschild (co)homology. Let $G$ be a group and  let $M$ be a left $G$-module over the
ground ring $k$. Recall that the standard complex for computing group
cohomology \cite{lod} is the complex
$$M \overset{\delta}{\longrightarrow} C^1(G, M)\overset{\delta}{\longrightarrow} C^2(G, M)
\overset{\delta}{\longrightarrow}
\cdots,$$
where
$$C^n(G, M)=\{f: G^n \longrightarrow M \}.$$
The differential $\delta $ is defined by
$$(\delta m )(g)=gm-m,$$
\begin{eqnarray*}
\delta f(g_1, \cdots ,g_{n+1})&=&g_1f(g_2, \cdots, g_{n+1})
+\sum_{i=1}^n (-1)^if(g_1, \cdots g_ig_{i+1}, \cdots g_{n+1})\\
& & +(-1)^{n+1} f(g_1, g_2, \cdots, g_n).
\end{eqnarray*}

Let $A=kG$ denote the group algebra of the group $G$ over $k$. Then $ M$ is a
$kG$-bimodule via the two
actions
$$g.m =g(m), \quad m.g=m,$$
for all $g$ in $G$ and $m$ in $M$. It is clear that for all $n$,
$$ C^n(kG, M)\simeq C^n(G, M),$$
and the two differentials are the same. It follows that the cohomology of $G$ with
coefficients in $M$ coincides with the Hochschild cohomology of $kG$
with  coefficients in $M$.\\

\noindent{\bf Example} (Enveloping Algebras). \\
We show that Lie algebra
(co)homology is an example of Hochschild (co)homology. Let $\mathfrak{g}$ be a Lie algebra
and $M$  be a $\mathfrak{g}$-module. This simply means that we have a Lie algebra
morphism
$$\mathfrak{g}  \longrightarrow End_k(M).$$
The {\it Lie algebra homology} of $\mathfrak{g}$  with
coefficients in $M$ is the homology of the {\it
Chevalley-Eilenberg complex} defined by
$$M \longleftarrow M\otimes \wedge^1 \mathfrak{g}  \longleftarrow M\otimes \wedge^2 \mathfrak{g} \longleftarrow M \cdots,$$
where the differential is defined by
$$\delta (m\otimes X)= X(m),$$
\begin{eqnarray*}
\delta (m\otimes X_1 \wedge X_2\wedge \cdots \wedge X_n)&= & \sum_{i<j}
(-1)^{i+j}
m\otimes [X_i, X_j]\wedge X_1  \cdots \wedge \hat{X_i}\cdots \hat{X_j}\cdots \wedge X_n\\
&+ & \sum_i (-1)^iX_i(m) \otimes X_1 \wedge \cdots \hat{X_i}\wedge \cdots
\wedge X_n.
\end{eqnarray*}
One checks that $\delta^2 =0.$

Let $U(\mathfrak{g})$ denote the enveloping algebra of $\mathfrak{g}$. Given a $\mathfrak{g}$ module $M$ we
define a $U(\mathfrak{g})$-bimodule $M'=M$ with left and right $U(\mathfrak{g})$-actions
defined by
$$ X\cdot m= X(m), \quad m\cdot X=0.$$

Define a map
$$\varepsilon_n :  C_n^{Lie}(\mathfrak{g}, M) \longrightarrow C_n(U(\mathfrak{g}), M'),$$
$$\varepsilon_n (m\otimes X_1\wedge  \cdots \wedge X_n)=\sum_{\sigma \in S_n} sgn (\sigma)
 m\otimes X_{\sigma (1)}\otimes \cdots \otimes X_{\sigma (n)}.$$
 One checks that $\varepsilon$ is  a chain map (prove this!). We claim that it
 is a quasi-isomorphism, i.e. it induces isomorphism on homology. To prove this, we define a filtration on
 $(C_{\bullet}(U(\mathfrak{g}), M), \delta)$ using the Poincare-Birkhoff-Witt
 filtration on $U(\mathfrak{g})$. The associated $E^1$ term is the de Rham complex of the symmetric algebra
 $S(\mathfrak{g})$. The induced map is the antisymmetrization map
 $$\varepsilon_n : M\otimes \wedge^n \mathfrak{g} \to M\otimes S(\mathfrak{g})^{\otimes n}.$$
By Hochschild-Kostant-Rosenberg's theorem, this map is a quasi-isomorphism
hence the original map is a quasi-isomorphism. \\

\noindent{\bf Example} (Morita Invariance of Hochschild (Co)Homology)\\
Let $A$ and $B$ be unital Morita equivalent $k$-algebras. Let $X$ be
an equivalence $A-B$ bimodule and $Y$ its inverse bimodule. Let $M$ be
an $A-A$ bimodule and $N=Y\otimes_A M\otimes_A X$ the corresponding
$B$-bimodule. Morita invariance of Hochschild homology states that there
is a natural isomorphism
$$H_n(A, M) \simeq H_n(B,N),$$
for all $n\geq 0$ \cite{lod}. We sketch a proof of this result for the
special case where $B=M_k(A)$ is the algebra of $k$ by $k$ matrices over
$A$.

Let $M$ be an $A$-bimodule,
and let $M_k(M)$ be the space of $k$ by $k$ matrices with coefficients in
$M$. It is a bimodule over  $M_k(A)$. The
 {\it generalized trace map} is defined by
$$Tr: C_n(M_k(A), M_k(M))\longrightarrow C_n(A, M),$$
\begin{eqnarray*}
Tr( \alpha_0 \otimes m_0 \otimes \alpha_1\otimes a_1 \otimes \cdots
\otimes \alpha_n \otimes a_n )
&=&\\
tr (\alpha_0 \otimes \alpha_1 \otimes \cdots \alpha_n) m_0\otimes a_1
\otimes \cdots a_n, & &
\end{eqnarray*}
where $\alpha_i \in M_k(k), a_i\in A, m_0 \in M,$ and $tr:M_k(k)
\longrightarrow k$ is the standard trace of matrices. \\

\noindent Exercise:\\
 1. Show that $Tr$ is a chain map.\\
2. Let $i: A \rightarrow M_k(A)$ be the map that sends $a$ in $A$ to the
matrix with only one non-zero component in the upper left corner equal
to $a$. There is a similar map $M \rightarrow M_k(M)$. Define a map
$$i_*: C_n(A, M) \longrightarrow C_n(M_k(A), M_k(M)),$$
$$i_* (m\otimes a_1 \otimes \cdots \otimes a_n)=i(m)\otimes i(a_1)
\otimes \cdots \otimes i(a_n).$$
Show that
$$Tr\circ i_* =id.$$

It is however not true that $i_* \circ Tr =id.$ There is instead a
homotopy between $i_* \circ Tr$ and  $id$. The homotopy is given in \cite{lod} and
we won't reproduce it here.

As a special case of the Morita invariance theorem, we have an
isomorphism of Hochschild homology groups
$$ HH_n(A)=HH_n(M_k(A)),$$
for all $n$ and $k$.

We need to know, for example when defining the noncommutative
Chern character map, that inner automorphisms act by identity on
Hochschild homology and inner derivations act by zero. Let $A$ be
an algebra, let $u\in A$ be an invertible element and let $a\in A$
be any element. They induce the chain maps
$$ \Theta: C_n(A) \to C_n(A) \quad a_0\otimes \cdots \otimes a_n \mapsto
ua_0u^{-1}\otimes \cdots \otimes ua_nu^{-1},$$
$$L_a :C_n(A) \to C_n(A) \quad a_0\otimes \cdots \otimes a_n \mapsto
\sum_{i=0}^na_0\otimes \cdots \otimes [a, a_i]\otimes \cdots \otimes a_n.$$

\begin{lem} $\Theta$ induces the identity map on Hochschild homology and $L_a$ induces the
zero map.
\end{lem}
\begin{proof}
The maps \cite{lod}, $h_i: A^{\otimes n+1} \to A^{\otimes n+2}, \; i=0,
\cdots, n$
$$h_i(a_0\otimes \cdots \otimes a_n)=(a_0u^{-1}\otimes ua_1u^{-1},
\cdots, u\otimes a_{i+1}\cdots \otimes a_n)$$
define a homotopy
$$h=\sum_{i=0}^n (-1)^i h_i$$
betweem $id$ and $\Theta$.

For the second part one  checks again  that the maps $h_i: A^{\otimes n+1} \to A^{\otimes n+2}, \; i=0,
\cdots, n$,
$$h_i(a_0\otimes \cdots \otimes a_n)=(a_0\otimes \cdots a_i\otimes a\cdots \otimes a_n),$$
define a homotopy between $L_a$ and $0$ \cite{lod}. 

\end{proof}

\subsection{Cyclic (co)homology}

\subsubsection{Connes' cyclic complex}
Cyclic cohomology was first defined by Connes \cite{ac81, ac85} through  a remarkable
subcomplex of the Hochschild complex called the {\it cyclic complex}.
Let $k$ be a field of characteristic zero and let 
$(C^{\bullet}(A), b)$ denote the Hochschild complex of a $k$-algebra $A$ with
coefficients in the $A$-bimodule $A^*$. We have
$$C^n(A)=Hom (A^{\otimes (n+1)}, k), \quad n=0, 1, \cdots,$$
\begin{eqnarray*}
(b f)(a_0, \cdots ,a_{n+1})&=&\sum_{i=0}^n(-1)^if(a_0, \cdots,a_ia_{i+1},\cdots ,a_{n+1})\\
& &+(-1)^{n+1}f(a_{n+1}a_0, \cdots,\cdots ,a_{n}).
\end{eqnarray*}

An $n$-cochain $f \in C^n(A)$ is called {\it cyclic} if
$$f(a_n, a_0, \cdots ,a_{n-1})=(-1)^n f(a_0, a_1, \cdots, a_n)$$
for all $a_0, \cdots ,a_n$ in $A$. We denote the space of cyclic
cochains on $A$ by $C^n_{\lambda}(A).$

\begin{lem} The space of cyclic cochains is invariant under the action of $b$, i.e. for all $n$,
$$b \,C^n_{\lambda} (A) \subset C^{n+1}_{\lambda} (A).$$
\end{lem}
\begin{proof} Define  the operators $\lambda :C^n(A) \rightarrow
C^n(A)$ and $b': C^n(A) \rightarrow C^{n+1}(A)$ by
\begin{eqnarray*}
(\lambda f)(a_0, \cdots ,a_n)&=&(-1)^n f(a_n, a_0, \cdots ,a_{n-1}),\\
 (b' f)(a_0, \cdots ,a_{n+1})&=&\sum_{i=0}^n(-1)^if(a_0, \cdots,a_ia_{i+1},\cdots ,a_{n+1}).
\end{eqnarray*}

One checks that
$$(1-\lambda)b=b'(1-\lambda).$$
Since
$$ C^n_{\lambda} (A)=Ker (1-\lambda),$$
the lemma is proved.
\end{proof}

We therefore have a subcomplex of the Hochschild complex, called the
{\it cyclic complex} of $A$:
$$ C_{\lambda}^0 (A)\overset{b}{\longrightarrow} C_{\lambda}^1 (A)
\overset{b}{\longrightarrow}C^2_{\lambda}(A) \overset{b}{\longrightarrow}\cdots .$$

The cohomology of this complex is called the {\it cyclic cohomology} of
$A$ and will be denoted by $HC^n (A)$, $n=0, 1, 2,\cdots.$
A cocycle for cyclic cohomology is called a {\it cyclic cocycle}. It
satisfies the two conditions:
$$ (1-\lambda )f=0, \quad {\text and} \quad bf=0.$$

The inclusion of complexes
$$ (C^{\bullet}_{\lambda}(A), b) \longrightarrow (C^{\bullet}(A), b),$$
induces a map $I$ from the cyclic cohomology of $A$ to the
Hochschild cohomology of $A$ with coefficients in the $A$-bimodule
$A^*$:
$$ I: HC^n(A) \longrightarrow HH^n(A), \quad n=0, 1, 2, \cdots.$$
We shall see that this map is part of a long exact sequence relating
Hochschild and cyclic cohomology, called Connes' long exact sequence.
For the moment we mention that $I$ need
not be injective (see example below).\\

\noindent{\bf Examples}\\
1. Let $ A=k$ be a field of characteristic zero. We have
$$C_{\lambda}^{2n}(k)\simeq k, \quad C_{\lambda}^{2n+1}(k)=0.$$
The cyclic complex reduces to
$$ 0\longrightarrow k\longrightarrow 0\longrightarrow k \longrightarrow \cdots.$$

It follows that for all $n\geq 0,$
$$HC^{2n}(k)=k, \quad HC^{2n+1}(k)=0.$$
Since $HH^n(k)=0$ for $n\geq 1$, we conclude that the map $I$ need not
be injective and the cyclic complex is not a retraction of the
Hochschild complex.\\

\noindent 2. It is clear that $HC^0(A)=HH^0(A)$ is the space of traces
on $A$. \\

\noindent 3. Let $A=C^{\infty}(M)$ be the algebra of smooth
functions on a closed smooth manifolds $M$ of dimension $n$. One
checks that
$$\varphi (f_0, f_1, \cdots , f_n)=\int_M f_0df_1 \cdots df_n,$$
is a cyclic $n$-cocycle on $A$. In fact the Hochschild cocycle
property $b\varphi =0$ is a consequence of the graded
commutativity of the algebra of differential forms and the cyclic
property $(1-\lambda)\varphi =0$ follows from Stokes formula.

This example can be generalized in several directions. For
example, Let $V$ be an $m$-dimensional closed singular chain (a
cycle) on $M$, e.g. $V$ can be a closed $m$-dimensional submanifold of $M$.  Then integration on $V$ defines
 an $m$-dimensional
cyclic cocycle on $A$:
$$ \varphi (f_0, f_1, \cdots , f_m)=\int_V f_0df_1 \cdots df_m.$$
We obtain a map
$$ H_m(M, \mathbb{C}) \longrightarrow HC^m(A), \quad m=0, 1, \cdots, $$
from singular homology (or its equivalents) to cyclic cohomology.

More generally, let  $C$ be a {\it closed $m$-dimensional de Rham current} on $M$. Thus $C: \Omega^m M \to
\mathbb{C}$ is a continuous linear functional on  $\Omega^m M$ such that
$dC (\omega):= C(d\omega)=0$ for all $\omega \in \Omega^{m-1} M.$ Then
one checks that the cochain $\varphi$ defined by
$$\varphi (f_0, f_1, \cdots , f_m)=< C, f_0df_1 \cdots df_m>,$$
is a cyclic cocycle on $A$. 

A noncommutative generalization of this procedure involves the notion of
a {\it cycle on an algebra} due to Connes \cite{ac85} that we recall now. It gives a  geometric
and intuitively appealing presentation
for cyclic cocycles. It also leads to a definition of cup product in
cyclic cohomology and the $S$ operator.
Let
$$\Omega =\Omega^0 \oplus \Omega^1 \oplus \Omega^1 \oplus  \Omega^2 \cdots $$
be a differential graded algebra, where we assume that $\Omega^0 =A$ and
the differential $d: \Omega^i \rightarrow \Omega^{i+1}$ increases the
degree. $d$ is a graded derivation in the sense that
$$ d(\omega_1 \omega_2)= d(\omega_1)\omega_2 +(-1)^{deg (\omega_1)}\omega_1
d(\omega_2), \quad \text{and} \quad d^2=0.$$ A {\it closed graded
trace} of dimension $n$ on $\Omega$ is a linear map
$$\int : \Omega^n \longrightarrow k$$
such that
$$ \int d\omega =0, \quad \text{and} \int (\omega_1
\omega_2-(-1)^{deg (\omega_1) deg (\omega_2)} \omega_2 \omega_1)=0,$$
for all $\omega$ in $\Omega^{n-1}$, $\omega_1$ in $\Omega^i$,  $\omega_2$ in $\Omega^j$ and $i+j=n.$
A triple of the form $(A, \Omega, \int)$ is called a {\it cycle} over
the algebra $A$.

Given a closed graded trace $\int$ on $A$, one defines a cyclic $n$-cocycle on $A$ by
$$\varphi (a_0, a_1, \cdots ,a_n)=\int a_0da_1 \cdots da_n.$$
Exercise: Check that $\varphi$ is a cyclic $n$-cocycle.

Conversely, one can show that   any cyclic cocycle on $A$ is
obtained in this way. To do this we introduce the   algebra
$(\Omega A, d)$ of {\it noncommutative differential forms} on $A$
as follows. $\Omega A$  is the universal (nonunital) differential
graded algebra generated by $A$ as a subalgebra. We have $\Omega^0
A=A$, and $\Omega^n A$ is linearly generated over $k$ by
expressions $a_0 da_1 \cdots da_n$ and $ da_1 \cdots da_n$ for
$a_i \in A$ (cf. \cite{ac85} for details). The differential $d$ is
defined by
$$d(a_0 da_1 \cdots da_n)=da_0 da_1 \cdots da_n, \quad \text{and} \quad
d ( da_1 \cdots da_n)=0.$$

Now it is easily checked that the relation
$$\varphi (a_0, a_1, \cdots ,a_n)=\int a_0da_1 \cdots da_n,$$
defines a 1-1 correspondence
$$\{\text{cyclic n-cocycles on $A$}\} \simeq \{\text{ closed graded traces on $\Omega^n A$}\}.$$

\noindent Exercise: Give a similar description for Hochschild cocycles $\varphi
\in C^n(A, A^*).$\\

\noindent 3. (From group cocycles to cyclic cocycles)\\
Let G be a discrete group
and let $c(g_1, \cdots, g_n)$ be a  group $n$-cocycle on $G$. Assume $c$ is
{\it normalized} in the sense that
$$c(g_1, \cdots , g_n)=0,$$
if $g_i=e$ for some $i$, or if $g_1g_2\cdots g_n=e$. One checks that
\begin{eqnarray*}
\varphi_c(g_0, \cdots , g_n)&=&c(g_1, \cdots ,g_n) \quad if \quad
g_0g_1\cdots g_n=e,\\
&=& 0 \quad \mbox{otherwise},
\end{eqnarray*}
is a cyclic n-cocycle on the group algebra $kG$ \cite{ac83, acb}. In
this way one obtains a map
$$H^n(G, k) \longrightarrow HC^n(kG), \quad c\mapsto \varphi_c.$$

By a theorem of D. Burghelea \cite{bur} (see below), the cyclic cohomology group $HC^n
(kG)$ decomposes over the conjugacy classes of $G$ and the component
corresponding to the conjugacy class of the identity  is exactly the
group cohomology $H^n(G, k)$. \\

\noindent 4. (From Lie algebra cocycles to cyclic cocycles)\\
Let $\mathfrak{g}$  be a Lie algebra acting by derivations on an
algebra $A$. This means we have a Lie algebra map
$$\mathfrak{g}\rightarrow Der(A, A),$$
from $\mathfrak{g}$ into the Lie algebra of derivations on $A$.
Let $\tau : A\rightarrow k$ be an {\it invariant trace} on $A$. Thus
$\tau$ is a trace on $A$ and
$$\tau (X(a))=0 \quad \text{ for all $X \in \mathfrak{g}, \, a \in A$}.$$

  For each $n\geq 0$, define  a map
$$\wedge^n g \longrightarrow C^n(A), \quad c\mapsto \varphi_c$$
$$\varphi_c(a_0, a_1, \cdots a_n)=\sum_{\sigma \in S_n}sgn (\sigma)\tau (a_0 X_{\sigma (1)}(a_1)\cdots
X_{\sigma (n)}(a_n)),$$
where $c=X_1 \wedge \cdots \wedge X_n.$\\
Exercise: Check that 1) For any $c$, $\varphi_c$ is a Hochschild cocycle, i.e. $b\varphi_c =0$. 2) If
$c$ is a Lie algebra cycle (i.e. if $\delta (c)=0$), the  $\varphi_c$ is a
cyclic cocycle.

We therefore obtain, for each $n\geq 0$,  a map
$$\chi_{\tau}: H_n^{Lie}(\mathfrak{g}, k) \longrightarrow HC^n (A), \quad c\mapsto \varphi_c,$$
from the Lie algebra homology of $\mathfrak{g}$ with trivial coefficients to the
cyclic cohomology of $A$ \cite{ac83}.

In particular if $\mathfrak{g}$ is abelian then $H^{Lie}_n(\mathfrak{g})\simeq
\wedge^n(\mathfrak{g})$ and we obtain a
well defined map
$$\wedge^n(\mathfrak{g}) \longrightarrow HC^n(A), \quad n=0, 1, \cdots.$$

Here is an example of this construction, first appeared in \cite{ac80}.
 Let $A=\mathcal{A}_{\theta}$ denote the ``algebra of smooth functions" on  the
 noncommutative torus. Let $X_1 =(1, 0), X_2=(0, 1)$. There is an action of
the abelian Lie algebra $\mathbb{R}^2$ on  $\mathcal{A}_{\theta}$
defined on generators of $\mathcal{A}_{\theta}$ by
\begin{eqnarray*}
X_1(U)=U, &\quad & X_1(V)=0,\\
X_2(U)=0, & \quad & X_2(V)=V.
\end{eqnarray*}
The induced derivations on $\mathcal{A}_{\theta}$ are given by
\begin{eqnarray*}
X_1(\sum a_{m, n}U^m V^n)&=& \sum ma_{m, n}U^m V^n,\\
X_2(\sum a_{m, n}U^m V^n)&=& \sum na_{m, n}U^m V^n.
\end{eqnarray*}
It is easily checked that the trace $\tau$ on  $\mathcal{A}_{\theta}$ defined by
$$\tau (\sum a_{m, n}U^mV^n)=a_{0,0},$$
is invariant under the above action of $\mathbb{R}^2$. The generators of
$H^{Lie}_{\bullet}(\mathbb{R}^2, \mathbb{R})$ are: $1, X_1, X_2, X_1\wedge X_2$.

We therefore obtain the following 0-dimensional,
1-dimensional and 2-dimensional cyclic cocycles on
$\mathcal{A}_{\theta}$:
\begin{eqnarray*}
\varphi_0 (a_0)&=& \tau (a_0),\\
\varphi_1(a_0, a_1)&=& \tau(a_0 X_1 (a_1)), \quad \varphi'_1 (a_0,
a_1)=\tau (a_0 X_2(a_1)),\\
\varphi_2(a_0, a_1,
a_2)&=&\tau (a_0( X_1(a_1)X_2(a_2)-X_2(a_1)X_1(a_2))).
\end{eqnarray*}

It is shown by Connes \cite{ac85} that these classes generate
the continuous periodic cyclic cohomology of $\mathcal{A}_{\theta}$.\\

\noindent 5. (Cup product and the $S$-operation on cyclic cohomology)\\
Let $(A, \Omega, \int)$ be an $m$-dimensional cycle on an algebra $A$ and $(B, \Omega', \int')$
and $n$-dimensional cycle on an algebra $B$. Let $\Omega \otimes
\Omega'$ denote the (graded) tensor product of the differential graded algebras $\Omega$ and $\Omega'$.
By definition, we have
$$(\Omega \otimes \Omega')_n =\bigoplus_{i+j=n}\Omega_i \otimes
\Omega'_j,$$
$$d(\omega \otimes \omega')=(d\omega) \otimes \omega' + (-1)^{deg
(\omega)}\omega \otimes (d\omega'),$$
$$\int'' \omega \otimes \omega'= \int \omega \int' \omega', \quad
\text{if} \quad \deg (\omega)=m, \; \deg (\omega')=n.$$
It is easily checked that $\int'' $ is a closed graded trace of dimension
$m+n$ on $\Omega \otimes \Omega'$.

Using the universal property of noncommutative differential forms,
applied to the identity map $A\otimes B \longrightarrow  \Omega_0
\otimes \Omega_0'$, one obtains a morphism of differential graded
algebras
$$ (\Omega (A\otimes B), d) \longrightarrow (\Omega \otimes \Omega',
d).$$
We therefore obtain a closed graded trace of dimension $m+n$  on
$(\Omega(A\otimes B), d).$ In \cite{ac85}, it is shown that the resulting map, called
{\it cup product} in cyclic cohomology,
$$\#:  HC^m(A) \otimes HC^n(B) \to HC^{m+n}(A\otimes B)$$
is well defined.\\

\noindent Exercise: By following the steps in the definition of the cup product, find an ``explicit''
 formula for $\varphi \# \psi$, when $\varphi$ is an $m$-dimensional cyclic cocycle and $\psi$ is an
 $n$-dimensional cyclic cocycle \cite{ac85}. 

Let $\beta$ denote the generator of $HC^2 (k)$ defined by $\beta (1, 1,
1)=1$. Using the cup product with $\beta$ we obtain the $S$-map
$$ S: HC^n(A) \to HC^{n+2}(A), \quad \quad \varphi \mapsto \varphi \# \beta.$$

\noindent Exercise: Find explicit formulas for $S \varphi$ when $\varphi$ is a $0$ or $1$ dimensional
cyclic cocycle. \\

In the next section we give a different approach to $S$ via the cyclic
bicomplex.

The {\it generalized trace map}
$$Tr : C_n^{\lambda} (A) \to C_n^{\lambda} (M_p (A)), $$
defined by 
$$ (Tr \varphi )(a_0\otimes m_0,  \cdots, a_n\otimes m_n)= tr (m_0\cdots m_n) \varphi (a_0, \cdots, a_n),$$
can be shown to be an example of cup product as well \cite{ac85}.
Indeed we have
$$ Tr (\varphi)= \varphi \# tr.$$

So far we studied the cyclic cohomology of algebras. There is a
``dual'' theory called {\it cyclic homology} that we introduce now.
Let $A$ be an algebra and for $n\geq 0$ let $C_n(A)= A^{\otimes
(n+1)}$. For each $n\geq 0$, define the operators
 $b: C_n(A) \longrightarrow C_{n-1}(A)$, $b': C_n(A) \longrightarrow C_{n-1}(A),$ and
 $ \lambda: C_n(A) \longrightarrow C_{n}(A)$  by
\begin{eqnarray*}
b(a_0\otimes \cdots \otimes a_n) & =& \sum_{i=0}^{n-1} (-1)^i( a_0\otimes
\cdots a_i a_{i+1} \otimes a_n)\\
& +& (-1)^n(a_na_0\otimes a_1\cdots \otimes a_{n-1}),\\
b'(a_0\otimes \cdots \otimes a_n) & =& \sum_{i=0}^{n-1} (-1)^i( a_0\otimes
\cdots a_i a_{i+1} \otimes a_n),\\
\lambda (a_0\otimes \cdots \otimes a_n) & =& (-1)^n( a_n\otimes
a_{0}
\cdots  \otimes a_{n-1}).
\end{eqnarray*}

The relation
$$ (1-\lambda )b'=b(1-\lambda)$$
can be easily established. Let $(C_{\bullet}(A), b)$ denote the
Hochschild complex of $A$ with coefficients in the $A$-bimodule
$A$ and let
$$ C_n^{\lambda}(A):= C_n(A)/Im (1-\lambda).$$
The relation
$ (1-\lambda )b'=b(1-\lambda)$ shows that the operator $b$ is well-defined on
$C_{\bullet}^{\lambda}(A)$. The quotient complex
$$(C_{\bullet}^{\lambda}(A), b)$$
is called
{\it cyclic complex} of $A$. Its homology, denoted by $HC_n (A), n=0, 1, \cdots$,  is
called the {\it cyclic homology} of $A$.\\

\noindent {\bf Example} For $n=0$, $HC_0(A)=A/[A, A]$ is the {\it
commutator quotient} of $A$. Here $[A, A]$ denotes the subspace of $A$
generated by the commutators $ab-ba$, for $a$ and $b$ in $A$.

\subsubsection{Connes' long exact sequence}
Our goal in this section is to establish  the long exact sequence
of Connes relating Hochschild and cyclic homology groups. There is
a similar sequence relating Hochschild and cyclic cohomology
groups. Connes' original proof in \cite{ac85} is based on the
notion of {\it cobordism of cycles}. This leads to an operator $B:
HH^n(A) \longrightarrow HC^{n-1}(A).$ He then shows that the three
operators $I: HC^n(A) \longrightarrow HH^n(A)$, $S: HC^n (A)
\longrightarrow HC^{n+2}(A)$ and $B$ fit into a long exact
sequence. An alternative approach is based on the following
bicomplex, called the {\it cyclic bicomplex} of $A$ and denoted by
$\mathcal{C}(A)$ \cite{ac83, lq, lod}:

$$\begin{CD}
\vdots @.\vdots @.\vdots @.\\
A^{\otimes 3}@<1-\lambda<< A^{\otimes 3} @<N<< A^{\otimes 3}@<1-\lambda<< \dots  \\
@VV bV @VV-b'V @VV bV  \\
A^{\otimes 2}@<1-\lambda<< A^{\otimes 2} @<N<< A^{\otimes 2}@<1-\lambda<< \dots \\
@VV bV @VV-b'V @VV bV  \\
A@<1-\lambda<< A @<N<< A@<1-\lambda<< \dots
\end{CD}
$$

Here the operator $N: A^{\otimes (n+1)} \longrightarrow A^{\otimes (n+1)}$ is defined by
$$N=1+\lambda +\lambda^2 +\cdots +\lambda^n.$$
The relations
\begin{eqnarray*}
N(1-\lambda)=(1-\lambda)N=0, & \quad &Nb=b'N, \\
(1-\lambda )b'=b(1-\lambda),& &
\end{eqnarray*}
can be easily verified. Coupled with
relations $b^2=0,\,  b'^2=0, $ it follows that  $\mathcal{C}(A)$ is a
bicomplex.

In particular we can consider its total complex $(Tot \mathcal{C}(A),
 \delta).$ We show that there is  a quasi-isomorphism of complexes
$$(Tot \mathcal{C}(A), \delta) \overset{ q. i.}{\sim} (C^{\lambda}(A), b).$$

To this end, for $n\geq 0$, define a map
\begin{eqnarray*}
(Tot_n \mathcal{C}(A), \delta) &\longrightarrow & (C^{\lambda}_{n}(A), b)\\
(x_0, x_1, \cdots ,x_n)&\mapsto & [x_n],
\end{eqnarray*}
where $[x_n]$ denotes the class of $x_n \in A^{\otimes (n+1)} $ in
$C^{\lambda}_n(A)=A^{\otimes (n+1)}/Im (1-\lambda).$
One checks that this is a morphism of complexes.
Now assume that $k$ is a field of characteristic zero. We show that  the rows of 
$\mathcal{C}(A)$, i.e. the complexes
$$A^{\otimes n+1}\overset{ 1-\lambda }{\longleftarrow} A^{\otimes n+1}
\overset{ N }{\longleftarrow} A^{\otimes n+1} \cdots $$
 are exact for each $n\geq 0$. The relation
$$Ker (1-\lambda)=Im N$$
is obvious. To show that $Ker N= Im (1-\lambda)$, define the
operator
$$ H=\frac{1}{n+1} (1+2 \lambda +3 \lambda^2  +\cdots + (n+1)\lambda^n): A^{\otimes n+1}
\longrightarrow A^{\otimes n+1}.$$
We have
$$ (1-\lambda) H= N+1,$$
which shows that $Ker N= Im (1-\lambda).$ It follows that the spectral
sequence converging to the total homology of $\mathcal{C}(A)$ collapses and the total homology
is the homology of its $E^2$-term. But the $E^2$-term is exactly the cyclic complex
$(C^{\lambda}_{\bullet}(A), b)$ of $A$. Alternatively, one can simply 
 apply a ``tic-tac-toe" argument to finish the proof.

We can utilize the 2-periodicity of the cyclic  bicomplex
$\mathcal{C}(A)$ to define a short exact sequence of complexes
$$ 0\longrightarrow Tot' \mathcal{C}(A) \longrightarrow Tot
\mathcal{C}(A) \longrightarrow Tot
\mathcal{C}(A) \longrightarrow 0,$$
where $Tot' \mathcal{C}(A)$ denotes the total complex of the first two columns and $Tot
\mathcal{C}(A)[2]$ is the shifted by two 2 complex. The last map is
defined by  truncation:
$$ (x_0, x_1, \cdots ,x_n) \mapsto (x_0, \cdots, x_{n-2}).$$
The kernel of this map is the total complex of the first two columns of
$\mathcal{C}(A)$:
$$Tot'_n \mathcal{C}(A)= A^{\otimes n}\oplus A^{\otimes (n-1)}.$$

Now when $A$ is unital  the $b'$-complex is acyclic. To prove this we define an operator
$s: A^{\otimes n} \longrightarrow A^{\otimes n+1}$,
$$ s(a_0\otimes \cdots \otimes a_{n-1})=1\otimes a_0\otimes \cdots
\otimes a_{n-1}.$$ One checks that
$$ b's+sb'=id,$$
which shows that the $b'$-complex is acyclic. We conclude that the
complex
  $Tot'_{\bullet}
\mathcal{C}(A)$ is homotopy equivalent to the Hochschild complex
$(C_{\bullet}(A), b)$.
Therefore   the  long exact sequence associated to the above   short
exact sequence  is of the form:
$$ \cdots \longrightarrow  HC_n(A) \overset{S}{\longrightarrow} HC_{n-2}(A)
\overset{B}{\longrightarrow}
 HH_{n-1}(A) \overset{I}{\longrightarrow}
HC_{n-1}(A) \longrightarrow \cdots.$$
This exact sequence was  first obtained by Connes in 1981 \cite{ac81,
ac85} in its cohomological form:
$$ \cdots \longleftarrow  HC^n(A) \overset{S}{\longleftarrow} HC^{n-2}(A)
\overset{B}{\longleftarrow}
 HH^{n-1}(A) \overset{I}{\longleftarrow}
HC^{n-1}(A) \longleftarrow \cdots.$$

The {\it periodicity operator}
$$S: HC_n(A) \rightarrow HC_{n-2}(A)$$ is induced by the
truncation map
$$ (x_0, x_1, \cdots ,x_n) \mapsto (x_0, \cdots, x_{n-2}).$$
 The operator $B$ is the {\it connecting homomorphism}
of the long exact sequence. It can therefore be represented on the level of
chains by the formula
$$B=Ns(1-\lambda): C_n(A) \rightarrow C_{n+1}(A).$$

\noindent Exercise: Recall how the connecting homomorphism is
defined for the long exact sequence associated to a short exact
sequence of complexes and drive the above formula for $B$.

\begin{remark}
 We defined, using the cup product with the generator
of $HC^2 (k)$, an operation $S: HC^n(A) \longrightarrow HC^{n+2}(A)$, for
$n=0,1, \cdots $. It can be shown that this definition coincides with the cohomological
form of the above  definition of $S$.
\end{remark}

Typical applications of Connes' long exact sequence involve extracting
information on cyclic homology from Hochschild homology. We list some of them:\\

\noindent 1. Let $f:A \rightarrow A'$ be an algebra homomorphism and suppose that
the induced maps on Hochschild groups
$$f_*: HH_n(A) \longrightarrow
HH_n(A'),$$
 are isomorphisms for all $n\geq 0$. Then
 $$f_*: HC_n(A)
\longrightarrow HC_n(A')$$
 is an isomorphism for all $n\geq 0.$ This simply follows by comparing the $SBI$ sequences for $A$ and $B$ and
applying the ``five lemma". In particular, it follows that inner
automorphisms act as identity on cyclic homology,  while inner
derivations act like zero on cyclic homology.\\

\noindent 2. (Morita invariance of cyclic homology). Let $A$ and $B$ be Morita
equivalent unital algebras. The Morita invariance property of
cyclic homology states that there is a natural isomorphism
$$ HC_n(A)\simeq HC_n (B), \quad n=0, 1, \cdots$$
For a proof of this fact in  general  see \cite{lod}. In the
special case where
 $B=M_k(A)$ a simple proof can be given as follows.  Indeed, by Morita invariance property
 of Hochschild homology, we know that
 the inclusion $i: A \rightarrow M_k(A)$ induces
isomorphism on Hochschild groups and therefore on cyclic groups by 1) above. \\

The bicomplex $\mathcal{C}(A)$ can be extended to the left. We obtain a
bicomplex, in the upper half plane,

$$\begin{CD}
\vdots @.\vdots @.\vdots @.\\
\dots @<N<< A^{\otimes 3}@<1-\lambda<< A^{\otimes 3} @<N<< A^{\otimes 3}@<1-\lambda<< \dots  \\
@VV bV @VV-b'V @VV bV  @VV -b'V\\
\dots @<N<<A^{\otimes 2}@<1-\lambda<< A^{\otimes 2} @<N<< A^{\otimes 2}@<1-\lambda<< \dots \\
@VV bV @VV-b'V @VV bV  @VV-b'V\\
\dots @<N<< A@<1-\lambda<< A @<N<< A@<1-\lambda<< \dots \\
\end{CD}
$$\\

The total homology (where direct products instead of direct sums is used in the
definition of the total complex) of this bicomplex is by definition the {\it periodic
cyclic homology} of $A$ and is denoted by $HP_{\bullet}(A)$. Note that because of the
2-periodicity of this bicomplex $HP_{\bullet}(A)$ is 2-periodic.\\

\noindent Exercise: Show that the resulting homology is trivial if
instead of direct products we use direct sums. \\

Similarly one defines the {\it periodic cyclic cohomology }
$HP^{\bullet}(A)$, where this time one uses direct sums for  the definition of the total complex.
 Since cohomology and direct limits commute, these periodic groups are indeed direct limits
of cyclic cohomology
groups under the $S$-map:
$$ HP^n(A)=\varinjlim HC^{n+2k}(A).$$

One checks, using the relations
$$b's+sb'= 1,\, bN=Nb', \, (1-\lambda)b=b'(1-\lambda)$$ that
$$B^2=0, \quad bB+Bb=0.$$
Indeed, we have
\begin{eqnarray*}
B^2=Ns(1-\lambda) Ns(1-\lambda)&=&0,\\
bB+Bb=b Ns(1-\lambda) +Ns(1-\lambda)b&=&\\
Nb's(1-\lambda)+Nsb'(1-\lambda)&=&\\
N(1)(1-\lambda)&=&0.
\end{eqnarray*}

The relations $b^2=B^2=bB+Bb=0$ suggest the following new bicomplex for
unital algebras, called {\it Connes' (b, B)-bicomplex}. It was
first defined in \cite{ac81, ac85}. As we shall see, it leads to a third
definition of cyclic (co)homology. Let $A$ be a unital algebra. The $(b,
B)$-bicomplex of $A$ is the following bicomplex:

$$\begin{CD}
\vdots @.\vdots @.\vdots \\
A^{\otimes 3}@<B<< A^{\otimes 2}@<B<< A \\
@VVbV @VVbV\\
A^{\otimes 2}@<B<<A\\
@VVbV\\
A
\end{CD}
$$
\begin{lem} The complexes $Tot \mathcal{B}(A)$ and $Tot \mathcal{C}(A)$
are homotopy equivalent.
\end{lem}
\begin{proof}
We define explicit chain maps between these complexes and show that they
are chain homotopic via explicit homotopies.
Define
\begin{eqnarray*}
I: Tot \mathcal{B}(A) &\rightarrow & Tot \mathcal{C}(A), \quad I=id +sN\\
J: Tot \mathcal{C}(A) &\rightarrow & Tot \mathcal{B}(A), \quad  J=id
+Ns.
\end{eqnarray*}
One checks that $I$ and $J$ are  chain maps.

Now consider the operators
\begin{eqnarray*}
g: Tot \mathcal{B}(A) &\rightarrow & Tot \mathcal{B}(A),  \quad g=B_0s^2N\\
h: Tot \mathcal{C}(A) & \rightarrow & Tot \mathcal{C}(A), \quad h=s,
\end{eqnarray*}
where $B_0=(1-\lambda)s.$

We have, by direct computation:
\begin{eqnarray*}
I\circ J&=& id + h \delta +\delta h,\\
I\circ J&=& id + g \delta' +\delta' g,\\
\end{eqnarray*}
where $\delta$ (resp. $\delta'$) denotes the differential of $Tot
\mathcal{C}(A)$ (resp. $Tot \mathcal{B}(A)$). (cf. \cite{kha2} for
details). The operator $I$ appears in \cite{lq}.  The operator $h$ and the second
homotopy formula was first
defined in \cite{kha2}.
\end{proof}

It is clear that the above result extends to  periodic
cyclic (co)homology.\\

\noindent {\bf Examples}\\
1. (Algebra of smooth functions). Let $M$ be a closed smooth manifold,
$A=C^{\infty}(M)$  denote the algebra of smooth complex valued functions
on $M$, and let
$(\Omega^{\bullet} M, d)$ denote the
de Rham complex of $M$. We saw that, by  a theorem of Connes, the map
$$\mu : C_n(A) \rightarrow  \Omega^n M, \quad \mu (f_0 \otimes \cdots \otimes f_n)= \frac{1}{n!}f_0df_1 \cdots df_n,$$
induces an isomorphism between the continuous Hochschild homology of $A$
and differential forms on $M$:
$$HH_n(A)\simeq \Omega^n M.$$

To compute the continuous  cyclic homology of $A$, we first show that under the map $\mu$ the
operator $B$ corresponds to the de Rham differential $d$. More
precisely, for each integer $n\geq 0$ we have a commutative diagram:
$$\begin{CD}
C_{n}(A)@> \mu >> \Omega^{n} M \\
@VVBV @VVdV\\
C_{n+1}(A)@>\mu >> \Omega^{n+1} M\\
\end{CD}
$$\\

We have
\begin{eqnarray*}
\mu B (f_0\otimes \cdots \otimes f_n)&=&\mu \sum_{i=0}^n (-1)^{ni}(1\otimes f_i \otimes
\cdots \otimes f_{i-1}-(-1)^n f_i \otimes
\cdots f_{i-1} \otimes 1)\\
&= & \frac{1}{(n+1)!}\sum_{i=0}^n (-1)^{ni}df_i \cdots df_{i-1}\\
&= &\frac{1}{(n+1)!}(n+1)df_0 \cdots df_n\\
&=& d \mu (f_0\otimes \cdots \otimes f_n).
\end{eqnarray*}

It follows that $\mu$ defines a morphism of bicomplexes
$$ \mathcal{B}(A) \longrightarrow \Omega (A),$$
Where $\Omega (A)$ is the bicomplex
$$\begin{CD}
\vdots @.\vdots @.\vdots \\
\Omega^2M@<d<< \Omega^1M@<d<< \Omega^0M \\
@VV0V @VV0V\\
\Omega^1M@<d<<\Omega^0M\\
@VV0V\\
\Omega^0M
\end{CD}
$$

 Since $\mu$ induces isomorphisms on row homologies, it induces isomorphisms on
 total homologies as well. Thus we have \cite{ac85}:
 $$HC_n(A)\simeq \Omega^n M/Im d \oplus H^{n-2}_{dR}(M) \oplus \cdots \oplus H^k_{dR}
 (M),$$
 where k=0 if $n$ is even and $k=1$ if $n$ is odd.

Using the same map $\mu$ acting between the corresponding periodic
complexes, one concludes that the periodic cyclic homology of $A$
is given by
$$HP_k (A)\simeq \bigoplus_i H^{2i+k}_{dR}(M), \quad k=0, 1.$$

2. (Group algebras). Let $kG$ denote the group algebra of a discrete
group $G$ over a field $k$ of characteristic zero. By a theorem of Burghelea \cite{bur},
 Hochschild
and cyclic homology groups of $kG$ decompose over the set of conjugacy
classes of $G$ where each summand is the group homology (with trivial coefficients) of a group
associated  to a conjugacy class. We recall this result.

Let $\widehat{G}$ denote the set of
conjugacy classes of $G$, $G'$ be the set of conjugacy classes
of elements of finite order, and let $G''$ denote the set of
conjugacy classes of elements of infinite order. For an element $g \in
G$, let $C_g=\{h\in G; \; hg=gh\}$ denote the centralizer of $g$, and let $W_g= C_g/<g>$, where $<g>$
is the group generated by $g$. Note that these groups depend, up to isomorphism,
 only on the conjugacy class of $g$. We denote
the group homology with trivial coefficients of a group $K$ by $H_*(K)$.

The Hochschild homology of $kG$ is given by \cite{bur, acb, lod}:
$$HH_n(kG)\simeq \bigoplus_{g\in \widehat{G}}H_n(C_g).$$

There is a similar, but more complicated, decomposition for the cyclic homology of $kG$:
$$HC_n(kG)\simeq \bigoplus_{g\in G'}(\bigoplus_{i\geq
0}H_{n-2i}(W_g))   \bigoplus_{g\in G''}H_n(W_g).$$
 In particular, the Hochschild group has $H_n(G)$ as a direct summand, while the cyclic
 homology group has $\oplus_i H_{n-2i}(G)$ as a direct summand  (corresponding to the conjugacy
 class  of the identity element of
 $G$).

\section{  Chern-Connes character}

Recall that the classical {\it Chern character} is a natural transformation
from $K$-theory to ordinary cohomology theory with rational coefficients \cite{ms}. More precisely
for each compact Hausdorff
space $X$ we have a natural homomorphism
$$Ch: K^0(X)\longrightarrow \bigoplus_{i\geq 0} H^{2i}(X, \mathbb{Q}),$$
where $K^0$ (resp. $H$) denote the $K$-theory (resp. Cech cohomology
with rational coefficients). It satisfies
 certain axioms and  these axioms completely
characterize $Ch$. But we won't recall these axioms  here since they are not very useful for  finding the noncommutative
analogue of $Ch$.  It is furthermore known that $Ch$ is a rational isomorphism
in the sense that upon tensoring it with $\mathbb{Q}$  we obtain an isomorphism
$$Ch_{\mathbb{Q}}: K^0(X) \otimes \mathbb{Q}
\overset{\sim}{\longrightarrow}
\bigoplus_{i\geq 0} H^{2i}(X, \mathbb{Q}).$$

When $X$ is a smooth manifold there is an alternative construction of
$Ch$, called the {\it Chern-Weil construction}, that uses  the differential
geometric notions of connection and curvature on vector bundles
\cite{ms}. It goes as follows. Let $E$ be a complex vector bundle on $X$
and let $\nabla$ be a connection on $E$. Thus
$$\nabla : C^{\infty}(E) \longrightarrow C^{\infty}(E)\otimes_A \Omega^1 X $$
is a $\mathbb{C}$-linear map satisfying the Leibniz condition
$$\nabla (fs)=f \nabla (s) +s\otimes df,$$
for all smooth sections $s$ of $E$  and smooth functions $f$ on $X$.
  Let
$$\hat{\nabla} : C^{\infty}(E)\otimes_A \Omega^{\bullet}X\longrightarrow C^{\infty}(E)\otimes_A
\Omega^{\bullet +1} X, $$
denote the natural extension of $\nabla$ satisfying a graded Leibniz
property. It can be easily shown that the {\it curvature operator}
$\hat{\nabla}^2$ is an $\Omega^{\bullet}X$-linear map. Thus it is
completely determined by its restriction to $C^{\infty}( E)$. This gives us  the {\it curvature form}
$$R \in C^{\infty}(End (E))\otimes \Omega^2 X$$
 of $\nabla$. Let
$$Tr: C^{\infty}( End (E))\otimes_A \Omega^{ev} X\to \Omega^{ev} X,$$
denote the canonical trace.    The Chern character of $E$ is then 
defined to be the class of the non-homogeneous even form
$$Ch(E)= Tr (e^{R}).$$
(We have omitted the normalization factor of $\frac{1}{2\pi i}$ to be multiplied by $R$.) One shows that $Ch(E)$ is a closed form and its cohomology class is
independent of the choice of  connection.

In \cite{ac80, ac85, acb},
Connes shows that
this Chern-Weil theory admits a vast generalization. For example, for an
algebra $A$
and  each integer $n\geq 0$ there are natural maps, called {\it
Chern-Connes character} maps,
$$Ch^{2n}_0 : K_0(A) \longrightarrow HC_{2n}(A),$$
$$Ch^{2n+1}_1 : K_1(A) \longrightarrow HC_{2n+1}(A),$$
compatible with $S$-operation.

Alternatively the noncommutative Chern character can be defined as
a pairing between cyclic cohomology groups and $K$-theory called
{\it Chern-Connes pairing}:
$$ HC^{2n}(A) \otimes K_0(A) \longrightarrow k,$$
$$HC^{2n+1}(A) \otimes K_1(A) \longrightarrow k.$$
These pairings are shown to be compatible with the periodicity
operator $S$ in the sense that
$$<[\varphi], [e]>=<S[\varphi], [e]>,$$
and thus induce a pairing between periodic cyclic cohomology and
$K$-theory.

We start by recalling the definition of these pairings. Let $\varphi
(a_0,\cdots ,a_{2n})$ be an {\it even cyclic cocycle} on an algebra $A$. For each integer $k\geq 1$, let
$$\tilde{\varphi}=tr\# \varphi \in C_{\lambda}^{2n}(M_k(A))$$
denote
the extension of $\varphi$ to the algebra of $k\times k$ matrices on
$A$. Note that $\tilde{\varphi}$ is a cyclic cocycle as well and is given by the formula
$$\tilde{\varphi} (m_0\otimes a_0, \cdots, m_{2n}\otimes a_{2n})=tr(m_0\cdots
m_{2n})\varphi (a_0, \cdots, a_{2n}).$$
 Let $e \in
M_k(A)$ be an idempotent. Define a bilinear map by the formula
$$<[\varphi], [e]>= \tilde{\varphi} (e, \cdots,e).$$
Let us first check that the value of the pairing depends only on the cyclic cohomology class of
$\varphi$ in $HC^{2n}(A)$. Suffices to assume $k=1$ (why?). Let
$\varphi =b\psi$ with $\psi \in C_{\lambda}^{2n-1}(A)$, be a
coboundary. Then we have
\begin{eqnarray*}
\varphi (e, \cdots,e)&=&b\psi (e, \cdots, e)\\
&=&\psi (ee, e,
\cdots,e)
-\psi (e, ee, \cdots, e)+ \cdots +(-1)^{2n} \psi(ee, e,
\cdots, e)\\
&=&\psi (e, \cdots,e)\\
&=&0,
\end{eqnarray*}
where the last relation follows from the cyclic property of $\psi$. 

To verify that the value of $<[\varphi], [e]>$, for fixed
$\varphi$, only depends on the class of $[e] \in K_0(A)$ we have
to check that for $u\in GL_k(A)$ an invertible matrix, we have
$<[\varphi], [e]>=<[\varphi], [ueu^{-1}]>$. Again  suffices to
show this for $k=1$. But this is exactly the fact, proved in the last section, that inner
automorphisms act by identity on  cyclic cohomology.

The formulas in the {\it odd case} are as follows. Given an invertible
matrix $u\in M_k(A)$ and an odd cyclic cocycle $\varphi (a_0, \cdots
,a_{2n+1})$ on $A$, we have
$$<[\varphi], [u]>= \tilde{\varphi}(u^{-1}-1, u-1, \cdots, u^{-1}-1, u-1).$$

\noindent Exercise: Show that the above formula defines a pairing
$K_1(A) \otimes HC^{2n+1}(A) \to k.$

There are also formulas for Chern-Connes pairings when the cyclic cocycle is in the  
 $(b, B)$ or cyclic bicomplex; but  we won't recall them here (cf. \cite{acb}, \cite{lod}). 

There is an alternative \textquoteleft \textquoteleft
infinitesimal proof" of the well-definement of these pairings which
works for Banach (or certain classes of topological) algebras
where elements of $K_0(A)$ can be defined as smooth homotopy
classes of idempotents \cite{acb}:

\begin{lem} Let $e_t, 0\leq t\leq 1,$ be a smooth family of idempotents in a
Banach  algebra $A$. There exists an smooth family $x_t, 0\leq t\leq 1$ of elements of $A$  such that
$$\overset{.}{e_t}:= \frac{d}{dt} (e_t)=[x_t, e_t], \quad \text{for} \quad 0\leq t\leq 1.$$
\end{lem}
\begin{proof} Let
$$x_t=[\overset{.}{e_t}, e_t]=\overset{.}{e_t}e_t-e_t\overset{.}{e_t}.$$
Differentiating the idempotent condition $e_t^2=e_t$ with respect to $t$
we obtain
$$\frac{d}{dt}(e_t^2)=\overset{.}{e_t}e_t +e_t\overset{.}{e_t}=\overset{.}{e_t}.$$
Multiplying this last relation on the left by $e_t$ we obtain
$$e_t\overset{.}{e_t}e_t=0.$$
Now we have
$$[x_t, e_t]=[\overset{.}{e_t}e_t-e_t\overset{.}{e_t}, e_t]=
\overset{.}{e_t}e_t +e_t\overset{.}{e_t}=\overset{.}{e_t}.$$
\end{proof}

It follows that if $\tau: A \rightarrow \mathbb{C}$ is a trace (= a cyclic zero
cocycle), then
$$ \frac{d}{dt} <\tau , e_t>=  \frac{d}{dt} \tau (e_t)=\tau (\overset{.}{e_t})=
\tau ([x_t, e_t])=0.$$
So  that the value of the pairing, for  a fixed $\tau$, depends only on the homotopy class of
the idempotent. This shows that the pairing
$$\{\mbox{ traces on $A$}\} \times K_0 (A) \longrightarrow \mathbb{C}$$
is well-defined.

 This is generalized in

\begin{lem} Let $\varphi (a_0, \cdots ,a_{2n})$ be a cyclic $2n$-cocycle on
$A$ and let $e_t$ be a smooth family of idempotents in $A$. Then the number
$$<[\varphi], [e_t]>= \varphi (e_t, \cdots , e_t)$$
is constant in $t$.
\end{lem}
\begin{proof}
Differentiating with respect to $t$ and using the above Lemma, we obtain
\begin{eqnarray*}
\frac{d}{dt}\varphi (e_t, \cdots, e_t)&=&\varphi (\overset{.}{e_t},
\cdots, e_t) +\varphi (e_t, \overset{.}{e_t},
\cdots, e_t)\\
& & \cdots +\varphi (e_t,
\cdots, e_t, \overset{.}{e_t})\\
&=& \sum_{i=0}^{2n} \varphi (e_t, \cdots, [x_t, e_t],
\cdots, e_t)\\
&=&L_{x_t} \varphi (e_t, \cdots, e_t).
\end{eqnarray*}
We saw that inner derivations act trivially on Hochschild and cyclic cohomology. This means that
for each $t$ there is a cyclic cocchain $\psi_t$ such that the {\it Lie
derivative} $L_{x_t} \varphi =b\psi_t$. We then have 
$$\frac{d}{dt}\varphi (e_t, \cdots, e_t)= (b\psi_t)( e_t, \cdots,
e_t)=0.$$
\end{proof}

Exercise: Repeat the above proof in the odd case.

The formulas for the even and odd {\it Chern-Connes character} maps
$$Ch^{2n}_0 : K_0(A) \longrightarrow HC_{2n}(A),$$
$$Ch^{2n+1}_1 : K_1(A) \longrightarrow HC_{2n+1}(A),$$
are as follows. In the even case, given an idempotent $e=(e_{ij})\in M_k(A)$, we have
$$ Ch^{2n}_0 (e)= Tr (\underbrace{e\otimes e\cdots \otimes e}_{2n+1})=\sum_{i_0,i_1, \cdots
i_{2n}} e_{i_0
i_1}\otimes e_{i_1i_2}\otimes \cdots \otimes e_{i_{2n}i_{0}}.$$
In low dimensions we have
\begin{eqnarray*}
 Ch^0_0(e)&=&\sum_{i=1}^k e_{ii},\\
 Ch^2_0(e)&=&\sum_{i_0=1}^k \sum_{i_1=1}^k
 \sum_{i_2=1}^ke_{i_0i_1}\otimes
 e_{i_1i_2}\otimes e_{i_2i_0}.
\end{eqnarray*}

In the odd case, given an invertible matrix $u\in M_k(A)$, we have
$$Ch^{2n+1}_1 ([u])=Tr (\underbrace{(u^{-1}-1)\otimes (u-1)\otimes\cdots \otimes
(u^{-1}-1)\otimes (u-1)}_{2n+2}).$$\\

\noindent {\bf Examples:}\\

\noindent 1. For $n=0$,  $HC^0(A)$ is the space of traces on $A$.
Therefore the Chern-Connes pairing reduces to the map
$$\{\mbox{traces on $A$}\} \times K_0(A) \longrightarrow k,$$
$$ <\tau, [e]>=\sum_{i=1}^n \tau (e_{ii}),$$
where $e=[e_{ij}] \in M_n(A)$ is an idempotent. The induced function on $K_0 (A)$ is
called the {\it dimension function} and denoted by $dim_{\tau}$. Here is an slightly different approach
to this dimension function.

 Let $E$ be a finite
 projective right $A$-module. A trace $\tau$ on $A$ induces  a trace on the endomorphism algebra of $E$,
$$Tr: End_A(E) \longrightarrow k$$
as follows.
First assume that $E=A^n$ is a free module. Then $End_A (E)\simeq M_n(A)$ and our
trace map is defined by
$$Tr (a_{i, j})=\sum a_{ii}.$$
It is easy to check that the above map is a trace. In general, there is an $A$-module $F$ such that
 $E\oplus F\simeq A^n$ is a free module and $End_A(E)$ embeds in
 $M_n(A)$. One can check that the induced trace on $End_A(E)$ is
 independent of the choice of splitting.\\
 Exercise: Since $E$ is finite and projectice, we have $End_A(E)\simeq
 E\otimes_A E^*$. The induced trace is simply the canonical pairing
 between $E$ and $E^*$. 

\begin{definition} The {\it dimension function} associated to a trace
$\tau$ on $A$ is the additive map
$$dim_{\tau}: K_0(A) \longrightarrow k,$$
induced by the map
$$dim_{\tau}(E)=Tr (id_E),$$
for any finite projective $A$-module $E$.
\end{definition}

It is clear that if $E$ is a vector bundle on a connected topological
space $X$ and $\tau (f) =f(x_0)$, where $x_0 \in X$ is a fixed point,
then  $dim_{\tau} (E)$ is the rank of the vector bundle $E$ and is an
integer. One of the striking features of noncommutative geometry is the existence of noncommutative
vector bundles with non integral dimensions. A beautiful example of this
phenomenon is shown by {\it Rieffel's idempotent} $e\in
\mathcal{A}_{\theta}$ with $\tau (e)=\theta$, where $\tau$ is the canonical
trace on the noncommutative torus \cite{acb}.\\

\noindent 2. Let $A=C^{\infty}(S^1)$ denote the algebra of smooth
complex valued functions on the circle. One knows that $K_1(A)\simeq K^1
(S^1)\simeq \mathbb{Z}$ and  $u(z)=z$ is a generator of this group. Let
$$\varphi (f_0, f_1)=\int_{S^1} f_0df_1$$
denote the cyclic cocycle on $A$ representing the fundamental class
of $S^1$ in de Rham homology. We have
$$ <[\varphi], [u]>=  \varphi (u, u^{-1})=\int_{S^1} udu^{-1}=-2\pi i.$$
Alternatively    the Chern character
$$Ch_1^1 ([u]) =u\otimes u^{-1}\in HC_1(A)\simeq H^1_{dR} (S^1),$$
is the class of the differential form $\omega = z^{-1}dz$, representing the
fundamental class of $S^1$ in de Rham cohomology. \\

\noindent 3. Let $A=C^{\infty} (S^2)$ and let $e\in M_2(A)$ denote the
idempotent representing the Hopf line bundle on $S^2$:
$$
e= \frac{1}{2}\left(
\begin{matrix}
1 + x_3 & x_1 + ix_2 \\
x_1-ix_2 & 1-x_3
\end{matrix}
\right).
$$
 Let us check that under the map
 $$HC_2(A)\to \Omega^2 S^2, \quad a_0\otimes a_1 \otimes a_2 \mapsto
 a_0da_1 da_2,$$
 the Chern-Connes character of $e$ corresponds to the fundamental class
 of $S^2$. We have
 $$ Ch_0^2(e)=Tr (e\otimes e\otimes e) \mapsto Tr (edede)=$$
 $$ \frac{1}{8}Tr \left(
\begin{matrix}
1 + x_3 & x_1 + ix_2 \\
x_1-ix_2 & 1-x_3
\end{matrix}
\right) \left(
\begin{matrix}
dx_3 & dx_1 + idx_2 \\
dx_1-idx_2 & -dx_3
\end{matrix}
\right) \left(
\begin{matrix}
dx_3 & dx_1 + idx_2 \\
dx_1-idx_2 & -dx_3
\end{matrix}
\right).
$$
Performing the computation one obtains
$$ Ch_0^2(e)\mapsto \frac{-i}{2}(x_1dx_2 dx_3-x_2dx_1 dx_3 +x_3dx_1 dx_2).$$
One can then integrate this 2-form on the two sphere $S^2$. The result is $-2\pi i$. In particular
the class of $e$ in $K_0(A)$ is non-trivial, a fact which can not be proved using just $Ch_0^0(e)=Tr(e)$. \\

\noindent 4. For smooth commutative algebras, the noncommutative Chern character reduces to the classical
Chern character. We verify this only in the even case.

Let $M$ be a
smooth closed manifold and let $E$
be a complex vector bundle on $M$. Let $e\in C^{\infty}(M,
M_n(\mathbb{C})$ be an idempotent representing the vector bundle $E$.
One can check that the following formula defines a connection on $E$,
called the Levi-Civita or Grassmanian connection:
$$ \nabla (eV)=edV,$$
where $V: M \to \mathbb{C}^n$ is a smooth function and $eV$
represents an arbitrary smooth section of $E$. Computing the
curvature form we obtain
$$R (eV)=\hat{\nabla}^2 (eV)= ede dV= edede. eV,$$
which shows that the curvature form is the ``matrix valued 2-form"
$$R=edede.$$
>From $e^2=e$, one easily obtains $ede.e=0.$ This implies that
$$R^n=(edede)^n=e\underbrace{dede\cdots dede}_{2n}.$$
Under the canonical map
$$HC_{2n}(A) \to H^{2n}_{dR} (M), \quad a_0\otimes \cdots \otimes
a_{2n}\mapsto \frac{1}{(2n)!}a_0da_1 \cdots da_{2n},$$
we have
$$Ch_{n}^0(e):=Tr(e\otimes \cdots \otimes e)\mapsto \frac{1}{(2n)!}Tr (edede\cdots de).$$
The classical Chern-Weil formula for $Ch(E)$ is
$$ Ch (E)= Tr (e^R)=Tr (\sum_{n\geq 0} \frac{R^n}{n!}).$$
So that its  $n$ th component is given by
$$Tr \frac{R^n}{n!}= \frac{1}{n!}Tr ((edede)^n)= \frac{1}{n!}Tr (ede\cdots de). $$

\appendix
\section{Banach and $C^*$-algebras}

By an {\it algebra} in this book we mean an {\it associative algebra} over a commutative unital
ground ring $k$.
An algebra is called {\it unital} if there is a (necessarily unique) element $1 \in A$ such that
$1a = a1 = a$ for all $a\in A$. It is called {\it commutative}  if $ab=ba$ for all
$a,b \in A$.

Now let $k=\mathbb{R}$ or $\mathbb{C}$ be the field of real
or complex numbers. A {\it norm} on a real or complex algebra $A$ is a map
$$\| \; \| : A\to \mathbb{R},$$
such that for all $a, b$ in $A$ and $\lambda$ in $k$ we have:
\begin{enumerate}
\item[1)] $\| a \| \ge 0 , \;  \text{and}\; \| a \| = 0 \; \text{iff} \; a=0,$
\item[2)] $\| a+b \| \le \| a \| + \| b \|, $
\item[3)] $\| \lambda a \| = |\lambda | \| a \|,$
\item[4)] $\| ab \| \le \| a \|\| b \|.$
\end{enumerate}
If $A$ is unital, we assume that $\| 1 \| = 1$. An algebra endowed with a norm is called a
{\it normed algebra}.

A {\it Banach algebra} is a  normed algebra which is {\it complete}. Recall that a normed
vector space $A$ is called complete  if any Cauchy sequence in $A$ is
convergent.  One of the main consequences of completeness is that
absolutely convergent series are convergent, i.e. if $ \sum_{n=1}^\infty \| a_n \|$ is
convergent, then $\sum_{n=1}^\infty a_n$ is convergent in
$A$. In particular the geometric series $\sum_{n=1}^{\infty}a_n $ is
convergent if $\| a \|<1$.  From this it easily follows that the group
of invertible elements in a unital Banach algebra $A$ is an open subset of $A$.

An {\it $\ast$-algebra} is a complex algebra  endowed with an $\ast$-operation, i.e. a map
$$\ast: A\to A, \hspace{20pt} a\mapsto a^\ast,$$
which is {\it anti-linear} and {\it involutive}:
\begin{enumerate}
\item[1.] $(a+b)^\ast = a^\ast + b^\ast , \; (\lambda a)^\ast =\bar{\lambda} a^\ast, $
\item[2.] $(ab)^\ast = b^\ast a ^\ast,$
\item[3.] $(a^\ast )^\ast = a, $
\end{enumerate}
for all $a, b$ in $A$ and $\lambda $ in $\mathbb{C}$.

A {\it Banach $\ast$-algebra} is a complex Banach algebra endowed with an $\ast$-operation such that
for all $a\in A, \| a^\ast \| = \| a \|$. In particular for all $a$ in $A$ we have
$$\| a^\ast a \| \le \| a^\ast \|\| a \| = \| a \|^2. $$

\begin{definition}
A {\it $C^\ast$-algebra} is a Banach $\ast$-algebra $A$ such that for all $a
\in A$,
$$\| aa^\ast \| = \| a \|^2. $$
We refer to this last identity as the {\it $C^*$-identity}.
\end{definition}

For reasons that will become  clear  later in this section,
$C^\ast$-algebras occupy a very special place among all Banach
algebras. This should be compared with the role played by Hilbert
spaces among all Banach spaces. In fact, as we shall see
 there is an intimate relationship between Hilbert spaces and $C^*$-algebras thanks to the
 GNS construction and the Gelfand-Naimark embedding theorem.

A {\it morphism} of  $C^*$-algebras is an algebra homomorphism
$$ f: A \longrightarrow B $$
between $C^*$-algebras $A$ and $B$  such that $f$  preserves the $*$ structure, i.e.
  $$ f(a^*)=f(a)^*, \quad \text{for all} \; a\in A.$$

  It can be shown that morphisms of $C^*-$ algebras are {\it contractive}
  in the sense that for all $a \in A$,
  $$\| f(a)\| \le \| a\|.$$
  In particular they are automatically continuous. It follows from this
  fact that the norm of a $C^*$-algebra is unique in the sense that if $(A, \|\, \|_1)$ and
   $(A, \|\, \|_2)$ are both $C^*-$algebras then
$$\|a \|_1 =\|a \|_2,$$
for all $a\in A$. Note also that a morphism  $f$ of $C^*-$algebras
is an {\it isomorphism} if an only if $f$ is one to one and onto.
Isomorphisms of $C^*$-algebras are necessarily {\it isometric}.

In sharp distinction from $C^*$-algebras, one can have different
Banach algebra norms on the same algebra. For example on the
algebra of $n\times n$ matrices one can have different Banach
algebra norms and only one of
them is a $C^*$-norm.\\

\noindent {\bf Examples:}\\
\begin{enumerate}
\item[1.] (commutative $C^\ast$-algebras). Let $X$ be a locally compact Hausdorff space. We associate
to $X$ several classes of algebras of functions on $X$ which are $C^\ast$-algebras.

\item[1.a] Let
$$C_0 (X) = \{f:X\to \mathbb{C} ; f\mbox{ is continuous and $f$ vanishes at $\infty$} \} \; .$$
By definition, $f$ vanishes at $\infty$ if for all $ \epsilon >0$
there exists  a compact subset $K\subset X$ so that $|f(x) | <
\epsilon $ for all  $x \in X\setminus K$.  Under pointwise
addition and scalar multiplication $C_0(X)$ is obviously an
algebra over the field of complex numbers $\mathbb{C}$. Endowed
with the sup-norm
$$\| f \| = \| f \|_\infty = \sup \{ |f(x) | ; \; x\in X \}, $$
and $\ast$-operation
$$f\mapsto f^\ast,\; f^\ast (x) = \bar{f}(x),$$
one  checks that
$C_0(X)$ is a commutative $C^\ast$-algebra. It is unital if and only if $X$ is compact. If $X$ is compact,
we simply write $C(X)$ instead of $C_0(X)$.

By  a theorem of Gelfand and Naimark, any commutative $C^*$-algebra is of the type $C_0(X)$ for some locally
compact Hausdorff space $X$ (see below).
\item[1.b] Let
$$C_b (X) = \{ f:X\to \mathbb{C} ; \; f \mbox{ is continuous and bounded} \}. $$
Then with the same operations as above, $C_b (X)$ is a unital $C^\ast$-algebra.  Note
that $C_0 (X) \subset
C_b (X)$  is an {\it essential ideal} in $C_b (X)$. (An ideal $I$ in an algebra
$A$ is an  essential ideal if for all $a$ in $A$, $aI = 0 \Rightarrow a=0.) $

\item[2.] (commutative Banach algebras). It is easy to give examples of Banach algebras which are not
$C^\ast$-algebras. For an integer $ n\geq 1$, let
$$C^{n} [0,1] = \{ f: [0,1] \to \mathbb{C} ; \; f \in C^n \}, $$
be the space of functions with  continuous $n$-th derivative.  We denote the $i$-th derivative
of $f$ by $f^{(i)}$. With the norm
$$\| f \|_n = \sum_{i=0}^n \frac{\| f^{(i)} \|_\infty}{i!}$$
and the $\ast$-operation $f^\ast (x) = \bar{f}(x)$, one checks that $C^{n}[0, 1]$ is a
Banach $\ast$-algebra.   It
is, however,  not a $C^\ast$-algebra as one can easily show that the $C^\ast$-identity  fails.
Note that for all $f \in C^{n} [0,1]$,
$$\| f \|_\infty \le \| f \|_n \; .$$

\item[3.] (noncommutative $C^\ast$-algebras). By a theorem of Gelfand and Naimark recalled below,
any
 $C^\ast$-algebra can be realized as a closed $\ast$-subalgebra of
the algebra of bounded operators on a complex Hilbert space. We start
with the simplest examples: the algebra of complex $n$ by $n$ matrices.

\item[3.a] Let $A=M_n(\mathbb{C})$ be the algebra of $n$ by $n$ matrices over the field of complex
numbers $\mathbb{C}$. With operator norm and the standard adjoint
operation
$T\mapsto T^*$, $A$ is a $C^*$-algebra (see below).

A direct sum of matrix algebras
$$A=M_{n_1}(\mathbb{C})\oplus M_{n_2}(\mathbb{C})\oplus \cdots
\oplus M_{n_k}(\mathbb{C})$$
is a $C^*$-algebra as well. It can be shown that any finite dimensional
$C^*$-algebra is unital and is a direct sum of matrix algebras as above \cite{D}. In other words,
finite dimensional $C^*$-algebras are semi-simple.

\item[3.b] The above example can be generalized as follows.
Let $H$ be a complex Hilbert space and let $ A= \mathcal{L}(H)$ denote the set of bounded linear operators
$H\to H$. For a bounded  operator  $T$,  we define $T^\ast$ to be  the {\it adjoint} of $T$ defined by
$$<Tx, y > = <x, T^\ast y>,  \hspace{20pt} x,y \in H.$$

Under the usual algebraic operations
of addition and
multiplication of operators and the {\it operator norm}
$$\| T \| = \sup \{ \| T(x) \|; \; \| x \| \le 1 \}, $$
$\mathcal{L}(H)$ is a $C^\ast$-algebra.

It is clear that any subalgebra $A\subset \mathcal{L}(H)$ which is {\it self-adjoint} in the sense that
$$T\in A \Rightarrow T^\ast \in A,$$
 and is {\it norm closed}, in the sense that
 $$T_n \in A,\; \|T_n- T\| \rightarrow 0 \Rightarrow T\in A,$$
  is a $C^\ast$-algebra.
\item[3.c] (group $C^{\ast}-$algebras). Let $G$ be  a discrete group and
let $H=l^2 G$ denote the Hilbert space of square summable functions on
$G$;
$$H=\{f:G\to \mathbb{C}; \sum_{g\in G} |f(g)|^2<\infty \}.$$
The {\it left regular representation} of $G$ is the unitary
representation $\pi$ of $G$ on $H$, defined by
$$(\pi g)f (h)=f(g^{-1}h).$$
It has a linear extension to an (injective) algebra homomorphism
$$\pi :\mathbb{C}G \longrightarrow \mathcal{L}(H),$$
from the group algebra of $G$ to the algebra of bounded operators on
$H$. Its image $\pi (\mathbb{C}G)$ is a $\ast$-subalgebra of
$\mathcal{L}(H)$.

The {\it reduced group $C^*$-algebra} of $G$, denoted by $C^*_rG$,
is the norm closure of $\pi (\mathbb{C}G)$ in $\mathcal{L}(H)$. It
is obviously a $C^*$-algebra.

There is second $C^*$-algebra associated to any discrete group $G$ as
follows. The (non-reduced) {\it group $C^*$-algebra} of $G$ is the norm completion
of the $\ast$-algebra $\mathbb{C}G$ under the norm
$$\|f\|=\mbox{sup}\; \{\|\pi (f)\|; \pi \, \mbox{ is a $\ast$-representation of $\mathbb{C}G$}\},$$
where by a $\ast$-representation we mean a  $\ast$-representation on a
Hilbert space. Note that $\|f\|$ is finite since for $f=\sum_{g \in
G}a_g g$ (finite sum) and any *-representation $\pi$ we have
$$\| \pi (f)\|\leq \sum \| \pi (a_gg)\| \leq  \sum |a_g| \|\pi (g)\| \leq \sum |a_g|.$$

By its very definition it is clear that there is a 1-1
correspondence between unitary representations of $G$ and $C^*$
representations of $C^* G$.

Since the identity map $id: (\mathbb{C}G, \| \,\|) \to (\mathbb{C}G, \| \,\|_r)$ is continuous,
we obtain  a surjective $C^*-$algebra homomorphism
$$C^*G \longrightarrow C^*_rG.$$
It is known that this map is an isomorphism if and only if $G$ is
an amenable group \cite{D }. Abelian groups are amenable.

We give a few examples of reduced group $C^*$-algebras. Let $G$ be an
abelian group and $\hat{G}=Hom (G, \mathbb{T})$ the group of characters
of $G$. It is a locally compact Hausdorff space. Moreover
it is easily seen that $\hat{G}$ is in fact homeomorphic
with the space of characters, or the maximal ideal space, of the $C^*$-algebra  $C^*_rG$.
Thus the Gelfand transform defines an isomorphism  of $C^*$-algebras
$$ C^*G \simeq C_0(\hat{G}).$$

In general one should think of the  group $C^*$-algebra of a group $G$  as the ``algebra of
functions" on the noncommutative
space representing the unitary dual of $G$. Note that, by the above paragraph, this is justified
  in the commutative case.
 In the noncommutative case, the unitary dual is a badly
behaved space in general but the noncommutative dual is a perfectly legitimate noncommutative space (see
 the unitary dual of the infinite dihedral group in \cite{acb} and its noncommutative replacement).

Let $G$ be a finite group. Since $G$ is finite the group $C^*$-algebra coincides with the group
algebra of $G$. From basic representation theory we know that
the group algebra $\mathbb{C}G$ decomposes as a sum of matrix algebras
$$ C^*G\simeq \mathbb{C}G\simeq \oplus M_{n_i}(\mathbb{C}),$$
where the summation is over the set of conjugacy classes of $G$.

\end{enumerate}

 It is generally believed that the classic paper of Gelfand and
Naimark \cite{gelnai} is the birth place of the theory of
$C^*$-algebras. The following two
results  on the structure of $C^*$-algebras are proved in this paper:\\

\begin{theorem}[Gelfand-Naimark \cite{gelnai}] a) Let $A$ be a commutative
$C^*$-algebra and let $\Omega (A)$ denote the maximal ideal space of $A$.  Then the Gelfand transform
$$ A \rightarrow C_0(\Omega (A)), \quad a\mapsto \hat{a},$$
is an isomorphism of $C^*$-algebras. \\
b) Any $C^\ast$-algebra is  isomorphic to a $C^\ast$-subalgebra of the algebra $\mathcal{L}(H)$
of bounded linear operators on some Hilbert
space $H$.
\end{theorem}

In the remainder of this appendix we sketch the proofs of statements a)
and b) above. They are based on Gelfand's theory of commutative Banach
algebras, and the Gelfand-Naimark-Segal (GNS) construction of
representations of $C^*-$algebras from states, respectively.

\subsection{Gelfand's theory of commutative Banach algebras}
The whole theory is based on the notion of spectrum of an element of a Banach algebra and the fact
that the
spectrum   is non-empty. The notion of
spectrum can be defined for elements of an arbitrary algebra and it can
be easily shown that for finitely generated complex algebras the spectrum is non-empty. As is shown in
\cite{cg}, this latter fact leads to an easy proof of Hilbert's Nullstellensatz. This makes the proofs
of the two major duality theorems remarkably similar. We use this
approach in this book.

Let $A$ be a unital  algebra over a field $\mathbb{F}$. The {\it spectrum}
of an element $a\in A$ is defined by
$$\mbox{sp}(a) = \{ \lambda \in \mathbb{F} ; \; a-\lambda 1 \; \mbox{is {\it not}
invertible} \} \; .$$

We should think of the spectrum as the noncommutative analogue of the
set of values of a function. This is justified in Example 1
below.\\

\noindent {\bf Examples:} \\
\begin{enumerate}
\item[1.] Let $A=C(X)$ be the algebra of continuous complex  valued functions
on a compact space $X$. For any $f\in A$,
$$\mbox{sp} (f) = \{ f(x) ;\; x\in X \},$$
is the range of $f$.
\item[2.] Let $A=M_n (\mathbb{F} )$ be the algebra of $n\times n$ matrices with
coefficients in $\mathbb{F}$. For any matrix $a\in A$
$$\mbox{sp}(a) = \{ \lambda \in \mathbb{F} ; \; \det (a-\lambda 1 ) = 0 \} \; , $$
is the set of eigenvalues of $a$.
\end{enumerate}

\noindent Exercise: \\
1) Show that if $a$ is nilpotent then sp $(a)=\{0\}$.\\
\noindent 2) Show that
$$\mbox{sp} (ab)\setminus \{0\}=\mbox{sp} (ba)\setminus \{0\}.$$
 \noindent 3) Let $T: H \to H$ be a Fredholm operator on a Hilbert space
 $H$. Let $Q$ be an operator such that $1-PQ$ and $1-QP$ are trace class
 operators. Show that
 $$Index (T)= Tr(1-PQ)-Tr(1-QP).$$\\

In general, the spectrum may be empty.  We give two general results that guarantee the spectrum is
non-empty. They are at the
foundation of Gelfand-Naimark theorem and Hilbert's Nullstellensatz. Part b) is in
\cite{cg}.

\begin{theorem}
\begin{enumerate}
\item[(a)] (Gelfand) Let $A$ be a unital Banach algebra over $\mathbb{C}$. Then for any
$a\in A,\; \mbox{sp}(a) \ne \emptyset . $
\item[(b)] Let $A$ be a unital algebra over $\mathbb{C}$. Assume $\dim_{\mathbb{C}} A$ is
countable. Then for any $a\in A,\; \mbox{sp} (a) \ne \emptyset .$
Furthermore, an element $a$ is nilpotent if and only if sp $(a)=\{0\}.$
\end{enumerate}
\end{theorem}
\begin{proof} We sketch a proof of both statements. For a) assume the
spectrum of an element $a$ is empty. Then the function
$$R: \mathbb{C}
\to A, \quad \lambda \mapsto (a-\lambda 1)^{-1},$$ is holomorphic
(in an extended sense), non-constant, and bounded. This is easily
shown to contradict the classical   Liouville's theorem from
complex analysis.

For b), again assume the spectrum of a is empty. Then it can be
shown that the uncountable set
$$\{(a-\lambda 1)^{-1}; \; \lambda \in \mathbb{C} \}$$
is a linearly independent set. But this contradicts the fact that $dim_{\mathbb{C}}A$ is
countable.

For the second part of b), assume sp $(a)=\{0\}.$ Since
$dim_{\mathbb{C}}(A)$ is countable, any element $a\in A$ satisfies a polynomial equation. Let
$$p(a)=a^k(a-\lambda_1)\cdots (a-\lambda_n)=0$$
be the minimal polynomial of $A$. Then $n=0$ since otherwise an
element $a-\lambda_i$ is not invertible with $\lambda_i \neq 0$.
But this contradicts our assumption that sp $(a)=\{0\}.$  The
other direction is true in general and is easy.
\end{proof}

The first part of the following corollary is known as Gelfand-Mazur
theorem.

\begin{cor}
Let $A$ be either a unital complex Banach algebra or a unital   complex algebra
such that  $\dim_{\mathbb{C}} A$ is
countable. If $A$ is a division algebra, then $A\simeq\mathbb{C}.$
\end{cor}

Let $A$ be an algebra. By a {\it character} of $A$ we mean  a non-zero algebra
homomorphism
$$\varphi :A\to
\mathbb{F}.$$
 Note that if $A$ is unital, then $\varphi (1)=1$. We establish the link
between characters and maximal ideals of $A$. For the following
result $A$ is either a commutative unital complex Banach algebra,
or is a commutative unital algebra with $\dim_\mathbb{C} A$
countable.

\begin{cor}
The relation $I=\ker\varphi$ defines a 1-1 correspondence between the set of maximal
ideals of $A$ and the set of characters of $A$.
\end{cor}

Before embarking on the proof of Gelfand-Naimark theorem, we sketch a
proof of Hilbert's Nullstellensatz, following \cite{cg}.

Let
$$ A=\mathbb{C} [x_1, \cdots ,x_n]/I$$
be a finitely generated commutative reduced algebra. Recall that reduced
means if $a^n=0$ for some $n$ then $a=0$ (no nilpotent elements). Equivalently the ideal $I$ is radical.
Let
$$ V=\{ z \in \mathbb{C}^n ; \; p(z)=0, \mbox{ for all $p$ in $I$}\},$$
let $J(V)$ be the ideal of functions vanishing
on $V$,  and  let
$$\mathbb{C}[V]=\mathbb{C}[x_1, \cdots  ,x_n]/J(V)$$
be  the algebra of regular functions on $V$. Since $I\subset
J(V)$, we have an algebra homomorphism
$$\pi: A \longrightarrow
\mathbb{C}[V].$$ One of the original forms of  Hilbert's
Nullstellensatz states that this map is an isomorphism. It is clearly
surjective. For its injectivity, let $a\in A$ and let $\pi (a)=0$, or
equivalently $\pi (a)\in J(V)$. Since $a$ vanishes on all points of $V$,
it follows that $a$ is in the intersection of all the maximal ideals of $A$. This shows that
its spectrum sp ($A$)=$\{0\}$. By Theorem A.2 (b), it follows that $a$
is nilpotent and since $A$ is reduced, we have $a=0$.  

The rest of this section is devoted
to sketch a proof of the Gelfand-Naimark theorem on the structure of
commutative $C^*$-algebras. . Let $A$ be a unital Banach algebra. It is easy to see that any character of $A$ is
continuous of norm 1. To prove this, note that if this is not the case then there
exists an  $a\in A$, with $\| a \| < 1$ and $\varphi (a) = 1$. Let $ b = \sum_{n\ge 1}
a^n$. Then from $a+ab = b$, we have
$$\varphi (b) = \varphi (a) + \varphi (a) \varphi (b) = 1 + \varphi (b),  $$
which is impossible. Therefore $\| \varphi \| \le 1$, and since $\varphi (1) = 1, \;
\| \varphi \| =1$.

Let $A$ be a complex Banach algebra and let $\Omega (A)$ denote
the set of characters of $A$. Thus if $A$ is unital, then $\Omega
(A)$= set of maximal ideals of $A$. It is clear that a  pointwise
limit of characters is again a character. Thus $\Omega (A)$ is a
closed subset of the unit ball of the dual space $A^\ast$. Since
the latter space is a compact Hausdorff space in the
$\mbox{weak}^\ast$ topology, we conclude that  $\Omega (A)$ is
also a compact Hausdorff space.

If $A$ is not unital, let $A^+ = A\oplus \mathbb{C}$ be the unitization of $A$. It is
clear that
$\Omega (A) = \Omega (A^+) \setminus \{ \varphi_0 \},$
 where
$\varphi_0$ is the trivial character $\varphi_0 (a)=0$ for all  $a \in A$. Since
$\Omega (A^+ )$ is compact, we conclude that $\Omega (A)$ is a locally compact Hausdorff
space. We have thus proved the lemma:

\begin{lem}
Let $A$ be a Banach algebra. Then $\Omega (A)$ is a locally compact Hausdorff
space. $\Omega (A)$ is compact if and only if $A$ is unital.
\end{lem}

Let $f: A\to B$ be a continuous homomorphism of commutative unital Banach algebras. Define
a map $\Omega (f)=f^\ast : \Omega (B) \to \Omega (A)$ by
$$f^\ast (\varphi ) = \varphi \circ f.$$
It is clear that $f^\ast \varphi$ is multiplicative, linear and continuous.
Thus we have defined the {\it spectrum functor} $\Omega$
from the category of commutative unital complex Banach algebras to the category of  compact
Hausdorff spaces.

Next we define the {\it Gelfand transform}. Let $A$ be a commutative
Banach algebra. The Gelfand transform is the map $\Gamma$ defined by
$$ \Gamma:  A\to C_0 ( \Omega (A)),$$
$$a\mapsto \hat{a} \; , \; \hat{a} (\varphi ) = \varphi (a). $$
This map is obviously an algebra homomorphism. It is also clear that  $\| \Gamma \| \le 1$, i.e.
$\Gamma$ is contractive.

The spectrum of an element is easily seen to be a closed and
bounded subset of $\mathbb{C}$ (this follows from: $\|1-a\|<1
\Rightarrow a$ is invertible with inverse $a^{-1}=\sum_{n\geq 0}
(1-a)^n.$) Let
$$r(a)=\text{max}\; \{|\lambda |; \; \lambda \in \text{sp}(a)\},$$
denote the {\it spectral radius} of an element $a \in A.$ Note that  $r(a)
\leq \|a\|.$

Now Corollary 6.2 tells us that
$$\mbox{sp}(a)=\{\varphi (a); \: \varphi \in \Omega (A)\}.$$
Then following result is then immediate:

\begin{prop}
The Gelfand transform $A\to C_0 (\Omega (A))$ is a norm decreasing algebra homomorphism
and its image separates points of $\Omega (A)$. Moreover, for all $a\in A$,
$\| \hat{a} \|_{\infty} = r(a) \; .$
\end{prop}

The kernel of the Gelfand transform is called the {\it radical} of the Banach algebra
$A$. It consists of elements $a$ whose spectral radius $r(a)=0$, or  equivalently, sp$(a)=\{
0 \}$. Hence the radical contains all the nilpotent elements. But it may contain more.
An element $a$ is called
{\it quasi-nilpotent} if sp$(a)=0$. $A$ is said to be {\it semi-simple} if its
radical is zero, i.e. the only quasi-nilpotent elements of $A$ is 0.\\

\noindent{\bf Example.} We give an example of a commutative Banach
algebra for which the Gelfand transform is injective but not surjective.
 Let $H(D)$ be the space of continuous functions on the unit disk
$D$ which are holomorphic in the interior of the disk. With the  sup-norm $\| f\|=\| f\|_\infty$
it is a Banach algebra. It is, however, not  a
$C^\ast$-algebra (why?). Show that $\Omega (A)\simeq D$ and the Gelfand transform is an isometric embedding $H(D)\to C(D)$.

We are now ready to prove the first main theorem of Gelfand and
Naimark in \cite{gelnai}: for commutative $C^*$-algebras $\Gamma$
is an isometric $*$-isomorphism. We need a few simple facts about
$C^*$-algebras first.

Let $f: A \to B$ be a morphism of $C^*$-algebras. It is easily seen that
sp $(f(a))\subset $ sp ($a$). Hence, using the $C^*$-identity, we have
$$\|f(a)\|^2=\|f(a) f(a^*)\| =\|f(aa^*)\| =r(a^*a)\leq
\|a^*a\|=\|a\|^2.$$

Let $a\in A$ be a standpoint element ($a=a^*$).  Then sp ($a$)
$\subset \mathbb{R}$ is real. Indeed since $e^{ia}$ is unitary its
spectrum is located on the unit circle (why?). Hence for $\lambda
\in sp (a)$, $e^{i\lambda}$ is located on the unit circle which
shows that $\lambda$ is real. From this it follows that if $f:
A\to  \mathbb{C}$ is an algebra homomorphism (a character), then
$f (a^*)=\overline{f(a)},$ for all $a$.

\begin{theorem}
Let $A$ be a commutative $C^\ast$-algebra. The Gelfand transform $A\to C_0 (\Omega (A))$
is an isomorphism of $C^\ast$-algebras.
\end{theorem}

\begin{proof}
We prove the unital case. The non-unital case follows with minor modifications \cite{D}.
What we have shown so far amounts to the fact that $\Gamma$ is an
isometric $*$-algebra map whose image separates points of $\Omega (A)$.
Thus 
 $\Gamma (A)$ is a closed $*$-subalgebra of $C(\Omega (A))$ that
 separates points of $\Omega (A)$. By Stone-Weierstrass theorem, $\Gamma (A)=C(\Omega (A)).$ 

\end{proof}

The above theorem is one of the landmarks  of Gelfand's theory of commutative Banach
algebras.  While a complete classification of all commutative Banach algebras seems to be
impossible, this result classifies all commutative $C^\ast$-algebras.

Care must be applied in dealing with non-compact spaces and non-unital
algebras. First of all it is clear that if $X$ is not compact then the pull-back $f^*: C_0(Y)\to C_0(X), \quad
f^*(g)=g\circ f$
of a  continuous map $f:X \to Y$ is well defined if and only if $f$ is a
proper map. Secondly, one notes that not all $C^*$-maps $C_0(Y)\to C_0(X)$ are obtained in this way.
For
example, for $X=(0, 1)$ an open interval and $Y$  a single point,  the
zero map $0: C_0(Y) \to C_0(X)$,
which is always a $C^*$-morphism, is not the pull-back of any proper map. It
turns out that one way to single out the appropriate class of morphisms
is as follows. A morphism of $C^*$-algebras $\varphi: A \to B$ is
called {\it proper} if for an approximate unit $(e_i)$ of $A$, $\varphi
(e_i)$ is an approximate unit of $B$. Recall that \cite{D,gvf} an {\it
approximate unit} of a $C^*$-algebra $A$ is an increasing net of positive elements $(e_i), i\in I,$ of
 $A$ such that for all $a \in A$, $e_ia\to a$, or equivalently, $ae_i
 \to a$.  Now it can be shown that any proper map $C_0(Y) \to C_0(X)$ is the   pull-back of a proper
 map $X\to Y$.  Similarly, in the other direction, if $\varphi : A \to B$ is a proper morphism
 of $C^*$-algebras, then $\Omega (\varphi): \Omega (B) \to \Omega (A)$ is a proper continuous map.
 We refer to \cite{gvf} and references therein for  more details. 

We are now half-way through showing that the functors $C_0$ and $\Omega$ are quasi-inverse to
each other. But the proof of the other half is much simpler and is left
to the reader.

\begin{lem}
Let $X$ be a locally compact Hausdorff space. Then the {\it evaluation map}
$$X\to \Omega (C_0 (X)), \quad \; x\mapsto \varphi_x,$$
$$ \varphi_x (f)=f(x),$$
is a homeomorphism.
\end{lem}

\begin{example} We give a few elementary applications of the
Gelfand-Naimark correspondence between commutative $C^*$-algebras and
locally compact spaces.
\begin{enumerate}
\item[1.] (Idempotents and connectedness). Let $A$ be a unital
commutative $C^\ast$-algebra. Then
 $\Omega (A)$ is disconnected iff $A$ has a non-trivial idempotent (i.e. an element $e\ne 0, 1$
such that $e^2 = e$).
\item[2.] (Ideals and closed subsets). Let $X$ be a compact Hausdorff space. The Gelfand-Naimark
 duality
 shows that there is a 1-1 correspondence
$$\{ \mbox{closed subsets of } X \} \leftrightarrow \{ \mbox{closed ideals of } C(X) \}, $$
where to each closed subset $Y\subset X$, we associate the ideal
$$I= \{ f\in C(X) ; f(y) = 0 \quad \forall\; y\in Y \}$$
of all functions vanishing on $Y$. We have natural isomorphisms
$$C_0 (X\setminus Y) \simeq I,$$
$$C_0(X)/ I \simeq C_0(Y), $$
$$C_0(X/ Y)\simeq C_0 (X\setminus Y)^+.$$

\item[3.] (Essential ideals and compactification)
Let $X$ be a locally compact Hausdorff space. Recall that a {\it Hausdorff compactification} of $X$ is
a compact Hausdorff space $Y$ where $X$ is homeomorphic to a dense subset of $Y$. We consider $X$
as a subspace of $Y$. Then $X$ is open in $Y$ and its boundary $Y\setminus X$ is compact. We have
an exact sequence
$$0\to C_0(X) \to C(Y) \to C(Y\setminus X) \to 0,$$
where $C_0 (X)$ is an essential ideal of $C(Y)$. Conversely, show that any extension
$$0\to C_0(X) \to A \to B \to 0,$$
Where $A$ is a unital $C^*$-algebra and $C_0(X)$ is an essential
ideal of $A$, defines a Hausdorff compactification of $X$. Thus,
we have a 1-1 correspondence between Hausdorff compactifications
of $X$ and (isomorphism classes of) essential extensions of
$C_0(X)$. In particular,   the 1-point compactification and the
Stone-Cech compactifications correspond to
$$0\to C_0(X) \to C_0 (X)^+ \to \mathbb{C} \to 0,$$ and
$$0\to C_0(X) \to C_b (X) \to C(\beta X) \to 0$$
\end{enumerate}
\end{example}

\subsection{ States and the GNS construction}

Our goal in this section is to sketch a proof of the second main
result of Gelfand and Naimark in \cite{gelnai} to the effect that
any $C^*$-algebra can be embedded in the algebra of bounded
operators on a Hilbert space. The main idea of the proof is an
adaptation of the idea of {\it left regular representation} of
algebras to the context of $C^*$-algebras, called
Gelfand-Naimark-Segal (GNS) construction.

A {\it positive linear functional} on a $C^*$-algebra $A$ is a
$\mathbb{C}$-linear map $\varphi : A \longrightarrow \mathbb{C}$ such
that for all $a$ in $A$,
$$\varphi (a^* a)\geq 0.$$
A {\it state} on $A$ is a positive linear functional $\varphi$ with
$\|\varphi\|=1.$ It can be shown that if $A$ is unital then
this last condition is equivalent to $\varphi (1)=1.$

If $\varphi_1$ and $\varphi_2$ are states then  for any $t\in [0, 1]$,
t$\varphi_1$ +(1-t) $\varphi_2$ is a state as well. Thus the set of
states of $A$, denoted by $\mathcal{S}(A)$,
form a convex subset of the unit ball of $A^*$. The extreme points of
$\mathcal{S}(A)$ are called {\it pure states}. \\

\noindent {\bf Examples:}\\

\noindent 1. States are noncommutative analogues of probability measures. This
idea is corroborated by the Riesz representation theorem: For  a locally compact Hausdorff
space $X$ there is a 1-1
correspondence between states on $C_0(X)$ and Borel probability measures
on $X$. To a probability measure $\mu$ is associated the stats $\varphi$ defined by
$$ \varphi (f)=\int_X fd\mu. $$
 $\varphi$ is a pure  state if and only if $\mu =\delta_x$ is a Dirac
 measure for a point $x\in X$.

\noindent 2. Let $A=M_n(\mathbb{C})$ and $p\in A$ be a positive matrix with $tr(p)=1$. (Such matrices, and
their infinite dimensional analogues, are called {\it density matrices} in quantum statistical
 mechanics.) Then
$$\varphi (a)=tr (ap)$$
defines a state on $A$.

\noindent 3. Let $\pi : A \longrightarrow \mathcal{L}(H)$ be a {\it representation
} of a unital $C^*$-algebra $A$ on a Hilbert space $H$. This simply means that $\pi$ is a
morphism of unital $C^*$-algebras. Let $x\in H$ be a vector of length
one . Then
$$\varphi (a)=<\pi (a)x, x>$$
defines a state on $A$, called a {\it vector state}. In the following we show that, conversely, any
state on $A$ is a vector state with respect to a suitable representation called the GNS  representation.

Let $\varphi$ be a positive linear functional on $A$. Then
$$<a, b>:=\varphi (b^*a)$$
is linear in the first variable and anti-linear in the second variable.
It is also semi-definite in the sense that $<a, a>=\varphi (a^*a)\geq 0$
for all $a$ in $A$. Thus it satisfies the {\it Cauchy- Schwartz}
inequality: for all $a, b$
$$ |\varphi(b^*a)|^2\leq \varphi(a^*a)\varphi(b^*b).$$

Let
$$N=\{ a\in A; \quad f(a^*a)=0\}.$$
It is easy to see, using the above Cauchy-Schwarz inequality, that $N$
is a closed left ideal  of $A$ and the following positive definite inner product is
well-defined on the quotient space $A/N$:
$$<x+N, y+N>:= <x, y>.$$

Let $H_{\varphi}$ denote the Hilbert space completion of $A/N$
under the above inner product. The {\it left regular
representation} $A \times A \to A, (a, b)\mapsto ab$ of $A$ on
itself induces a bounded linear map $A \times A/N \to A/N, (a,
b+N)\mapsto ab+N$. We denote its unique extension to
$H_{\varphi}$ by
$$\pi_{\varphi} : A \longrightarrow \mathcal{L}(H_{\varphi}).$$

The representation $(\pi_{\varphi}, H_{\varphi})$  is called the
GNS representation defined by the state $\varphi$. The state
$\varphi $ can be recovered from the representation
$(\pi_{\varphi}, H_{\varphi})$ as a vector state as follows. Let
$x=\pi_{\varphi}(1)$. Then for all $a$ in $A$,
$$\varphi (a)=<(\pi_{\varphi}a)(x), x>.$$

The representation $(\pi_{\varphi}, H_{\varphi})$ may not be
faithful. It can be shown that it is irreducible if and only if
$\varphi$ is a pure state \cite{D}. To construct a faithful
representation, and hence an embedding of $A$ into the algebra of
bounded operators on  a Hilbert space, one first shows that there
are enough pure states on $A$. The proof of the following result
is based on Hahn-Banach and Krein-Milman theorems.

\begin{lem} For any selfadjoint element $a$ of $A$, there exists a pure
state $\varphi$ on $A$ such that $\varphi (a)=\|a\|.$
\end{lem}

Using the GNS representation associated to $\varphi$, we can then
construct, for any $a \in A$,  an irreducible representation $\pi$
of $A$ such that $\|\pi (a) \|=|\varphi (a)|=\|a\|.$

We can now prove the second theorem of Gelfand and Naimark.
\begin{theorem} Every $C^*$-algebra is isomorphic  to a
$C^*$-subalgebra of the algebra of bounded operators on a Hilbert space.
\end{theorem}
\begin{proof}
Let $\pi=\sum_{\varphi \in \mathcal{S}(A)} \pi_{\varphi}$ denote
the direct sum of all GNS representations for  all  states of $A$.
By the above remark $\pi$ is faithful.
\end{proof}

\section{Idempotents and finite projective modules}

Let $A$ be a unital algebra over a commutative  ring $k$ and let
$\mathcal{M}_A$ denote the category of right $A$-modules. We assume our
modules $M$
are {\it unitary} in the sense that the unit of the algebra acts as the
identity of $M$. A morphism  of this category is a
right $A$-module map  $f: M \rightarrow N$, i.e. $f(ma)=f(m)a$, for all $a$ in $A$ and $m$ in
$M$.

A {\it free} module, indexed by a set $I$, is a module of the type
$$M=A^I=\bigoplus_IA,$$
where the action of $A$ is by component-wise right multiplication.
Equivalently, $M$ is free if and only if there are elements $m_i
\in M, i\in I$, such that any $m\in M$ can be uniquely expressed
as a finite sum $m=\sum_i m_i a_i$. A module $M$ is called {\it
finite} (= {\it finitely generated}) if there are elements $m_1,
m_2, \cdots ,m_k$ in $M$  such that every element of $m\in M$ can
be expressed as $m=m_1a_1 +\cdots +m_ka_k$, for some $a_i\in A.$
Equivalently $M$ is finite if there is a surjective $A$-module map
$A^k \rightarrow M$ for some integer $k$.

Free modules correspond to trivial vector bundles. To obtain a
more interesting class of modules we consider the class of {\it
projective modules}. A module $P$ is called projective if it is a
direct summand of a free module. That is there exists
 a module $Q$ such that
$$P\oplus Q \simeq A^I.$$

A module is said to be {\it finite  projective} (= {\it finitely generated projective}),
if it is both finitely generated
and projective.

\begin{lem}
Let $P$ be an $A$-module. The following conditions on $P$ are
equivalent:
\begin{enumerate}
\item[1.] $P$ is projective.
\item[2.] Any surjection
$$M\overset{f}{\to} P \to 0,$$
splits in the category of $A$-modules.
\item[3.]For all $A$-modules $N$ and $M$ and morphisms $f,g$ with
$g$ surjective in the following diagram, there exists a morphism $\tilde{f}$ such that the diagram
commutes:
$$
\xy
\xymatrix@=9ex{
 & P \ar@{-->}[dl]_{\exists \tilde{f}} \ar[d]^f\\
N \ar[r]^g & M \ar[r] & 0
}
\endxy
$$
We say that $\tilde{f}$ is a lifting of $f$ along $g$.
\item[4.] The functor
$$\mbox{Hom}_A (P,-):\mathcal{M}_A\to \mathcal{M}_k$$
 is exact in the sense that
for any short exact sequence of $A$-modules
$$0 \to R\to S\to T \to 0,$$
the sequence of $k$-modules
$$0\to \mbox{Hom}_A (P,R)\to \mbox{Hom}_A (P,S) \to \mbox{Hom}_A (P,T) \to 0$$
is exact.
\end{enumerate}
\end{lem}

\begin{example}

\begin{enumerate}
\item[1.] Free modules are projective.
\item[2.] If $A$ is a division ring, then any $A$-module is free, hence projective.
\item[3.] $M=\mathbb{Z}/n\mathbb{Z}, \; n\ge 2$, is not projective as a $\mathbb{Z}$-module.
\item[4.] A direct sum $P=\oplus_i P_i$ of modules is projective iff each summand $P_i$
is projective.
\item[5.] (idempotents) Let
$$e\in M_n (A) = \mbox{End}_A (A^n), \quad e^2 = e,$$
be an idempotent.
Multiplication by $e$ defines a right $A$-module map
$$P_e : A^n \to A^n,\; \quad x\mapsto ex.$$
Let $P=\mbox{Im} (P_e)$ and $Q=\mbox{Ker} (P_e)$ be the image and kernel
of $P_e$.  Then,
using the idempotent condition $e^2=e$, we obtain a direct sum
decomposition
$$P\oplus Q = A^n,$$
which shows that  both  $P$ and $Q$ are projective modules.
Moreover they are obviously finitely generated. It follows that
both $P$ and $Q$ are finite projective  modules.

Conversely, given any finite projective  module $P$, let $Q$ be a module
such that $P\oplus Q \simeq A^n$ for some integer $n$. Let $e: A^n
\rightarrow A^n$ be the right $A$-module map that corresponds to the
projection map
$$ (p, q)\mapsto (p, 0).$$
Then it is easily seen that we have an isomorphism of $A$-modules
$$ P\simeq P_e.$$
This shows that any finite projective module is obtained from an
idempotent in some matrix algebra over $A$.
\end{enumerate}
\end{example}

 The idempotent $e \in M_n(A)$  associated to a finite projective $A$-module $P$ depends of course
 on the choice of the splitting $P\oplus Q\simeq A^n.$ Let $P\oplus
 Q'\simeq A^m$ be  another  splitting and $f\in M_m(A)$ the corresponding
 idempotent.  Define the operators $u\in Hom_A(A^m, A^n)$, $v\in Hom_A(A^n, A^m)$ as compositions
 $$ u: A^m \overset{\sim}{\longrightarrow} P\oplus Q\longrightarrow P\longrightarrow
 P\oplus Q' \overset{\sim}{\longrightarrow} A^n,$$
$$ v: A^n \overset{\sim}{\longrightarrow} P\oplus Q'\longrightarrow P\longrightarrow
 P\oplus Q \overset{\sim}{\longrightarrow} A^m.$$
 We have
 $$ uv=e, \quad vu=f.$$
 In general, two idempotents satisfying the above relations are called
 Murray-von Neumann equivalent. Conversely it is easily seen that Murray-von Neumann equivalent
 idempotents  define isomorphic finite projective modules.

Defining finite projective modules through idempotents is certainly very
convenient since both finiteness and projectivity of the module are automatic in this case.
In some cases however this is not very useful.  For example, modules over quantum tori are
directly defined and
checking directly that they are finite and projective is difficult. In this case
the following method, due to Rieffel \cite{rie2},  is very useful. (I am grateful  to Henrique
Bursztyn for bringing   this method to my attention).

Let A and B be unital algebras over a ground ring $k$.
let X be an (A,B)-bimodule endowed with $k$-bilinear  pairings (``algebra-valued" inner
products):
$$ <-,->_A : X \times X \longrightarrow A, $$
$$ <-,->_B : X \times X  \longrightarrow B,$$
satisfying the ``associativity" condition
$$ <x, y>_Az = x<y, z>_B  \quad \text{ for all x, y, z in $X$},$$
and the ``fullness" conditions:
$$ <X, X>_A=A \quad \text{ and}\quad  <X, X>_B=B.$$

We claim that $X$ is a finite projective left $A$-module.
Let $1_B$  be the unit of $B$. By fullness of $<-,->_B$,
we can find $x_i, y_i $ in $X, i=1, \cdots, k$ such that
$$1_B = \sum_{i=1}^k <x_i, y_i>_B.$$

Let $e_i, i=1, \cdots, k$, be a basis for $A^k$.
Define the map
$$P: A^k \longrightarrow X,  \quad  P(e_i)=y_i.$$
We claim that this map splits and hence $X$ is finite and projective.
Consider the map
$$I: X \longrightarrow A^k, \quad I(x)= \sum_i <x,x_i>_Ae_i. $$
We have
$$ P\circ I(x)= \sum_i <x,x_i>_Ay_i = \sum_i x<x_i,y_i>_B  \quad \text{(by
associativity)}$$
       $$= x \quad   \text{(since $\sum_i<x_i,y_i>_B=1_B$)}.$$

\section{ Equivalence of categories}

There are at least two ways to compare categories: isomorphism and  equivalence.
{\it Isomorphism} of categories is a very strong requirement and is hardly useful.
 {\it Equivalence} of categories, on the other hand, is a much more flexible
concept and is very useful.

Categories $A$ and $B$ are said to be {\it equivalent} if there is a
functor $F:A\to B$ and
 a functor $G:B\to A$, called a {\it quasi-inverse} of $F$, such that
$$F\circ G \simeq 1_B,\quad G\circ F\simeq 1_A \;, $$
where $\simeq$ means isomorphism, or natural equivalence, of functors.
This means,
for every $X\in \mbox{obj} A, \; Y\in \mbox{obj} B,$
$$FG(Y)\sim Y, \quad  \text{and} \quad  GF(X)\sim
X,$$
where $\sim$ denotes  natural isomorphism of objects.

If $F\circ G = 1_B$ and $G\circ F=1_A$ (equality
of functors), then we say that categories $A$ and $B$ are {\it
isomorphic}, and
$F$ is an isomorphism.

Categories $A$ and $B$ are said to be {\it antiequivalent} if the opposite
category $A^{op}$ is equivalent to $B$.

Note that a functor $F:A \to B$ is an isomorphism if and only if $F:\mbox{obj} A \to \mbox{obj} B$ is
1-1, onto and  $F$ is
{\it full and faithful} in the sense  that  for all $X,Y\in \mbox{obj}
A,$
$$F:\mbox{Hom}_A (X,Y) \to \mbox{Hom}_B (F(X),F(Y))$$
is 1-1 (faithful) and onto (full).

It is easy to see that an equivalence $F:A\to B$ is full and faithful, but it may not
be 1-1, or  onto on the class of objects. As a result an equivalence may have many
quasi-inverses. The following concept clarifies the situation with objects of equivalent categories.

A subcategory
$A'$ of a category $A$ is called {\it skeletal}
if 1) the embedding $A' \to A$  is full, i.e.
$$\mbox{Hom}_{A'} (X,Y) = \mbox{Hom}_A (X,Y)$$
for all
$X,Y \in \mbox{obj} A'$ and 2) for any object $X\in \mbox{obj} A$, there is a
unique object $X^\prime \in \mbox{obj} A^\prime$ isomorphic to $X$. Any skeleton of
$A$ is equivalent to $A$ and it is not difficult to see that two categories $A$ and $B$
are equivalent if and only if they have isomorphic skeletal subcategories $A^\prime$
and $B^\prime$.

In some examples, like the Gelfand-Naimark theorem, there is a canonical
choice for a
quasi-inverse for a given equivalence functor $F$ ($F=C_0$ and $G=\Omega$). There are instances,
however, like the Serre-Swan
theorem, where there is no canonical choice for a quasi-inverse. The following proposition gives a
necessary and sufficient condition for a functor $F$ to be an
equivalence of categories. We leave its simple proof to the reader.

\begin{prop}
A functor $F:A\to B$ is an equivalence of categories if and only if
\begin{enumerate}
\item[a)] $F$ is full and faithful, and
\item[b)] Any object $Y\in \mbox{obj} B$ is isomorphic to an object of the form $F(X)$, for some
$X \in obj A$.
\end{enumerate}
\end{prop}

\end{document}